\documentclass[oneside,american,english]{amsart}
\usepackage[T1]{fontenc}
\usepackage{textcomp}
\usepackage[latin9]{inputenc}
\usepackage{mathtools}
\usepackage{amsbsy}
\usepackage{amstext}
\usepackage{amsthm}
\usepackage{amssymb}
\usepackage{geometry}
\makeatletter
\numberwithin{equation}{section}
\numberwithin{figure}{section}
\theoremstyle{plain}
\newtheorem{thm}{\protect\theoremname}[section]
\theoremstyle{plain}
\newtheorem{conjecture}[thm]{\protect\conjecturename}
\theoremstyle{plain}
\newtheorem*{thm*}{\protect\theoremname}
\theoremstyle{definition}
\newtheorem{defn}[thm]{\protect\definitionname}
\theoremstyle{remark}
\newtheorem{rem}[thm]{\protect\remarkname}
\theoremstyle{definition}
\newtheorem{example}[thm]{\protect\examplename}
\theoremstyle{plain}
\newtheorem{lem}[thm]{\protect\lemmaname}
\theoremstyle{plain}
\newtheorem{prop}[thm]{\protect\propositionname}
\theoremstyle{plain}
\newtheorem{cor}[thm]{\protect\corollaryname}

\def\makebbb#1{
    \expandafter\gdef\csname#1\endcsname{
        \ensuremath{\Bbb{#1}}}
}\makebbb{R}\makebbb{N}\makebbb{Z}\makebbb{C}\makebbb{H}\makebbb{E}\makebbb{H}\makebbb{P}\makebbb{B}\makebbb{Q}\makebbb{E}\makebbb{E}
\makebbb{H}
\makebbb{F}
\makebbb{A}
\makebbb{G}
\makebbb{K}

\usepackage{babel}

\usepackage{babel}

\makeatother

\usepackage{babel}

\makeatother

\usepackage{babel}
\addto\captionsamerican{\renewcommand{\conjecturename}{Conjecture}}
\addto\captionsamerican{\renewcommand{\corollaryname}{Corollary}}
\addto\captionsamerican{\renewcommand{\definitionname}{Definition}}
\addto\captionsamerican{\renewcommand{\examplename}{Example}}
\addto\captionsamerican{\renewcommand{\lemmaname}{Lemma}}
\addto\captionsamerican{\renewcommand{\propositionname}{Proposition}}
\addto\captionsamerican{\renewcommand{\remarkname}{Remark}}
\addto\captionsamerican{\renewcommand{\theoremname}{Theorem}}
\addto\captionsenglish{\renewcommand{\conjecturename}{Conjecture}}
\addto\captionsenglish{\renewcommand{\corollaryname}{Corollary}}
\addto\captionsenglish{\renewcommand{\definitionname}{Definition}}
\addto\captionsenglish{\renewcommand{\examplename}{Example}}
\addto\captionsenglish{\renewcommand{\lemmaname}{Lemma}}
\addto\captionsenglish{\renewcommand{\propositionname}{Proposition}}
\addto\captionsenglish{\renewcommand{\remarkname}{Remark}}
\addto\captionsenglish{\renewcommand{\theoremname}{Theorem}}
\providecommand{\conjecturename}{Conjecture}
\providecommand{\corollaryname}{Corollary}
\providecommand{\definitionname}{Definition}
\providecommand{\examplename}{Example}
\providecommand{\lemmaname}{Lemma}
\providecommand{\propositionname}{Proposition}
\providecommand{\remarkname}{Remark}
\providecommand{\theoremname}{Theorem}

\begin{document}
\title{Gibbs polystability of Fano manifolds, stability thresholds and symmetry
breaking}
\begin{abstract}
We extend the probabilistic approach for constructing Kähler--Einstein
metrics on log Fano manifolds $(X,\Delta)$ -- involving random point
processes -- to the case of non-discrete automorphism groups, by
breaking the symmetry using a moment map constraint. In particular,
an algebraic notion of Gibbs polystability is introduced, ensuring
that the corresponding point processes on $X$ are well-defined. We
conjecture that the Gibbs polystability of $(X,\Delta)$ is equivalent
to the existence of a Kähler--Einstein metric and that the unique
such metric with vanishing moment emerges when sampling a large number
of $N$ points on $X.$ The definition of Gibbs polystability is formulated
in terms of a new algebro-geometric invariant, defined as a limit
of certain log canonical thresholds on the GIT semistable locus of
the $N$-fold products $X^{N}.$ We conjecture that it coincides with
a previously studied analytic polystability threshold, encoding the
coercivity of the K-energy functional modulo automorphisms. These
conjectures follow from an overarching conjectural Large Deviation
Principle for the large $N-$limit. The conjectured equivalence between
Gibbs polystability and K-polystability is established when $(X,\Delta)$
is a log Fano curve. For Fano manifolds of any dimension we establish
an effective one-sided version of the conjecture, which yields an
effective criterion for the existence of a Kähler-Einstein metric
and, more generally, an effective lower algebraic bound on the analytic
polystability threshold. In companion papers applications will be
given to optimal stability results for the logarithmic HLS inequality
on the two-sphere, to Onsager's point vortex model on the two-sphere
and the AdS/CFT correspondence. 
\end{abstract}

\author{Rolf Andreasson, Robert J. Berman, Ludvig Svensson}
\maketitle

\section{Introduction}

The quest for canonical metrics on complex manifolds is a classical
theme in geometry, tracing back to the uniformization theorem. According
to the seminal \emph{Yau--Tian--Donaldson conjecture} a compact
complex manifold $X$ admits a Kähler--Einstein metric with positive
Ricci curvature, i.e. a Kähler metric with constant positive Ricci
curvature, iff $X$ is a \emph{K-polystable }Fano manifold --- this
is now a theorem \cite{c-d-s,ti2}. The Fano condition requires that
the anti-canonical line bundle $-K_{X}$ of $X$ is ample, while K-polystability
is a finer algebro-geometric condition that amounts to a positivity
property of all polarized $\C^{*}$-equivariant deformations of $(X,-K_{X}),$
called test configurations. This algebro-geometric characterization
has also found striking connections to the AdS/CFT correspondence
in theoretical physics, where Kähler--Einstein metrics on Fano manifolds
induce Sasaki--Einstein metrics on links that are dual to superconformal
gauge theories \cite{m-s-y,C-S,c-z-y}. More precisely, Fano manifolds
yield regular links, while Fano orbifolds give rise to quasi-regular
links (and, upon introducing weights, Sasaki--Einstein metrics on
general links \cite{a-c,li0}). 

Even though the existence problem for Kähler--Einstein metrics is
now settled, it should be emphasized that the problem of explicitly
describing the Kähler--Einstein metric in terms of algebro-geometric
data of $X$ is, in general, intractable, and only very few explicit
examples are known. This problem is one of the motivations for the
probabilistic approach to the construction of Kähler--Einstein metrics
on $X$ introduced in \cite{berm8,berm8 comma 5}, that was further
applied to the AdS/CFT correspondence in \cite{b-c-p}. In the case
of positive Ricci curvature, it requires, however, that the Fano manifold
$X$ does not admit any holomorphic vector fields, i.e. that the group
$\text{Aut}_{0}(X)$ of automorphisms of $X$ homotopic to the identity
be trivial. In the present work the framework is extended to Fano
manifolds with \emph{nontrivial} $\text{Aut}_{0}(X)$. This is a first
paper in a series. In companion papers we will present applications
to optimal stability results for the logarithmic HLS inequality on
the two-sphere \cite{a-b-s3}, to Onsager's point vortex model on
the two-sphere \cite{a-b-s1} and the AdS/CFT correspondence \cite{a-b-s2}. 

\subsection{\label{subsec:Motivation}Recap of the case when $\text{Aut}_{0}(X)$
is trivial}

Before stating our main results, let us recall the probabilistic approach
on a Fano manifold $X$ with trivial group $\text{Aut}_{0}(X),$ introduced
in \cite{berm8 comma 5}. The starting point is the basic fact that
a Kähler--Einstein (KE) metric $\omega_{\text{KE}}$ with positive
Ricci curvature on $X$ can be readily be recovered from its (normalized)
volume form $dV_{\text{KE}}.$ The strategy of the probabilistic approach
is thus to construct the normalized volume form $dV_{KE}$ by sampling
a large number of $N$ points on $X.$ More precisely, the goal is
to construct a canonical symmetric probability measure on the $N$-fold
products $X^{N}$ --- solely defined in terms of explicit algebro-geometric
data of $X$ --- and show that $dV_{\text{KE}}$ emerges as $N\rightarrow\infty,$
when sampling $N$ points. The construction of the probability measure
on $X^{N}$ is based on the observation that $N$-fold products $X^{N}$
come with an algebraic measure $\nu^{(N)},$ canonically determined
up to scaling. More precisely, $N$ is the subsequence defined by
\begin{equation}
N=N_{k}\coloneqq\dim H^{0}(X,-kK_{X})=\frac{V}{n!}k^{n}+O(k^{n-1}),\,\,V:=\int_{X}c_{1}(X)^{n},\label{eq:def N intro}
\end{equation}
 where the \emph{level} $k$ ranges over all positive rational numbers
such that $-kK_{X}$ is a well-defined line bundle, using additive
notation for tensor products of line bundles. To simplify the notation
the subscript $k$ of $N_{k}$ will be suppressed. The measure $\nu^{(N)}$
on $X^{N}$ may be expressed as 
\begin{equation}
\nu^{(N)}\coloneqq\frac{i^{(Nn)^{2}}}{2^{Nn}}\big(\mathrm{det}S^{(N)}\big)^{-1/k}\wedge\overline{\left(\det S^{(N)}\right)^{-1/k},}\label{eq:def of nu N intro}
\end{equation}
where $n$ is the dimension of $X$ and $\det S^{(N)}$ is the holomorphic
section of the line bundle $-kK_{X^{N}}\rightarrow X^{N},$ defined
as the Slater determinant 
\begin{equation}
\left(\det S^{(N)}\right)(x_{1},x_{2},\dots,x_{N})\coloneqq\det\left(s_{i}(x_{j})\right),\label{eq:slater determinant intro}
\end{equation}
in terms of a given basis $s_{1},\dots,s_{N}$ in $H^{0}(X,-kK_{X})$
(changing the basis only has the effect of scaling $\nu^{(N)}).$
However, the total mass of $\nu^{(N)}$ may be infinite, since the
density of $\nu^{(N)}$ blows up along the zero-locus of $\det S^{(N)}.$
Accordingly, $X$ is dubbed \emph{Gibbs stable} \cite{berm8 comma 5},
if $\nu^{(N)}$ has finite total mass, for $N$ sufficiently large
(this is equivalent to an algebro-geometric criterion,  recalled in
the following section). Then
\begin{equation}
\mu^{(N)}\coloneqq\nu^{(N)}\left/\int_{X^{N}}\nu^{(N)}\right.\label{eq:def of mu N when G trivial intro}
\end{equation}
defines a \emph{canonical} probability measure on $X^{N},$ which
is symmetric (i.e. invariant under permutations of the factors of
$X^{N}).$ It is conjectured in \cite{berm8 comma 5} that, when sampling
$N$ points on $X$ according to the probability measure $\mu^{(N)},$
the corresponding \emph{empirical measure $\delta_{N}$ on $X,$}
i.e. the random discrete probability measure
\begin{equation}
\delta_{N}\coloneqq\frac{1}{N}\sum_{i=1}^{N}\delta_{x_{i}},\label{eq:emp measure}
\end{equation}
on $X,$ converges in probability as $N\rightarrow\infty,$ to the
volume form $dV_{\text{KE}}$ of a Kähler--Einstein metric $\omega_{\text{KE}}.$ 

However, a necessary condition for $X$ to be Gibbs stable is that
$\text{Aut}_{0}(X)$ be trivial. In fact, when $\text{Aut}_{0}(X)$
is nontrivial there is \emph{no} probability measure on $X^{N}$ that
is invariant under the diagonal action of $\text{Aut}_{0}(X)$ (under
a mild regularity assumption) \cite{berm14}. The probabilistic approach
in \cite{berm8 comma 5} is thus too canonical, in a certain sense.
This is in line with the fact that the group $\text{Aut}_{0}(X)$
acts nontrivially on any given Kähler--Einstein metric --- as a
consequence there is no canonical way of singling out a unique Kähler--Einstein
metric. In order to probabilistically construct a Kähler--Einstein
metric we will thus need to break the $\text{Aut}_{0}(X)$-symmetry.
This is explained in Section \ref{subsec:Symmetry-breaking-and} below,
but we first introduce an algebro-geometric notion of Gibbs \emph{poly}stability,
taking the action of the group $\text{Aut}_{0}(X)$ into account and
explore the relations to stability thresholds.

\subsection{\label{subsec:Gibbs-polystability-and}Gibbs polystability and microscopic
polystability thresholds }

We will consider the general framework of \emph{log Fano manifolds}
$(X,\Delta),$ consisting of a complex projective manifold $X$ and
an effective $\Q$-divisor $\Delta$ on $X,$ whose log anti-canonical
line bundle $-K_{(X,\Delta)}$ is an ample $\Q$-line bundle \cite{bbegz}.
This setting encompasses, in particular, Fano orbifolds \cite{g-k}.
Kähler--Einstein metrics on $(X,\Delta)$ are, in general, understood
as positive currents $\omega_{\text{KE}}$ in $c_{1}\left(-K_{(X,\Delta)}\right)$
with locally bounded potentials, that restrict to genuine Kähler--Einstein
metrics on the complement of the support of $\Delta.$ In particular,
$\omega_{\text{KE }}$ satisfies the following equation in the sense
of currents: 
\[
\text{Ric \ensuremath{\omega}}=\omega+[\Delta],
\]
 where $[\Delta]$ denotes the current of integration along the divisor
$\Delta.$ Given a log Fano manifold $(X,\Delta)$ and a positive
integer $k$ such that $-kK_{(X,\Delta)}$ defines a line bundle we
set 
\begin{equation}
N=N_{k}\coloneqq\dim H^{0}(X,-kK_{(X,\Delta)}),\,\,\,\,\G\coloneqq\text{Aut}_{0}\ensuremath{(X,\Delta)}\label{eq:def N k log intro}
\end{equation}
and we will assume that the complex Lie group $\G$ is\emph{ reductive}
(this is a necessary condition for the existence of a Kähler--Einstein
metric for $(X,\Delta)$ \cite{c-d-s}).

We start by introducing an algebro-geometric notion of \emph{Gibbs
polystability,} taking the group $\G$ into account. It is formulated
in terms of the Zariski open subvariety $(X^{N})_{\text{ss}}$ of
$X^{N},$ defined by the \emph{semistable locus} of $X^{N},$ with
respect to the ample $\Q$--line bundle $-K_{(X,\Delta)},$ in the
sense of Geometric Invariant Theory (GIT) \cite{ho}. We will denote
by $\Delta_{N}$ the divisor on $X^{N}$ induced by $\Delta$ that
is invariant under permutation of the factors of $X^{N}.$ Thus $(X^{N},\Delta_{N})$
is a\emph{ }log Fano manifold. Denote by\emph{ $\mathcal{D}^{(N)}$}
the effective $\Q$-divisor on $X^{N},$ linearly equivalent to $-K_{(X^{N},\Delta_{N})},$
cut out by the $k^{\mathrm{th}}$ root of the section $\det S^{(N)},$
defined by formula \ref{eq:slater determinant intro}, in terms of
a basis in $H^{0}(X,-kK_{(X,\Delta)}).$ To recall the algebraic definition
of Gibbs stability, introduced in \cite{berm8 comma 5}, consider
the\emph{ microscopic stability threshold} $\gamma^{(N)}(X,\Delta)$
and its limit $\gamma(X,\Delta)$ defined by
\[
\gamma^{(N)}(X,\Delta)\coloneqq\text{lct}\left(X^{N},\Delta_{N};\mathcal{D}^{(N)}\right),\,\,\,\,\gamma(X,\Delta)\coloneqq\liminf_{N\rightarrow\infty}\gamma^{(N)}(X,\Delta),
\]
 where $\text{lct}(Y,\Delta;D)$ denotes the \emph{log canonical threshold}
of a divisor $D$ on a quasi-projective log pair $(Y,\Delta)$ ---
this number is a standard algebro-geometric invariant of a divisor
--- whose valuative definition is recalled in Section \ref{subsec:Log-canonical-thresholds}.
A log Fano manifold $(X,\Delta)$ is called \emph{Gibbs stable at
level} $k$ if $\gamma^{(N_{k})}(X,\Delta)>1$ and \emph{Gibbs stable}
iff $\gamma^{(N)}(X,\Delta)>1,$ when $N$ is sufficiently large (which
is equivalent to the integrability criterion recalled in the previous
section). Furthermore, $(X,\Delta)$ is called \emph{Gibbs semistable}
if $\gamma(X,\Delta)\geq1$ and \emph{uniformly Gibbs stable} if $\gamma(X,\Delta)>1$
\cite{berm8 comma 5,f-o}. In order to define a notion Gibbs polystability
we first introduce the\emph{ microscopic polystability threshold}
$\gamma_{\text{PS}}^{(N)}(X,\Delta)$:
\[
\gamma_{\text{PS}}^{(N)}(X,\Delta)\coloneqq\text{lct}\left((X^{N})_{\text{ss}},\Delta_{N};\mathcal{D}^{(N)}\right),\,\,\,\gamma_{\text{PS}}(X,\Delta)\coloneqq\liminf_{N\rightarrow\infty}\gamma_{\text{PS}}^{(N)}(X,\Delta),
\]
where the limit $\gamma_{\text{PS}}(X,\Delta)$ is dubbed the \emph{Gibbs
polystability threshold}. We call $(X,\Delta)$ \emph{Gibbs polystable}
if it is Gibbs semistable and $\gamma_{\text{PS}}^{(N)}(X,\Delta)>1,$
when $N$ is sufficiently large and \emph{uniformly Gibbs polystable
}if the latter condition is replaced by $\gamma_{\text{PS}}(X,\Delta)>1.$
\begin{conjecture}
\label{conj:Gibbs poly intro}Let $(X,\Delta)$ be a log Fano manifold.
Then $(X,\Delta)$ is Gibbs polystable iff it is uniformly Gibbs polystable
iff $(X,\Delta)$ is K-polystable. 
\end{conjecture}

The first equivalence may, perhaps, be too optimistic, in general.
But we present arguments motivating the second one. They are based
on a conjectural Large Deviation Principle that we shall come back
to below, linking the invariant $\gamma_{\text{PS}}(X,\Delta)$ to
the \emph{analytic polystability threshold }defined by 
\begin{equation}
\delta_{\text{PS}}^{\mathrm{A}}(X,\Delta)\coloneqq1+\sup\left\{ \epsilon\in\R\colon\exists C_{\epsilon}\in\R\,\,\mathcal{M}\geq\epsilon(\mathcal{I}-\mathcal{J})^{\G}+C_{\epsilon}\right\} ,\label{eq:def of anal st thr}
\end{equation}
where $\mathcal{M}$ denotes the (log) \emph{Mabuchi functional} on
the space of Kähler metrics in $c_{1}\left(-K_{(X,\Delta)}\right),$
whose minimizers are Kähler--Einstein metrics for $(X,\Delta),$
and $(\mathcal{I}-\mathcal{J})^{\G}$ denotes the following $\G$-invariant
version of the classical non-negative functional $\mathcal{I}-\mathcal{J},$
defined as the difference of the classical energy functionals $\mathcal{I}$
and $\mathcal{J}$ in Kähler geometry 
\begin{equation}
(\mathcal{I}-\mathcal{J})^{\G}=\text{inf}_{g\in\G}g^{*}(\mathcal{I}-\mathcal{J})\label{eq:def of I min J G intro}
\end{equation}
(the definitions of $\mathcal{I}$ and $\mathcal{J}$ are recalled
in Remark \ref{rem:plurienergy is I minus J}. The analytic polystability
threshold is reminiscent of the $\G$-reduced weighted delta invariant
introduced in \cite{de-ju}, and we show that they are equal in case
of vanishing Futaki character (Prop \ref{prop: relation with reduced delta}).
In the case of a complex curve, all the three functionals $\mathcal{I}-\mathcal{J},$
$\mathcal{I}$ and $\mathcal{J}$ coincide and are proportional to
the Dirichlet norm of the Kähler potential. On a general Kähler manifold
they are all comparable and are therefore commonly used as measures
of \emph{coercivity} (or properness) in Kähler geometry \cite{ti,bbgz,d-r};
their $\G$-invariant versions appear when studying coercivity modulo
the action of $\G$ \cite{z-z,d-r}. In particular, by the solution
in \cite{d-r} of a conjecture of Tian 
\[
\delta_{\text{PS}}^{\mathrm{A}}(X,\Delta)>1\,\text{iff\,\ensuremath{(X,\Delta)\,\text{admits a KE metric}}}
\]
 and this equivalence holds regardless of which of the functionals
$\mathcal{I}-\mathcal{J},$ $\mathcal{I}$ and $\mathcal{J}$ are
used. In the present context, however, it will be important to work
with the functional $\mathcal{I}-\mathcal{J}$ as this choice leads
to a thermodynamical interpretation of the corresponding analytic
polystability threshold $\delta_{\text{PS}}^{\mathrm{A}}(X,\Delta)$
(see Prop \ref{prop:anal stab in terms of F }). It should be emphasized
that, while the definition of $\mathcal{I}-\mathcal{J}$ depends on
the choice of a reference Kähler metric $\omega_{0}\in c_{1}(X)$,
this does not affect the definition \ref{eq:def of anal st thr} of
$\delta_{\text{PS}}^{\mathrm{A}}(X,\Delta),$ since it only alters
the additive constant $C_{\epsilon}.$ 

When $\G$ is trivial the invariant $\delta_{\text{PS}}^{\mathrm{A}}(X,\Delta)$
equals the well-known analytic stability threshold $\delta^{\mathrm{A}}(X,\Delta),$
which, by \cite{zh}, coincides with the \emph{stability threshold}
$\delta(X,\Delta)$, introduced in \cite{f-o} and further developed
in \cite{bl-j}. This is an algebro-geometric (valuative) invariant
that has come to play a pivotal role in recent developments around
the Yau--Tian--Donaldson conjecture \cite{l-x-z}. It was conjectured
in \cite{berm8 comma 5} that $\delta^{\mathrm{A}}(X,\Delta)=\gamma(X,\Delta).$
Here we propose the following generalization:
\begin{conjecture}
\label{conj:lct is an stab intro}On any log Fano manifold $(X,\Delta)$
with vanishing Futaki character the invariants $\gamma_{\mathrm{PS}}(X,\Delta)$
and $\delta_{\mathrm{PS}}^{\mathrm{A}}(X,\Delta)$ coincide. 
\end{conjecture}

The validity of this conjecture would imply the second equivalence
in Conjecture \ref{conj:Gibbs poly intro}. Indeed, by \cite{d-r},
$\delta_{\text{PS}}^{\mathrm{A}}(X,\Delta)>1$ iff $(X,\Delta)$ admits
a Kähler--Einstein metric, which, by the solution of the YTD conjecture
is equivalent to K-polystability (see Theorem \ref{thm:KE equiv}).

Our first main result concerns the case of log Fano curves $(X,\Delta),$
i.e. the case when $n=1.$ In this case $X$ is the Riemann sphere,
$X=\P^{1}$ and $\Delta$ is any divisor of $r$ points on $X$ with
weights $w_{i}\in]0,1[$ satisfying 
\[
2-\sum_{i=1}^{r}w_{i}>0.
\]
When $\G$ is nontrivial it is well-known \cite{li000} that $(X,\Delta)$
is K-polystable iff $\Delta$ is the divisor $\Delta_{w}$ on $\P^{1}$
supported on two points with the same weight $w\in[0,1)$ (thus $\Delta_{w}=0$
when $w=0).$
\begin{thm}
\label{thm:gibbs poly and lct intro}Conjecture \ref{conj:Gibbs poly intro}
holds on any log Fano curve $(\P^{1},\Delta).$ More precisely, 
\[
\gamma^{(N)}(\P^{1},\Delta)=\frac{2}{V}\min\left\{ 1-1/N,\min_{l}\{1-w_{l}\}\right\} 
\]
and when the group $\G$ is nontrivial $(X,\Delta)$ is Gibbs polystable
iff $\Delta$ is of the form $\Delta_{w}$ for some $w\in[0,1[.$
Moreover, for $N>3$,
\[
\gamma_{\mathrm{PS}}^{(N)}(\P^{1},0)=\frac{N-1}{\lfloor N/2\rfloor},\,\,\,\,\gamma_{\mathrm{PS}}^{(N)}(\P^{1},\Delta_{w})=\min\left\{ \frac{2(1-1/N)}{V},\frac{N-1}{\lfloor N/2\rfloor-1}\right\} \,\text{if \,\ensuremath{w\neq0.}}
\]
For $N=2,3$, $\gamma_{\mathrm{PS}}^{(N)}(\P^{1},0)=\infty$, and
$\gamma_{\mathrm{PS}}^{(N)}(\P^{1},\Delta_{w})=\frac{2(1-1/N)}{V}$.
As a consequence, 
\begin{equation}
\gamma_{\mathrm{PS}}(\P^{1},0)=2,\,\text{and\,}\,\gamma_{\mathrm{PS}}(\P^{1},\Delta_{w})=\min\left\{ \frac{1}{1-w},2\right\} \,\text{if}\,w\neq0.\label{eq:formula for gamma in thm Gibbs poly intro}
\end{equation}
\end{thm}

The proof exploits the Fulton--MacPherson compactification of the
configuration space of $N$ points on $X.$ It also yields a description
of a sequence of valuations over the semi-stable locus $(X^{N})_{\text{ss}}$
of $X^{N}$ (for $X=\P^{1})$ asymptotically computing the log canonical
threshold $\gamma^{(N)}(\P^{1},\Delta)$ as $N\rightarrow\infty$.
These are obtained by blowing up a configuration of points $(x_{1}^{*},\dots,x_{N}^{*})$
in $(X^{N})_{\text{ss}}$ concentrating, as $N\rightarrow\infty,$
at a single point $x_{0}$ in $\P^{1}$ when $w\in]0,1/2[$ and on
two distinct points $x_{1}$ and $x_{2}$ when $w\notin]0,1/2[.$
More precisely, the following weak convergence of measures on $X$
holds:
\begin{equation}
\frac{\sum_{i=1}^{N}\delta_{x_{i}^{*}}}{N}\rightarrow\delta_{x_{0}}\,\,\text{when\, }w\in]0,1/2[,\,\,\,\frac{\sum_{i=1}^{N}\delta_{x_{i}^{*}}}{N}\rightarrow\frac{1}{2}\delta_{x_{1}}+\frac{1}{2}\delta_{x_{2}}\,\,\text{when\, }w\notin]0,1/2[\label{eq:conv optimizers intro}
\end{equation}

Coming back to Conjecture \ref{conj:lct is an stab intro} for a general
log Fano manifold, we introduce a slight modification of the log canonical
threshold $\gamma_{\text{PS}}^{(N)}(X,\Delta),$ denoted by $\gamma_{\epsilon_{N}}^{(N)}(X,\Delta),$
for a small parameter $\epsilon_{N}$ of the order $O(1/N).$ It is
obtained by replacing the semi-stable locus $(X^{N})_{\text{ss}}$
by a Zariski open subset $(X^{N})_{\text{ss},\epsilon_{N}}$ of $X^{N}$
containing slightly unstable configurations (see Definition \ref{def:strong}).
This leads to a notion of \emph{strong uniform Gibbs }\emph{polystability},
defined by replacing $(X^{N})_{\text{ss}}$ with $(X^{N})_{\text{ss},\epsilon_{N}}$
(which is expected to be equivalent to uniform Gibbs polystability,
as illustrated by the case of log Fano curves; see Example \ref{exa:strong curve}).
Our second main result yields, in particular, an effective criterion
for the existence of a Kähler--Einstein metric on Fano manifolds: 
\begin{thm}
\label{thm:uniform polyGibbs implies KE}Let $X$ be a Fano manifold
whose Futaki character vanishes. There exist effective constants $k_{1}$
and $C_{1}$ such that if
\[
\gamma_{\epsilon_{N}}^{(N)}(X)\geq1+C_{1}k^{-1},
\]
for some $k\geq k_{1},$ then $X$ admits a Kähler--Einstein metric
(and, as a consequence, $X$ is K-polystable). In particular, if $X$
is strongly uniformly Gibbs polystable, then $X$ admits a Kähler--Einstein
metric. More generally, there exist effective constants $k_{2}$ and
$C_{2}$ such that for any $k\geq k_{2}$
\[
-C_{2}k^{-1}+\frac{1+k^{-1}}{1/\gamma_{\epsilon_{N}}^{(N)}(X)+k^{-1}}\leq\delta_{\text{PS}}^{\mathrm{A}}(X).
\]
\end{thm}

Theorem \ref{thm:uniform polyGibbs implies KE} generalizes \cite[Cor 2.6]{f-o},
where it was shown --- using purely algebro-geometric techniques
--- that $\gamma^{(N)}(X)$ is bounded from above by Fujita--Odaka's
invariants $\delta_{k}(X)$ converging to the algebraic stability
threshold $\delta(X)$ (which, by \cite{zh}, coincides with $\delta^{A}(X)$)
and, as a consequence,
\[
\limsup_{N\rightarrow\infty}\gamma^{(N)}(X)\leq\delta(X)
\]
(and essentially the same arguments apply to log Fano varieties).
However, it seems challenging to extend the algebro-geometric definition
of $\delta_{k}(X)$ in \cite{f-o} to incorporate the action of the
group $\G,$ so that $\delta_{\text{PS}}^{\mathrm{A}}(X)$ is obtained
as $k\rightarrow\infty.$ \footnote{There exist reduced versions of $\delta(X,\Delta),$ defined by restricting
to torus invariant valuations \cite[Appendix A]{x-z} (or $\G$-invariant
valuations \cite{li}), but they do not coincide with $\delta_{\text{PS}}^{\mathrm{A}}(X,\Delta),$
as can be checked on $\P^{1}.$ } Instead, our proof of Theorem \ref{thm:uniform polyGibbs implies KE}
adapts the analytic approach developed in \cite{berm13}, by proving
an effective large deviation bound (Theorem \ref{thm:effective LD bound with Ding}),
motivated by the conjectural Large Deviation Principle alluded to
above. Extending Theorem \ref{thm:uniform polyGibbs implies KE} to
orbifold pairs $(X,\Delta)$ is rather straightforward, and we expect
that it can be extended to general log Fano manifolds $(X,\Delta),$
but a generalization seems to hinge on missing Bergman kernel asymptotics
for metrics with edge-cone singularities, of independent interest. 

An important feature of Theorem \ref{thm:uniform polyGibbs implies KE}
is that it yields an effective lower bound on $\delta_{\text{PS}}^{\mathrm{A}}(X)$
in terms of quantities that are, at least in principle, computable
by a finite algebraic procedure. Indeed, the invariant $\gamma_{\epsilon_{N}}^{(N)}(X)$
is a log canonical threshold, and hence can be computed in finitely
many algebraic steps. Moreover, the constants appearing in Theorem
\ref{thm:uniform polyGibbs implies KE} are effective in the same
sense. This is explained in \cite{a-b4}, where it is shown --- assuming
the strong form of Conjecture \ref{conj:Gibbs poly intro} --- that
the problem whether a Fano manifold admits a Kähler--Einstein metric
is \emph{decidable}. More precisely, there is an algorithm, involving
only standard algebraic operations, which terminates after finitely
many steps and returns the answer to the question of whether $X$
admits a Kähler--Einstein metric. Moreover, it is shown in \cite{a-b4}
that decidability holds unconditionally on any Fano manifold with
discrete automorphism group (building on the upper bounds on $\delta_{k}(X)$
in \cite{berm13}, refining bounds in \cite{zh}). 

In general, the log canonical threshold of a $\Q$-divisor can be
computed from a log resolution. Finding an explicit log resolution
of $\left((X^{N})_{\text{ss}},\Delta_{N};\mathcal{D}^{(N)}\right)$
seems, however, very challenging beyond the case of log Fano curves.
Recent reinforcement-learning approaches to resolution of singularities
\cite{b-f-z} nevertheless suggest a possible future role for machine-learning-guided
searches for efficient log resolutions, or for valuations computing
the corresponding log canonical thresholds.

In the case when $X$ is a K-polystable Fano curve, i.e. $X=\P^{1}$,
it follows directly from combining Theorem \ref{thm:uniform polyGibbs implies KE}
with the formula for $\gamma_{\text{PS}}^{(N)}(\P^{1})$ in Theorem
\ref{thm:gibbs poly and lct intro}, and using that $\gamma_{\epsilon_{N}}^{(N)}(\P^{1})=\gamma_{\mathrm{PS}}^{(N)}(\P^{1})+o(1)$
(see Example \ref{exa:strong curve}), that 
\[
\delta_{\text{PS}}^{\mathrm{A}}(\P^{1})\geq2,\,\,\,\text{i.e. \,\ensuremath{\mathcal{M}\geq\epsilon\mathcal{J}^{\G}}+\ensuremath{C_{\epsilon}}}\,\,\,\forall\epsilon\in[0,1[.
\]

\subsection{Results about analytic polystability thresholds established in \cite{a-b-s3}}

In the companion paper \cite{a-b-s3} we establish the following general
formula for $\delta_{\text{PS}}^{\mathrm{A}}(\P^{1},\Delta_{w})$,
answering a question posed in \cite{de} in case $w=0$. When combined
with formula \ref{eq:formula for gamma in thm Gibbs poly intro},
it shows that Conjecture \ref{conj:lct is an stab intro} holds on
any log Fano curve.
\begin{thm*}
\cite{a-b-s3} For any log Fano curve $(\P^{1},\Delta)$
\[
\delta_{\mathrm{PS}}^{\mathrm{A}}(\P^{1},0)=2,\,\,\,\text{and\,\,\,}\delta_{\text{PS}}^{\mathrm{A}}(\P^{1},\Delta_{w})=\min\left\{ \frac{1}{1-w},2\right\} \,\text{if}\,w\neq0.
\]
\end{thm*}
The proof of the upper bound on $\delta_{\text{PS}}^{\mathrm{A}}(\P^{1},\Delta_{w})$
in \cite{a-b-s3} is inspired by the convergence \ref{eq:conv optimizers intro}.
For example, when $w\notin]0,1/2[$ the proof shows that for any $\epsilon>1$,
$\mathcal{M}-\epsilon\mathcal{J}^{\G}$ is destabilized by the weak
geodesic ray in the $C^{1,1}$-closure of the space of Kähler metrics
in $c_{1}(\P^{1},\Delta_{w})$ corresponding to a test configuration
for $(\P^{1},\Delta_{w})$ whose central fiber has \emph{two} irreducible
components. As for the proof of the lower bound on $\delta_{\text{PS}}^{\mathrm{A}}(\P^{1},\Delta_{w})$
it exploits a generalization of Tian's $\alpha$-invariant taking
the action of the non-compact group $\G$ into account. Additionally,
in the case of $w=0,$ we establish in \cite{a-b-s3} the following
sharp inequality:
\begin{equation}
\mathcal{M}(\omega)\geq\mathcal{J}^{\G}(\omega),\label{eq:sharp ineq for M on P one intrp}
\end{equation}
when $\mathcal{M}$ and $\mathcal{J}$ are normalized so that they
vanish on the standard $\mathrm{SU}(2)$-invariant Kähler metric $\omega_{0}$
in $c_{1}(\P^{1}).$ Moreover, equality holds in \ref{eq:sharp ineq for M on P one intrp}
iff $\omega$ is Kähler--Einstein, i.e. $\omega\in\G\omega_{0}.$
The proof of the inequality \ref{eq:sharp ineq for M on P one intrp}
builds on the formula $\delta_{\text{PS}}^{\mathrm{A}}(\P^{1},0)=2$
and shows that the constants $C_{\epsilon}$ in formula \ref{eq:def of anal st thr}
in the definition of $\delta_{\text{PS}}^{\mathrm{A}}(\P^{1},0)$
can, in this case, be taken to vanish when $\epsilon\in[0,1].$ As
a consequence, as explained in \cite{a-b-s3}, the inequality \ref{eq:sharp ineq for M on P one intrp}
yields quantitative stability results for the sharp logarithmic Hardy--Littlewood--Sobolev
inequality on the two-sphere with optimal stability constants, improving
on \cite{c-f} and \cite{ca2}. 

\subsection{Higher dimensional examples}

Beyond complex curves, Gibbs stability has turned out to be notoriously
difficult to establish, except for the case of exceptional log Fano
manifolds, where it holds automatically at any level $k.$ In Section
\ref{subsec:Gibbs-stable-orbifolds} we provide the first non-exceptional
higher-dimensional examples of Gibbs stable log Fano varieties. More
precisely, we show that the Fano orbifold defined as the quotient
of the Fermat hypersurface $X_{n+1}$ of degree $(n+1)$ in $\P^{n+1},$
under the standard action of the group $(\Z/(n+1)\Z)^{n},$ is Gibbs
stable at \emph{any} level $k.$ This Fano orbifold may be identified
with a log pair $(\P^{n},\Delta),$ where $\Delta$ is a divisor that
assigns the same weight $1-1/(n+1)$ to $n+2$ normally crossing hyperplanes
in $\P^{n}.$ Although the group $\text{Aut}_{0}(\P^{n},\Delta)$
is trivial, the proof of the Gibbs stability of $(\P^{n},\Delta)$
exploits that the group $\text{Aut}_{0}(\P^{n})$ is nontrivial. Furthermore,
in the companion paper \cite{a-b-s2} we exhibit log Fano manifolds
$(\P^{2},\Delta)$ with nontrivial group $\text{Aut}_{0}(\P^{n},\Delta)$
that are Gibbs \emph{polystable} at the minimal level (i.e. for the
smallest $k\in]0,\infty[$ such that $-kK_{(X,\Delta)}$ defines a
line bundle): 
\begin{thm}
\label{thm: Gibbs polystable example}\cite{a-b-s2} Let $\Delta_{w}$
be the torus-invariant divisor on $\P^{2}$ assigning the same weight
$w\in[0,1[$ to each of the three torus-invariant hyperplanes. Then
the log Fano manifold $(\P^{2},\Delta_{w})$ is Gibbs polystable at
the minimal level iff $w>3/4$ or $w=0$. 
\end{thm}

The conclusion that $\mathbb{P}^{2}$ is Gibbs polystable at the minimal
level, but $(\mathbb{P}^{2},\Delta_{w})$ is not Gibbs polystable
for $w>0$ sufficiently small shows that the discontinuity phenomenon
observed for $(\mathbb{P}^{1},\Delta_{w})$ in Theorem \ref{thm:gibbs poly and lct intro}
persists in higher dimensions, at least for this example and at the
minimal level. The case of $w=0$ in Theorem \ref{thm: Gibbs polystable example}
is a special case of Proposition \ref{prop: Gibbs poly of minimal Pn},
stating that $\mathbb{P}^{n}$ is Gibbs polystable at the minimal
level for any $n$. When $w=1-1/m$ for a positive integer $m,$ the
log Fano manifold $(\P^{2},\Delta_{w})$ may be identified with the
Fano orbifold defined as the quotient of $\P^{2}$ under the standard
action of $(\Z/m\Z)^{2}.$ In this case the previous theorem applies
to the AdS/CFT correspondence, where the Fano orbifold in question
encodes a particular toric quiver gauge theory. 

\subsection{\label{subsec:Symmetry-breaking-and}Explicit symmetry breaking and
the emergence of a unique Kähler--Einstein metric}

Let $(X,\Delta)$ be a log Fano manifold and fix a maximal compact
subgroup $\K$ of the corresponding symmetry group $\G$ and a $\K$-invariant
metric on $-K_{(X,\Delta)}$ with positive curvature. The curvature
defines a Kähler form $\omega_{0}$ in $c_{1}(-K_{(X,\Delta)})$ and
thus a $\K$-invariant symplectic form on $X.$ As a consequence,
the natural action of $\K$ on $-K_{(X,\Delta)}$ gives rise to a
\emph{moment map }
\begin{equation}
\boldsymbol{m}\colon X\rightarrow\mathfrak{k}^{*}\label{eq:moment map intro}
\end{equation}
from $X$ into the dual of the Lie algebra $\mathfrak{k}$ of $\K.$
We will show that if $(X,\Delta)$ admits some Kähler--Einstein metric,
then there exists a \emph{unique }such metric whose normalized volume
form $\mu_{\text{KE}}$ has \emph{vanishing moment, }in the sense
that 
\begin{equation}
\mu_{\text{KE}}\in\mathcal{P}(X)_{0}\coloneqq\left\{ \mu\in\mathcal{P}(X)\colon\int_{X}\boldsymbol{m}\mu=0\right\} ,\label{eq:def of P X zero intro}
\end{equation}
where $\mathcal{P}(X)$ denotes the space of all probability measures
$\mu$ on $X$ (see Theorem \ref{thm:KE equiv}). The volume form
$\mu_{\text{KE}}$ is $\K$-invariant, and we will propose an essentially
explicit probabilistic construction of $\mu_{\text{KE}}.$ To this
end denote by $\boldsymbol{m}_{N}$ the moment map of the diagonal
action of the group $\text{\K}$ on $X^{N},$ normalized so that
\[
\boldsymbol{m}_{N}(x_{1},\dots,x_{N})=\frac{1}{N}\sum_{i=1}^{N}\boldsymbol{m}(x_{i})\in\mathfrak{k}^{*}.
\]
Consider first the case when $0$ is a regular value of $\boldsymbol{m}_{N}.$
We then denote by $\nu_{\boldsymbol{m}_{N}}^{(N)}$ the measure on
$\left\{ \boldsymbol{m}_{N}=0\right\} \subset X^{N},$ induced by
the measure $\nu^{(N)}$ on $X^{N}$ (defined by formula \ref{eq:def of nu N intro}
when $\Delta$ is trivial) --- see Section \ref{subsec:Probabilistic-setup}
for the general definition. The algebro-geometric definition of Gibbs
polystability of $(X,\Delta)$, introduced in the previous section,
is equivalent to the Gibbs semistability of $(X,\Delta)$ together
with the analytic assumption that the measure $\nu_{\boldsymbol{m}_{N}}^{(N)}$
has finite total mass (Prop \ref{prop:Gibbs poly iff integrable}),
i.e. that
\begin{equation}
\mu_{0}^{(N)}\coloneqq\frac{1}{Z_{N,0}}\nu_{0}^{(N)},\,\,\,\,Z_{N,0}\coloneqq\int\nu_{0}^{(N)}\label{eq:def of mu N broken intro}
\end{equation}
defines a probability measure on $X^{N},$ supported on $\left\{ \boldsymbol{m}_{N}=0\right\} .$
The probability measure $\mu_{0}^{(N)}$ is symmetric, i.e. invariant
under permutations of $X^{N}$ and $\K$-invariant, but not $\G$-invariant.
This is an instance of \emph{explicit symmetry breaking}, as the original
$\G$-symmetry has been reduced to a $\K$-symmetry. Accordingly,
the probability measure \ref{eq:def of mu N broken intro} on $X^{N}$
shall be called the\emph{ $\K$-reduced canonical probability measure}.\emph{ }

We can now view the empirical measure $\delta_{N}$ (formula \ref{eq:emp measure})
as a random variable on the probability space $\left(X^{N},\mu_{0}^{(N)}\right).$
Compared to the canonical probability measure \ref{eq:def of mu N when G trivial intro}
on $X^{N}$, this means that we are conditioning the $N$ points on
$X^{N}$ to have vanishing moment map $\boldsymbol{m}_{N}.$ The empirical
measure $\delta_{N}$ of any such configuration of points has vanishing
moment, i.e. $\delta_{N}$ maps $X^{N}$ into the space $\mathcal{P}(X)_{0}$
defined in formula \ref{eq:def of mu N broken intro}. 

Before stating our two main conjectures we note that in contrast to
the $\K$--reduced probability measure $\mu_{0}^{(N)},$ the measure
$\nu_{0}^{(N)}$ depends on a choice of bases in $H^{0}(X,-kK_{(X,\Delta)}).$
As a consequence, so does its total mass $Z_{N,0}.$ In what follows
we will take the basis to be orthonormal with respect to the $L^{2}$--norm
on $H^{0}(X,-kK_{(X,\Delta)})$ induced by the fixed metric on $-K_{(X,\Delta)}.$
The corresponding constant $Z_{N,0}$ shall then be called the\emph{
$\K$-reduced partition function, }in view of its statistical mechanical
interpretation described in Section \ref{subsec:Probabilistic-setup}.
The fixed metric on $-K_{(X,\Delta)}$ also determines the Mabuchi
functional $\mathcal{M}$ on the space of Kähler metrics in $c_{1}(-K_{(X,\Delta)})$
(which, a priori, is only defined up to an additive constant).
\begin{conjecture}
\label{conj:large N limtis intro}Let $(X,\Delta)$ be a log Fano
manifold and assume that $(X,\Delta)$ admits a Kähler--Einstein
metric. Then 
\begin{equation}
\text{(i)}-\lim_{N\rightarrow\infty}N^{-1}\log Z_{N,0}=\inf\mathcal{M}.\label{eq:conv towards Mab in conj intro}
\end{equation}
As a consequence, $Z_{N,0}<\infty$ for $N$ sufficiently large and,
viewing the empirical measure $\delta_{N}$ as a random variable on
$\left(X^{N},\mu_{0}^{(N)}\right),$
\begin{equation}
\text{(ii) }\lim_{N\rightarrow\infty}\delta_{N}=\mu_{\text{KE}},\label{eq:conv to KE conj intro}
\end{equation}
in probability, where $\mu_{\text{KE}}$ is the normalized volume
form of the unique KE metric for $(X,\Delta)$ with vanishing moment,
i.e. $\mu_{\text{KE}}\in\mathcal{P}(X)_{0}.$
\end{conjecture}

This conjecture is motivated by a heuristic mean field type argument,
outlined in Section \ref{subsec:Motivation:-A-conjectural}, suggesting
the following Large Deviation Principle (LDP) on the space $\mathcal{P}(X)_{0},$
as $N\rightarrow\infty$: 
\begin{equation}
(\delta_{N})_{*}\nu_{0}^{(N)}\sim e^{-NF(\mu)},\label{eq:formal LDP intro}
\end{equation}
 where $F$ is the free energy functional, introduced in \cite{berm6},
that may be identified with the Mabuchi functional on the space of
all Kähler metrics $\omega$ in $c_{1}(-K_{(X,\Delta)})$ via the
relation $F(\omega^{n}/V)=\mathcal{M}(\omega)$. We show that the
validity of such a LDP would imply both the asymptotics \ref{eq:conv towards Mab in conj intro}
and the convergence \ref{eq:conv to KE conj intro}, using that $\mu_{\text{KE}}$
is the unique minimizer of $F(\mu)$ on $\mathcal{P}(X)_{0}$. Moreover,
by introducing temperature into the picture the corresponding LDP
would also imply Conjecture \ref{conj:lct is an stab intro}, using
the observation that the infimum over $\G$ in formula \ref{eq:def of I min J G intro}
is attained at the unique Kähler metric in a given $\G$-orbit whose
volume form has vanishing moment (see Prop \ref{lem:inf of E over G}). 

To bypass the assumption that $0$ be a regular value of $\boldsymbol{m}_{N}$
(and to facilitate some analytic arguments) we will also consider
a setup where $\left\{ \boldsymbol{m}_{N}=0\right\} $ is replaced
by $\left\{ \left|\boldsymbol{m}_{N}\right|<\epsilon_{N}\right\} $
for an appropriate \emph{thickening parameter} $\epsilon_{N}$ tending
to zero. Using this setup with $\epsilon_{N}=O(1/N)$ we deduce Theorem
\ref{thm:uniform polyGibbs implies KE} from an effective large deviation
lower bound on the corresponding \emph{thickened partition function
}defined by
\[
Z_{N,\epsilon_{N}}\coloneqq\int_{\left\{ \left|\boldsymbol{m}_{N}\right|<\epsilon_{N}\right\} }\nu_{0}^{(N)}
\]
(see Theorem \ref{thm:effective LD bound with F}). This bound also
implies that, if $(X,\Delta)$ admits a Kähler--Einstein metric,
then 
\[
\limsup_{N\rightarrow\infty}-N^{-1}\log Z_{N,\epsilon_{N}}\leq\inf\mathcal{M},
\]
when $\epsilon_{N}\gg N^{-1/2},$ in line with the asymptotics in
Conjecture \ref{conj:large N limtis intro} $(i)$ (see Cor \ref{cor:asymptotic upper bound on minus log Z text}).
The latter asymptotics are, in the case of log Fano curves, completely
established in the companion paper \cite{a-b-s1} with an explicit
convergence rate:
\begin{thm}
\cite{a-b-s1} \label{thm:conv towards inf for n one intro}The convergence
\ref{eq:conv towards Mab in conj intro} holds on any K-polystable
log Fano curve $(X,\Delta)$ when $N$ is taken as the sequence of
odd positive integers. More precisely, 
\[
-N^{-1}\log Z_{N,0}=\inf\mathcal{M}+O\left(\frac{\log N}{N}\right).
\]
\end{thm}

The assumption that $N$ be odd is equivalent to $0$ being a regular
value of $\left\{ \boldsymbol{m}_{N}=0\right\} .$ As for the conjectured
convergence towards $\mu_{\text{KE}}$ in Conjecture \ref{conj:large N limtis intro}
$(ii)$ we will also show in \cite{a-b-s1}, that it holds up to an
arbitrarily small ``temperature deformation'', naturally appearing
in the context of Onsager's point vortex model on the two-sphere (see
Remark \ref{rem:vortex}). 

\subsection{Outlook on non-Archimedean aspects }

In the work in progress \cite{a-b-b-j} a non-Archimedean approach
to Gibbs polystability is introduced, generalizing the approach to
Gibbs stability outlined in \cite[Section 5]{ber1}. In particular,
it is shown that uniform Gibbs polystability implies K-polystability.
However, in contrast to Theorem \ref{thm:uniform polyGibbs implies KE},
it seems hard to obtain effective results using the non-Archimedean
approach, or proving that $\gamma_{\text{PS}}(X,\Delta)\leq\delta_{\text{PS}}^{\mathrm{A}}(X,\Delta).$

\subsection{Further relations to previous work}

The formula for $\gamma^{(N)}(\P^{1},\Delta)$ in Theorem \ref{thm:gibbs poly and lct intro}
was established in \cite[Section 3.1]{Fu}, when $\Delta=0,$ using
results about multiplier ideal sheaves of hyperplane arrangements.
The quantitative lower bound on $\delta_{\text{PS}}^{\mathrm{A}}(X)$
in Theorem \ref{thm:uniform polyGibbs implies KE} may be compared
with the asymptotics $\delta_{k}(X)(1-O(1/k))=\delta(X)$ for Fujita--Odaka's
invariants $\delta_{k}(X)$ established in \cite{j-r-t} under the
condition that $\delta(X)$ is computed by a divisorial valuation,
and subsequently proved unconditionally in \cite{pe}. However, the
corresponding lower bound on $\delta(X)$ in \cite{j-r-t,pe} is not
effective, as it depends on the valuation $v$ computing $\delta(X)$
(see \cite[Thm 5.1]{j-r-t}). On the other hand, the results in \cite{j-r-t,pe}
apply when $X$ is singular, while extending Theorem \ref{thm:uniform polyGibbs implies KE}
to the singular case would require extending the upper bound on the
Bergman measures in Prop \ref{prop:bound on Bergman} to singular
Fanos (for appropriate continuous metrics on $-K_{X}).$ This appears
to be an open problem of independent interest.

\subsection{Acknowledgements}

Thanks to Sébastien Boucksom and Mattias Jonsson for discussions and
to Thibaut Delcroix and Simon Jubert for pointing out their related
work \cite{de-ju}. AI (ChatGPT) was used in polishing this manuscript,
and aided in proving Lemma \ref{lem:A N is finite}. This work was
supported by a Wallenberg Scholar grant from the Knut and Alice Wallenberg
foundation.

\subsection{Organization of the paper }
\begin{itemize}
\item Section \ref{sec:Gibbs-polystability-of} is concerned with algebro-geometric
aspects of the $N$-fold product $X^{N},$ leveraging the Fulton--MacPherson
compactification. In particular, Theorem \ref{thm:gibbs poly and lct intro}
is established, giving explicit formulas for the microscopic stability
thresholds $\gamma^{(N)}(X,\Delta)$ and $\gamma_{\text{PS}}^{(N)}(X,\Delta)$
of log Fano curves. As a consequence, Conjecture \ref{conj:Gibbs poly intro}
(linking Gibbs polystability to K-polystability) holds on log Fano
curves. Some higher-dimensional examples are also provided in Section
\ref{subsec:Higher-dimensional-examples}.
\item In Section \ref{sec:Symmetry-breaking,-thermodynamic}, where the
role of $X^{N}$ is played by the space $\mathcal{P}(X)$ of all probability
measures on $X,$ the structure of the space $\mathcal{P}(X)_{0}$
of all probability measures on $X$ with vanishing moment is explored.
This leads to a thermodynamical formula for the analytic polystability
threshold $\delta_{\text{PS}}^{\mathrm{A}}(X,\Delta)$ and new analytic
characterizations of log Fano manifolds admitting Kähler--Einstein
metrics (Theorem \ref{thm:KE equiv}). 
\item In Section \ref{sec:The-reduced-probability} we study the $\K$-reduced
probability measures $\mu_{0}^{(N)}$ on $X^{N}$ and the corresponding
partition function $Z_{N,0}$ appearing in Conjecture \ref{conj:large N limtis intro}.
The asymptotics as $N\rightarrow\infty$ are described in an overarching
conjectural Large Deviation Principle, formulated in Conjecture \ref{conj:LDP},
that we show implies both Conjecture \ref{conj:Gibbs poly intro}
and Conjecture \ref{conj:lct is an stab intro}.
\item In Section \ref{sec:Uniform-Gibbs-polystability} we prove an effective
large deviation bound for the thickened partition functions $Z_{N,\epsilon_{N}}$
that implies Theorem \ref{thm:uniform polyGibbs implies KE}
\end{itemize}

\section{\label{sec:Gibbs-polystability-of}Gibbs polystability and microscopic
stability thresholds of log Fano manifolds}

In this section we will calculate the microscopic stability thresholds
of log Fano curves appearing in Theorem \ref{thm:gibbs poly and lct intro}.
This will prove Conjecture \ref{conj:Gibbs poly intro} on log Fano
curves, linking Gibbs polystability to K-polystability in this case. 

\subsection{\label{subsec:Setup Gibbs}Setup}

We start by describing in more detail the general setup, introduced
in Section \ref{sec:Gibbs-polystability-of}. Let $(X,\Delta)$ be
an $n$-dimensional\emph{ log Fano manifold}, i.e. $X$ is a complex
projective manifold and $\Delta$ is an effective $\Q$-divisor on
$X$ such that the log anticanonical line bundle $-K_{(X,\Delta)}$
of $(X,\Delta)$ is ample. Recall that $K_{(X,\Delta)}$ is the $\Q$-line
bundle defined as the tensor product of the canonical line bundle
$K_{X}$ of $X$ and the $\Q$-line bundle defined by $\Delta,$ i.e.
$K_{(X,\Delta)}\coloneqq K_{X}+\Delta,$ using the standard additive
notation for tensor products of line bundles that we shall adopt.
Fix a positive integer $k$ such that $k\Delta$ is a divisor ---
rather than a mere $\Q$-divisor --- so that $-kK_{(X,\Delta)}$
defines a line bundle and set
\[
N\coloneqq\dim H^{0}(X,-kK_{(X,\Delta}).
\]
Denote by $s_{\Delta}$ a (multi-valued) section cutting out the $\Q$-divisor
$\Delta$ and by $\Delta_{N}$ the $\Q$-divisor on the $N$-fold
product $X^{N}$ cut out by $\bigotimes_{i=1}^{N}\pi_{i}^{*}s_{\Delta}$,
where $\pi_{i}\colon X^{N}\rightarrow X$ is the projection to the
$i^{\mathrm{th}}$ factor. Thus $(X^{N},\Delta_{N})$ defines a log
Fano manifold. We will be concerned with the effective $\Q$-divisor
$\mathcal{D}^{(N)}$ on $X^{N},$ linearly equivalent to $-K_{(X^{N},\Delta_{N})},$
whose support consists of all configurations $(x_{1},\dots,x_{N})$
of $N$ points on $X$ which are in ``bad position'' with respect
to $H^{0}(X,-kK_{(X,\Delta)})$: 
\begin{equation}
\text{Supp}\left(\mathcal{D}^{(N)}\right)\coloneqq\left\{ (x_{1},\dots,x_{N})\colon\exists s\in H^{0}(X,-kK_{(X,\Delta)})-\{0\}\text{ s.t. }s(x_{1})=\cdots=s(x_{N})=0\right\} .\label{eq:support of anti-canon divisor}
\end{equation}
 Concretely, $k\mathcal{D}^{(N)}$ is the divisor on $X^{N}$ cut
out by the section $\det S^{(N)}$ of $-kK_{(X^{N},\Delta_{N})}\rightarrow X^{N}$
defined by
\begin{equation}
(\det S^{(N)})(x_{1},x_{2},\dots,x_{N})\coloneqq\det(s_{i}(x_{j})),\label{eq:slater determinant text}
\end{equation}
in terms of a fixed basis $s_{1},\dots,s_{N}$ in $H^{0}(X,-kK_{(X,\Delta)}).$
We will denote by $\G$ the complex Lie group defined by 
\[
\G\coloneqq\text{Aut}_{0}(X,\Delta),
\]
consisting of automorphisms of $X,$ homotopic to the identity, preserving
$\Delta.$ The group $\G$ naturally acts (diagonally) on $X^{N}$
and it also acts on $H^{0}(X,-kK_{(X,\Delta)}),$ using that the action
of $\G$ naturally lifts to $-kK_{(X,\Delta)}.$

\subsubsection{Microscopic stability thresholds and Gibbs polystability}

In \cite{berm8 comma 5}, the following \emph{microscopic stability
threshold} $\gamma^{(N)}(X,\Delta)$ was introduced: 
\begin{equation}
\gamma^{(N)}(X,\Delta)\coloneqq\text{lct}\left(X^{N},\Delta_{N};\mathcal{D}^{(N)}\right),\,\,\,\,\gamma(X,\Delta)\coloneqq\liminf_{N\rightarrow\infty}\gamma^{(N)}(X,\Delta)\label{eq:def of g amma N text}
\end{equation}
expressed as the log canonical threshold of $\mathcal{D}^{(N)}$ on
$(X^{N},\Delta_{N})$ (whose definition is recalled in Section \ref{subsec:Log-canonical-thresholds}
below). Following \cite{berm8 comma 5}, $(X,\Delta)$ is called \emph{Gibbs
stable at level} $k$ if $\gamma^{(N)}(X,\Delta)>1$ and \emph{Gibbs
stable} iff $\gamma^{(N)}(X,\Delta)>1,$ when $N$ is sufficiently
large. Furthermore, $(X,\Delta)$ is called \emph{Gibbs semistable}
if $\gamma(X,\Delta)\geq1$ and \emph{uniformly Gibbs stable} if $\gamma(X,\Delta)>1$
\cite{berm8 comma 5,f-o}. If $(X,\Delta)$ is Gibbs stable at some
level $k,$ then the group $\G$ is necessarily trivial \cite{berm14}.
Assuming that $\G$ is reductive, we define the \emph{microscopic
polystability threshold} $\gamma_{\text{PS}}^{(N)}(X,\Delta)$ by
\begin{equation}
\gamma_{\text{PS}}^{(N)}(X,\Delta)\coloneqq\text{lct}\left((X^{N})_{\text{ss}},\Delta_{N};\mathcal{D}^{(N)}\right),\,\,\,\,\gamma_{\text{PS}}(X,\Delta)\coloneqq\liminf_{N\rightarrow\infty}\gamma_{\text{PS}}^{(N)}(X,\Delta)\label{eq:def of gamma N G text}
\end{equation}
 where $(X^{N})_{\text{ss}}$ denotes the Zariski open subset of $X^{N}$
defined as the semistable locus of $X^{N}$ with respect to the natural
action of the complex reductive Lie group $\G$ on the polarized manifold
$(X^{N},-K_{(X^{N},\Delta_{N})})$ (whose definition is recalled in
Section \ref{subsec:Semi-stability-in-GIT} below). We propose the
following
\begin{defn}
\label{def:Gibbs polystability}A log Fano manifold $(X,\Delta)$
is called \emph{(uniformly) Gibbs polystable} if it is Gibbs semistable
and $\gamma_{\text{PS}}^{(N)}(X,\Delta)>1$ for $N$ sufficiently
large ($\gamma_{\text{PS}}(X,\Delta)>1$). 
\end{defn}

The reductivity of the group $\G$ is thus part of the definition
of (uniform) Gibbs polystability. We expect that, in the previous
definition, the condition of Gibbs semistability can be replaced by
the condition that the \emph{Futaki character} $\text{Fut}_{(X,\Delta)}$
of $(X,\Delta)$ vanishes (the definition of $\text{Fut}_{(X,\Delta)}$
and its quantizations $\text{Fut}_{(X,\Delta),k}$ is recalled in
Section \ref{subsec:The-Futaki-invariants}). This would follow from
the validity of the Large Deviation Principle appearing in Conjecture
\ref{conj:LDP} (and we will show that it is the case on log Fano
curves). Accordingly, it seems natural to also introduce the following
definition, involving the \emph{quantized Futaki character} $\text{Fut}_{(X,\Delta),k}$:
\begin{defn}
A log Fano manifold $(X,\Delta)$ is called \emph{Gibbs polystable
at level $k$} if $\gamma_{\mathrm{PS}}^{(N_{k})}(X,\Delta)>1$ and
the quantized Futaki character $\text{Fut}_{(X,\Delta),k}$ vanishes.
\end{defn}

This definition also applies to some positive \emph{rational} values
for $k,$ as long as $-kK_{(X,\Delta)}$ is defined as a line bundle.
For example, when $X=\P^{n}$, $-kK_{X}$ is defined as a line bundle
iff $k(n+1)$ is an integer, since $-K_{X}=\mathcal{O}(n+1).$ 
\begin{rem}
When $\G$ is nontrivial we do not expect that Gibbs polystability
is equivalent to Gibbs polystability at level $k$ for all sufficiently
large positive integers $k$ (or even suitable rational $k).$ Indeed,
there are examples of toric Fano manifolds, where the Futaki character
$\text{Fut}_{(X,\Delta)}$ vanishes, while its quantizations do not
vanish. However, we do expect that in the definitions of both $\gamma(X,\Delta)$
and $\gamma_{\text{PS}}(X,\Delta)$ the liminf can be taken either
over all positive integers $k$ such that $k\Delta$ is a divisor,
or over all rational $k$ for which $-kK_{(X,\Delta)}$ is defined
as a line bundle (for example, we will show that this is the case
on log Fano curves). Anyhow, in the definitions we will, for simplicity,
demand that $N$ only ranges over all $N_{k}$ for which $k$ is a
positive integer and $k\Delta$ is a divisor. 
\end{rem}

$ $
\begin{rem}
As an alternative to Definition \ref{def:Gibbs polystability}, one
could define Gibbs polystability without assuming that $\mathrm{Aut}_{0}(X)$
is reductive. To that end, let $\mathbb{G}$ be an arbitrary reductive
subgroup of $\mathrm{Aut}_{0}(X)$, for example a Levi subgroup, i.e.
it is a maximal reductive subgroup. Define (uniform) $\G$-Gibbs polystability
of a Fano manifold $X$ as in Definition \ref{def:Gibbs polystability}
but with possibly $\G\neq\mathrm{Aut}_{0}(X)$. If $\G$ is a Levi
subgroup, the condition is independent of such a choice since any
two Levi subgroups are conjugate and $\mathcal{D}^{(N)}$ is $\mathrm{Aut}_{0}(X)$-invariant
(Lemma \ref{lem:character}). In the upcoming work \cite{a-b-b-j},
it is shown that uniform $\G$-Gibbs polystability of a Fano manifold
implies K-polystability. A reasonable extension of Conjecture \ref{conj:Gibbs poly intro}
in this setting is that uniform $\G$-Gibbs polystability for a Levi
subgroup $\G$ is equivalent to K-polystability. Motivated by the
upcoming work \cite{b-j-r-s}, showing that vanishing of the quantized
Futaki invariant for a single $k$ implies reductivity of the automorphism
group for a toric Fano variety, another problem is to show that Gibbs
polystability at a single level $k$ implies reductivity of the automorphism
group. 
\end{rem}

It follows directly from the definitions that 
\begin{equation}
\gamma_{\text{PS}}^{(N)}(X,\Delta)\geq\gamma^{(N)}(X,\Delta)\geq\alpha(X,\Delta)\coloneqq\inf_{D\sim-K_{(X,\Delta)}}\text{lct}\left(X,\Delta;D\right),\label{eq:lower bound in terms of alpha}
\end{equation}
where $D$ ranges over all effective $\Q$-divisors $D$ linearly
equivalent to $-K_{(X,\Delta)}$. The invariant $\alpha(X,\Delta)$
coincides with Tian's\emph{ alpha invariant} of $(X,\Delta),$ which
is strictly positive for any log Fano manifold $(X,\Delta).$

\subsection{Preliminaries}

To compute the microscopic stability thresholds $\gamma^{(N)}(X,\Delta)$
and $\gamma_{\text{PS}}^{(N)}(X,\Delta)$ on log Fano curves $(X,\Delta)$
we first recall some preliminary notions and results.

\subsubsection{\label{subsec:Log-canonical-thresholds}Log canonical thresholds}

Let $(X,\Delta)$ be a \emph{log pair,} defined here as a normal projective
variety $X,$ that we for simplicity assume is non-singular, together
with an effective $\mathbb{Q}$-divisor $\Delta$ on $X$ such that
$K_{X}+\Delta$ defines a $\mathbb{Q}$-line bundle. Moreover, we
will always fix a (multivalued) holomorphic section $s_{\Delta}$
cutting out $\Delta$. The log canonical threshold of a $\Q$-divisor
$D$ on $X$ with respect to $(X,\Delta)$ may, in analytic terms,
be defined as the following \emph{integrability threshold}:
\begin{equation}
\mathrm{lct}(X,\Delta;D)=\sup\left\{ c\colon||s_{D}||_{D}^{-2c}||s_{\Delta}||_{\Delta}^{-2}\in L_{\mathrm{loc}}^{1}(X)\right\} ,\label{eq:lct anal}
\end{equation}
 where $s_{D}$ is a (multivalued) holomorphic section cutting out
$D$ and where $||\cdot||_{D}$ and $||\cdot||_{\Delta}$ denotes
fixed smooth metrics on the line bundles corresponding to $D$ and
$\Delta$, respectively. In order to calculate $\gamma^{(N)}(X,\Delta)$
and $\gamma_{\text{PS}}^{(N)}(X,\Delta)$ we will instead use the
well-known \emph{valuative} definition of the log canonical threshold
(see, for example, \cite[Thm B.5]{bbj}): 
\begin{equation}
\mathrm{lct}(X,\Delta;D)=\inf_{v}\mathrm{lct}(X,\Delta;D)_{v},\,\,\,\,\mathrm{lct}(X,\Delta;D)_{v}\coloneqq\frac{A_{(X,\Delta)}(v)}{v(D)},\label{eq:lct val}
\end{equation}
where the inf ranges over all divisorial valuations $v$ over $X$
and $A_{(X,\Delta)}(v)$ denotes the \emph{log discrepancy} of $v$
with respect to $(X,\Delta).$ We recall that a \emph{divisorial valuation}
$v$ over $X$ is an integer-valued function on the space $\C(X)$
of all rational functions $f$ on $X$ of the form
\[
v(f)=\mathrm{ord}_{E}(\pi^{*}f)
\]
where $\pi\colon Y\rightarrow X$ is a birational morphism from a
non-singular variety $Y,$ $E$ is a prime divisor on $Y$ and $\mathrm{ord}_{E}(g)$
denotes the order of vanishing of a given rational function $g$ along
$E.$ The subvariety $\pi(E)$ is called the \emph{center }of the
divisorial valuation. The action of a valuation on $\C(X)$ naturally
extends to sections of line bundles (since the action is invariant
under scaling by a non-zero holomorphic function) and thus to divisors
$D.$ Moreover, by imposing homogeneity it also extends to $\mathbb{Q}$-divisors,
so that $v$ takes rational values\footnote{These divisorial valuations are usually called $\emph{normalized}$
divisorial valuations, with a general divisorial valuation being a
positive multiple of a normalized divisorial valuation, but here we
will stick to the above notation. }. The \emph{log discrepancy} $A_{(X,\Delta)}(v)$ of a divisorial
valuation $v$ is defined as 
\begin{equation}
A_{(X,\Delta)}(v)=1+\mathrm{ord}_{E}(K_{Y}-\pi^{*}K_{(X,\Delta)}),\label{eq:def of log disc}
\end{equation}
 where $K_{Y}-\pi^{*}K_{(X,\Delta)}$ defines a divisor, supported
on the exceptional divisors of $\pi.$ 
\begin{example}
\label{exa: divisorial valuation ass to a blowup}Any smooth subvariety
$Z$ of $X$ determines a divisorial valuation that we will refer
to as $v_{Z}$, defined by $v=\mathrm{ord}_{E}(\pi^{*}(\cdot))$,
where $E$ is the exceptional divisor of the blowup $\pi\colon\mathrm{Bl}_{Z}(X)\rightarrow X$.
The center of $v$ is simply $Z$ and $A_{X}(v)$ is the codimension
of $Z.$
\end{example}

Let $\pi\colon Y\rightarrow X$ be a \emph{log resolution} of $(X,D+\Delta),$
i.e. a birational morphism that is an isomorphism over the complement
of $\text{supp}(D+\Delta)$ such that $\pi^{*}(D+\Delta)$ is a divisor
on $Y$ with simple normal crossings. Denoting by $E_{1},\dots,E_{m}$
its components, the inf in \ref{eq:lct val} can be restricted to
the divisorial valuations defined by $E_{1},\dots,E_{m}$. Using the
existence of log resolutions it then follows readily that the two
definitions \ref{eq:lct anal} and \ref{eq:lct val} of the log canonical
threshold agree. 

We will use the following basic\emph{ product property}: if $(X_{1},\Delta_{1})$
and $(X_{2},\Delta_{2})$ are log pairs, and $D_{1},D_{2}$ are effective
$\mathbb{Q}$--divisors on $X_{1}$ and $X_{2}$ respectively, then
\begin{equation}
\mathrm{lct}(X_{1}\times X_{2},\mathrm{pr_{1}^{*}(\Delta_{1})+\mathrm{pr}_{2}^{*}(\Delta_{2});\mathrm{pr_{1}^{*}(D_{1})+\mathrm{pr}_{2}^{*}(D_{2}))=\min_{i=1,2}\mathrm{lct}(X_{i},\Delta_{i};D_{i})}},\label{eq: outer product lct}
\end{equation}

where $\mathrm{pr}_{j}\colon X_{1}\times X_{2}$ are the natural projections.
This follows directly from the analytic representation \ref{eq:lct anal},
by decomposing the density $||s_{D}||^{-2c}||s_{\Delta}||^{-2}$ as
a product of densities on $X_{1}$ and $X_{2}.$ We will refer to
$\mathrm{pr_{1}^{*}(D_{1})+pr_{2}^{*}(D_{2})}$ as $D_{1}\boxplus D_{2}$
in this context. 

Next, if $U$ is an open subset in the Euclidean topology of $X$,
then $\mathrm{lct}(U,\Delta|_{U},D|_{U})$ designates the inf \ref{eq:lct val}
restricted to valuations whose centers have non-zero intersection
with $U$. With this definition, analogous versions of the above properties
hold, and moreover, if $X=\cup_{i}U_{i}$, then clearly, 
\[
\mathrm{lct}(X,\Delta;D)=\inf_{i}\mathrm{lct}(U_{i},\Delta|_{U_{i}},D|_{U_{i}}).
\]
A divisorial valuation $v$ over $X$ is said to \emph{compute }$\mathrm{lct}(X,\Delta;D)$
if the inf \ref{eq:lct val} is attained at $v.$ More generally,
when considering a sequence $\mathrm{lct}(X_{N},\Delta_{N};D_{N})$,
a sequence $v_{N}$ of valuations over $X_{N}$ will be said to \emph{asymptotically
compute} $\mathrm{lct}(X_{N},\Delta_{N};D_{N})$ if $\text{lct}_{v_{N}}\left(X^{N},\Delta_{N};D_{N}\right)$
converges towards $\liminf_{N\rightarrow\infty}\text{lct}\left(X^{N},\Delta_{N};D_{N}\right),$
as $N\rightarrow\infty.$ Finally, note that any $N$-tuple of divisorial
valuations $v_{1},\dots,v_{N}$ over $X$ induces a valuation over
$X^{N}$ by blowing up the variety $E_{1}\times\cdots\times E_{N}$
over $X^{N},$ where $E_{i}$ is the divisor over $X$ corresponding
to $v_{i}.$ Such a divisorial valuation over $X^{N}$ will be called
\emph{decomposable}. 

\subsubsection{\label{subsec:Semi-stability-in-GIT}Semistability in GIT}

Let $L\rightarrow Y$ be an ample line bundle over a complex projective
manifold $Y$ and $G$ a complex reductive group acting on $(Y,L).$
The corresponding \emph{semistable locus }(in sense of Geometric Invariant
Theory) is the Zariski open subset $Y_{\text{ss}}$ of $Y$ consisting
of all $y\in Y$ such that there exists a $G$-invariant section $s_{k}\in H^{0}(Y,kL)$
for some positive integer $k$ such that $s_{k}(y)\neq0$ \cite{ho}.
More generally, given an ample line bundle $L\rightarrow Y$ with
an action of $G$ on $(X,mL)$ for some positive integer $m,$ the
semistable locus $Y_{\text{ss }}$ is defined wrt $(X,mL).$ When
$Y=(\P^{1})^{N}$ and $G=\mathrm{PGL}(2,\C),$ acting diagonally,
$(X^{N})_{\text{ss}}$ consists of all $(x_{1},\dots,x_{N})$ such
that the fraction $I/N$ of points $x_{i}$ coinciding satisfies $I/N\leq1/2$,
i.e. 
\begin{equation}
\left(\frac{1}{N}\sum_{i=1}^{N}\delta_{x_{i}}\right)\{x\}\leq1/2\label{eq:SS cond}
\end{equation}
 for all $x\in\P^{1}$ (see \cite[Section 3]{m-f-k}). When $Y=(\P^{1})^{N}$
and $G=\C^{*}$ the semistable locus $(X^{N})_{\text{ss}}$ consists
of all $(x_{1},\dots,x_{N})$ such that the inequality \ref{eq:SS cond}
holds for the two fixed points $x$ of the action of $\C^{*}$ on
$\P^{1}.$ Just like the case $G=\mathrm{PGL}(2,\C$) this follows
readily from the Hilbert--Mumford criterion for one-parameter subgroups
of $\C^{*}\subset G.$ 

\subsubsection{\label{subsec:FM-compactification}The Fulton--MacPherson compactification
of configuration spaces}

In \cite{FuMa} a canonical algebraic compactification $X^{[N]}$
of the configuration space $X^{N}\setminus\cup_{i\neq j}\{x_{i}=x_{j}\}$
of a non-singular algebraic variety $X$ was constructed. It has the
following useful properties:
\begin{enumerate}
\item $X^{[N]}$ is a non-singular variety, projective if $X$ is \cite[Theorems 1 and 2]{FuMa}.
\item There is a canonical morphism $\pi_{\mathrm{FM}}\colon X^{[N]}\rightarrow X^{N}$,
which restricts to an isomorphism $\pi_{\mathrm{FM}}\colon X^{[N]}\setminus Z\rightarrow X^{N}\setminus\cup_{i\neq j}\{x_{i}=x_{j}\},$
where $Z$ defines a simple normal crossings divisor \cite[Theorem 3]{FuMa}.
\item There is a bijection between the set $\mathcal{W}$ of all submanifolds
$W$ of $X^{N}$ of the form 
\[
W=\left\{ x_{i_{1}}=x_{i_{2}}=\cdots=x_{i_{r}}\right\} \subset X^{N}
\]
 and a certain set of prime divisors $E_{W}^{\mathrm{}}$ on $X^{[N]}$,
such that for every $W\in\mathcal{W}$,
\[
\pi_{\mathrm{FM}}^{*}(\mathcal{I}(W))=\mathcal{I}\left(\sum_{W'\subset W}E_{W'}\right)
\]
 where $\mathcal{I}(\cdot)$ denotes the corresponding ideal sheaves
and the sum ranges over all $W'$ in $\mathcal{W}$ contained in $W$
\cite[Theorem 3]{FuMa}.
\end{enumerate}
Note that when $\mathrm{dim(X)=1}$, the first two properties imply
that the map $\pi_{FM}$ is a log resolution of the log pair $(X^{N},\sum_{i<j}\{x_{i}=x_{j}\})$. 

Given $W\in\mathcal{W},$ we denote by $v_{W}^{\mathrm{FM}}$ the
divisorial valuation corresponding to $E_{W}$. We can also consider
the valuation $v_{\mathrm{W}}$ as in Example \ref{exa: divisorial valuation ass to a blowup}.
The following Lemma shows that these two valuations associated to
a $W\in\mathcal{\mathcal{W}}$ are the same. 
\begin{lem}
\label{lem:equiv.  of the two val.  ass.  w W}For any $W\in\mathcal{W},$
$v_{W}^{\mathrm{FM}}=v_{W}$.
\end{lem}

\begin{proof}
First note that, to show that two valuations $v_{1}$ and $v_{2}$
over $X,$ associated to divisors $E_{1}$ on $Y_{1}\xrightarrow{\pi_{1}}X$
and $E_{2}$ on $Y_{2}\xrightarrow{\pi_{2}}X$, are the same it is
enough to show that there exist a birational model $Z$ dominating
$Y_{1}$ and $Y_{2}$ such that the strict transforms of $E_{1}$
and $E_{2}$ coincide on $Z.$ Now, by the universal property of blowups,
there is a birational morphism $f\colon X^{[N]}\rightarrow\mathrm{{Bl}}_{W}(X{}^{N})$
such that the following diagram commutes: 
\[
\begin{array}{rcl}
X^{[N]} & \xrightarrow{f} & \mathrm{{Bl}}_{W}(X^{N})\\
\pi_{\mathrm{FM}}\searrow &  & \swarrow\pi_{W}\\
 & X^{N}
\end{array}
\]
where $\pi_{W}\colon\mathrm{Bl}_{W}(X^{N})\rightarrow X^{N}$ is the
blow-up morphism, whose exceptional divisor we denote by $E_{\pi_{W}}.$
It will thus be enough to show that $f_{*}^{-1}(E_{\pi_{W}})=E_{W},$
where $f_{*}^{-1}(E_{\pi_{W}})$ denotes the strict transform of $E_{\pi_{W}}$
under $f.$ To this end let us first observe that there exists $W'$
in $\mathcal{W}$ such that
\begin{equation}
f_{*}^{-1}(E_{\pi_{W}})=E_{W'},\,\,\,W'\subset W\label{eq:pf of lemma equiv val}
\end{equation}
Indeed, $\pi_{W}^{*}(\mathcal{I}(W))=\mathcal{I}(E_{\pi_{W}})$ and
thus, by the commutativity of the diagram, $f^{*}(\mathcal{I}(E_{\pi_{W}}))=f^{*}\circ\pi_{W}^{*}(\mathcal{I}(W))=\pi_{\mathrm{FM}}(\mathcal{I}(E_{\pi_{W}})),$
which, by property 3 of the Fulton--MacPherson compactification above,
coincides with the ideal defined by the divisor $\sum_{W'\subset W}E_{W'}.$
Since the $E_{W}$'s are prime this proves the identity \ref{eq:pf of lemma equiv val}.
It follows from \ref{eq:pf of lemma equiv val} that 
\[
\pi_{FM}(E_{W'})=W
\]
(since $\pi_{\mathrm{FM}}(E_{W'})=\pi_{W}\circ f(E_{W'})=\pi_{W}(W)).$
But this implies that $W'=W.$ Indeed, by property 3 of the Fulton--MacPherson
compactification above, $W=\pi_{\mathrm{FM}}(E_{W'})\subset\pi_{\mathrm{FM}}(\bigcup_{W''\subset W'}E_{W''})=W'$,
and by assumption $W'\subset W$. 
\end{proof}
When $\mathrm{dim}(X)=1$, the subset $\cup_{i<j}\{x_{i}=x_{j}\}$
defines a divisor on $X^{N}$ that we denote by $D_{N}$. The following
lemma describes the log discrepancy and action of the valuations $v_{W}$
on $D_{N}$.
\begin{lem}
\label{lem: pullback divisors for Fulton-MacPherson}Let $\pi_{\mathrm{FM}}\colon(\mathbb{P}^{1})^{[N]}\rightarrow(\mathbb{P}^{1})^{N}$
be the log resolution described above. For any $W\in\mathcal{W}$
we have that
\[
A_{(\mathbb{P}^{1})^{N}}(v_{W})=\mathrm{codim}W.
\]
Additionally,
\[
v_{W}(D_{N})=\frac{\mathrm{codim}W(\mathrm{codim}W+1)}{2}.
\]
\end{lem}

\begin{proof}
Using Lemma \ref{lem:equiv.  of the two val.  ass.  w W} we realize
$v_{W}$ as the divisorial valuation associated to the exceptional
divisor $E_{\pi_{W}}$ of $\mathrm{Bl}_{W}(X{}^{N})$ along $W$.
Then the first formula already appeared in Example \ref{exa: divisorial valuation ass to a blowup}.
For the second formula, we use that the vanishing order of $\pi_{W}^{*}f$
along $E_{\pi_{W}}$, for any $f\in\C(X^{N})$, is the same as that
of $f$ along $W$, and this vanishing order is given by the number
of distinct subvarieties in $\mathcal{W}$ containing $W$. If $W=\{x_{i_{1}}=\cdots=x_{i_{r}}\}$
then the number of such subvarieties is precisely $r(r-1)/2=\mathrm{codim}W(\mathrm{codim}W+1)/2$. 
\end{proof}

\subsection{\label{subsec:Computation-of-stability}The microscopic stability
thresholds of log Fano curves}

We now turn to the case when $(X,\Delta)$ is a log Fano curve. This
means that $X=\P^{1}$ and $\Delta$ is any divisor of $r$ points
$p_{1},\dots,p_{m}$ on $X$ with weights $w_{i}\in]0,1[$ satisfying
\[
V\coloneqq\deg(K_{X}+\Delta)=2-\sum_{i=1}^{m}w_{i}>0.
\]
We can thus identify $-kK_{(\mathbb{P}^{1},\Delta)}$ with $kV\mathcal{O}(1),$
which shows that $N=kV+1.$ Now identifying $\P^{1}$ with $\C\cup\{\infty\}$
and taking the basis $e_{1},\dots,e_{N}$ appearing in the definition
\ref{eq:slater determinant text} of $\det S^{(N)}$ to consist of
monomials $1,z,z^{2},\dots,z^{kV}$ we may, in the affine chart $\C,$
express
\begin{equation}
\det S^{(N)}=\prod_{i<j}^{N}(z_{i}-z_{j}),\,\,\,\,\,N=Vk+1,\label{eq: vandermonde identity-1}
\end{equation}
 using the classical Vandermonde identity (and in general any point
in $(\mathbb{P}^{1})^{N}$ has a Zariski-open neighbourhood of this
type, by symmetry).

\subsubsection{The case $\Delta=0$}

As a warm-up, we start by computing $\gamma^{(N)}(\mathbb{P}^{1}),$
thus reproving the Gibbs semistability of $\mathbb{P}^{1}$, first
shown in \cite{Fu}, by algebraic means, and then in \cite{berm12},
by analytic means. 
\begin{prop}
\label{prop: alg proof of Gibbs semistab of P1}$\mathbb{P}^{1}$
is Gibbs semistable. More precisely, for any integer $N\geq2$
\[
\gamma^{(N)}(\mathbb{P}^{1})\coloneqq\mathrm{lct}((\mathbb{P}^{1})^{N};\mathcal{D}^{(N)})=1-1/N
\]
 and $\gamma^{(N)}(\mathbb{P}^{1})$ is computed by the divisorial
valuation induced by the diagonal of $X^{N}$ (viewed as a one-dimensional
submanifold). Moreover, for \emph{any} given point $p$ on $\P^{1}$
the sequence of decomposable valuations $v_{N}$ over $(\P^{1})^{N}$
induced by $(p,p,\dots,p)\in(\P^{1})^{N}$ asymptotically computes
$\mathrm{lct}((\mathbb{P}^{1})^{N};\mathcal{D}^{(N)}).$ 
\end{prop}

\begin{proof}
We have, 
\[
\mathrm{lct}((\mathbb{P}^{1})^{N};\mathcal{D}^{(N)})=\min_{v}\frac{A_{(\mathbb{P}^{1})^{N}}(v)}{v(\mathcal{D}^{(N)})}=\min_{v=v_{W}\colon W\in\mathcal{W}}\frac{A_{(\mathbb{P}^{1})^{N}}(v)}{v(\mathcal{D}^{(N)})}
\]
where the first minimum ranges over divisorial valuations on $X$,
which we can restrict to the ones of the form $v_{W},W\in\mathcal{W}$
since $\pi_{\mathrm{FM}}\colon(\mathbb{P}^{1})^{[N]}\rightarrow(\mathbb{P}^{1})^{N}$
is a log resolution of $((\mathbb{P}^{1})^{N},\mathcal{D}^{(N)})$,
and due to Lemma \ref{lem:equiv.  of the two val.  ass.  w W}. Lemma
\ref{lem: pullback divisors for Fulton-MacPherson} now yields
\begin{align*}
\mathrm{lct}((\mathbb{P}^{1})^{N},\mathcal{D}^{(N)}) & =k\min_{W\in\mathcal{W}}\frac{\mathrm{codim}W}{\mathrm{codim}W(\mathrm{codim}W+1)/2}=\min_{W\in\mathcal{W}}\frac{2k}{\mathrm{codim}W+1}\\
 & =\frac{2k}{2k+1}=\frac{N-1}{N}.
\end{align*}
For the second claim, note that the minimum is computed by any $W'$
such that $\mathrm{codim}(W')=N-1$, implying that $W'=(p_{1}=\cdots=p_{N})$.
Finally, for the valuation $v_{N}$ induced by $(p,p,\dots,p)\in(\P^{1})^{N}$,
since by Example \ref{exa: divisorial valuation ass to a blowup}
$A_{(\mathbb{P}^{1})^{N}}(v_{N})=N$ and $v_{N}(z_{i}-z_{j})=1,$
\[
\frac{A_{(\mathbb{P}^{1})^{N}}(v_{N})}{v_{N}(\mathcal{D}^{(N)})}=\frac{N}{N(N-1)/2k}=1,\,\,\,\,\,\,(2k=N-1).
\]
\end{proof}
In a similar way we next prove the uniform Gibbs polystability of
$\mathbb{P}^{1}$: 
\begin{prop}
\label{prop: alg proof of Gibbs polystab of P1}$\mathbb{P}^{1}$
is uniformly Gibbs polystable. More precisely, for an integer $N>3$,
\[
\gamma_{\mathrm{PS}}^{(N)}(\mathbb{P}^{1}\mathrm{)}=2,\text{ if }N\text{ is odd,\,\,\,\,}\gamma_{\mathrm{PS}}^{(N)}(\mathbb{P}^{1})=2(1-1/N)\text{ if }N\text{ is even}.
\]

For $N=2,3$, $\gamma_{\mathrm{PS}}^{(N)}(\mathbb{P}^{1}\mathrm{)}=\infty$. 
\begin{itemize}
\item In case $N>3$, $\gamma_{\mathrm{PS}}^{(N)}(\mathbb{P}\mathrm{^{1})}$
is computed by the divisorial valuations over $(X^{N})_{\mathrm{ss}}$
induced by the submanifolds $\{x_{i_{1}}=\cdots=x_{i_{\lfloor N/2\rfloor}}\}\subset X^{N}$
for any given index $I$ of length $\lfloor N/2\rfloor.$ 
\item A decomposable sequence $v_{N}$ of divisorial valuations over $(X^{N})_{\mathrm{ss}}$
asymptotically computing $\gamma_{\mathrm{PS}}^{(N)}(\mathbb{P}\mathrm{^{1})}$
is obtained by fixing three distinct points $p_{1},p_{2}$ and $p$
on $\P^{1}$ and letting $v_{N}$ be the valuation induced by $\{p_{1}\}^{\lfloor N/2\rfloor}\times\{p\}\times\{p_{2}\}^{\lfloor N/2\rfloor}\in(\P^{1})^{N}$
when $N$ is odd and by $\{p_{1}\}^{N/2}\times\{p_{2}\}^{N/2}\in(\P^{1})^{N}$
when $N$ is even.
\end{itemize}
\end{prop}

\begin{proof}
By Proposition \ref{prop: alg proof of Gibbs semistab of P1}, $\mathbb{P}^{1}$
is Gibbs semistable. The map $\pi_{FM}\colon Y\rightarrow(\mathbb{P}^{1})_{\mathrm{ss}}^{N}$
given by restricting $\pi_{\mathrm{FM}}$ to $Y\coloneqq\pi_{\mathrm{FM}}^{-1}((\mathbb{P}^{1})_{\mathrm{ss}}^{N})$,
is a log resolution of $((\mathbb{P}^{1})_{\mathrm{ss}}^{N},\mathcal{D}^{(N)})$.
Denote by $\mathcal{W}^{\mathrm{ss}}$ the subvarieties in $\mathcal{W}$
that have non-empty intersection with the semistable locus. By the
explicit description in Section \ref{subsec:Semi-stability-in-GIT},
$\mathcal{W}^{\mathrm{ss}}$ are the elements $W$ of $\mathcal{W}$
such that $\mathrm{codim}W\leq\lfloor N/2\rfloor-1$, which is non-empty
iff $N>3$. Consequently $\mathrm{lct}((\mathbb{P}^{1})_{\mathrm{ss}}^{N},\mathcal{D}^{(N)})=\infty$
if $N=2,3$ and otherwise
\[
\mathrm{lct}((\mathbb{P}^{1})_{\mathrm{ss}}^{N},\mathcal{D}^{(N)})=k\min_{W\in\mathcal{W}^{\mathrm{ss}}}\frac{\mathrm{codim}W}{\mathrm{codim}W(\mathrm{codim}W+1)/2}=\min_{W\in\mathcal{W}^{\mathrm{ss}}}\frac{2k}{\mathrm{codim}W+1}=\frac{N-1}{\lfloor N/2\rfloor}
\]
since the center of $v_{W}$ is precisely $W$. Finally, when $N$
is even, consider the valuation $v_{N}$ induced by $\{p_{1}\}^{N/2}\times\{p_{2}\}^{N/2}\in(\P^{1})^{N}$.
By Example \ref{exa: divisorial valuation ass to a blowup} $A(v_{N})=N$
and decomposing $\prod_{i<j}^{N}(z_{i}-z_{j})$ into the three groups
of indices, $(i,j)$ where $i,j\leq N/2,$ $i,j>N/2$, and the rest,
reveals that $kv_{N}(\mathcal{D}^{(N)})=2\left((N/2)\left((N/2)-1\right)/2\right).$
Hence, 
\[
\frac{A(v_{N})}{v_{N}(\mathcal{D}^{(N)})}=\frac{N}{(N/2)\left((N/2)-1\right)/k}=\frac{N(N-1)/2}{(N/2)\left((N/2)-1\right)}=2\frac{N-1}{N-2}
\]
converging to $2$ as $N\rightarrow\infty.$ The computation in the
case when $N$ is odd is similar.
\end{proof}

\subsubsection{The case $\Delta\protect\neq0$}

Instead of constructing an explicit log resolution of $((\mathbb{P}^{1})^{N},\mathcal{D}^{(N)}+\Delta_{N})$,
generalizing the Fulton--MacPherson compactification of $(\mathbb{P}^{1})^{N}$
(which, presumably, can be done using the results in \cite{D-CiPr}
or \cite{lili}) we will pursue a different route that reduces the
calculations to an application of the Fulton--MacPherson compactification
of $(\mathbb{P}^{1})^{N+1}.$ We begin by computing the invariant
$\gamma^{(N)}(X,\Delta).$ 
\begin{prop}
\label{prop:gamma N for log Fano curve}For any log Fano curve $(\mathbb{P}^{1},\Delta)$
and integer $N\ge2$,
\begin{equation}
\gamma^{(N)}(X,\Delta)=\mathrm{lct}(X^{N},\Delta_{N};\mathcal{D}^{(N)})=\frac{2}{V}\min\left\{ 1-1/N,\min_{l}\{1-w_{l}\}\right\} .\label{eq:detS lct log pair-1}
\end{equation}
Moreover, if $\frac{2(N-1)}{VN}\leq\min_{l}\frac{2(1-w_{l})}{V},$
then $v_{\{x_{1}=\cdots=x_{N}\}}$ computes $\mathrm{lct}(X^{N},\Delta_{N};\mathcal{D}^{(N)})$
and if $\frac{2(N-1)}{VN}\geq\min_{l}\frac{2(1-w_{l})}{V}$, then
$\mathrm{lct}(X^{N},\Delta_{N};\mathcal{D}^{(N)})$ is computed by
any divisorial valuation of the form $v_{\{x_{1}=\cdots=x_{N}=p_{l}\}}$,
where $l$ realizes the minimum in $\min_{l}\frac{2(1-w_{l})}{V}$
and recall that $p_{l}$ is in the support of $\Delta$. 
\end{prop}

\begin{proof}
Given $w\in]0,1[$ and $p\in\P^{1}$ denote by $\Delta_{N}^{p,w}$
the divisor on $(\P^{1})^{N},$ invariant under permutations, that
projects to the divisor $wp$ under any of the $N$ natural projections
to $X.$ 

\emph{Step 1: The following formula holds: }
\begin{equation}
\mathrm{lct}((\mathbb{P}^{1})^{N},\Delta_{N};\mathcal{D}^{(N)})=\min_{i}\mathrm{lct}((\mathbb{P}^{1})^{N},\Delta_{N}^{p_{i},w_{i}},\mathcal{D}^{(N)}).\label{eq:pf of Prop gamma N log Fano}
\end{equation}
 To prove this, consider the following covering of $(\mathbb{P}^{1})^{N}$
by open subsets, 
\begin{equation}
(\mathbb{P}^{1})^{N}=U\cup\bigcup_{q\in\{p_{1},\dots,p_{m}\}^{\tilde{N}},\tilde{N}=1,\dots,N}U_{q},\label{eq: log p1,  partition}
\end{equation}
where $U=\{(x_{1},\dots,x_{N})\in(\mathbb{P}^{1})^{N}\colon x_{i}\neq p_{j}\forall i,j\}$
and $U_{q}=\{(x_{1},\dots,x_{N})\in(\mathbb{P}^{1})^{N}\colon d(x_{i},q_{i})<d_{0}\forall i\}$
where $d$ is the round metric on $\mathbb{P}^{1}$ and $d_{0}<\min_{i,j}d(p_{i},p_{j})/2$.
Then 
\begin{equation}
\mathrm{lct}((\mathbb{P}^{1})^{N},\Delta_{N};\mathcal{D}^{(N)})=\min\Big\{\mathrm{lct}(U,\Delta_{N}|_{U};\mathcal{D}^{(N)}|_{U}),\min_{q\in\{p_{1},\dots,p_{m}\}^{\tilde{N}},\tilde{N}=1,\dots,N}\mathrm{\{lct}(U_{q},\Delta_{N}|_{U_{q}};\mathcal{D}^{(N)}|_{U_{q}})\}\Big\}.\label{eq: log p1, dividing lct by subsets}
\end{equation}

Since $\Delta_{N}|_{U}$ is trivial, we can use the Fulton--MacPherson
compactification as a log resolution on $U$, and we turn to the other
cases. Without loss of generality we may assume that $q=(p_{i_{1}},\dots,p_{i_{1}},p_{i_{2}},\dots,p_{i_{2}},\dots,p_{i_{l}},\dots,p_{i_{l}})$
for $1\leq l\leq m$ and some $p_{i_{1}},\dots,p_{i_{l}}$ distinct
and where $p_{i_{1}}$ is repeated $k_{1}$ times and $p_{i_{2}}$
is repeated $k_{2}$ times, and so on. Then $U_{q}=V_{1}\times\cdots\times V_{l}\times(\mathbb{P}^{1})^{r}$,
for some $0\leq r\leq N-1$, where $V_{j}=\{(x_{1},\dots,x_{k_{j}})\in(\mathbb{P}^{1})^{k_{j}}\colon d(x_{l},p_{i_{j}})<d_{0}\forall l\}$.
Observe that $\Delta_{N}|_{U_{q}}=\Delta_{k_{1}}^{p_{i_{1}},w_{i_{1}}}\boxplus\cdots\boxplus\Delta_{k_{l}}^{p_{i_{l}},w_{i_{l}}}\boxplus0$.
Furthermore, $\mathcal{D}^{(N)}|_{U_{q}}\eqqcolon\mathcal{D}_{k_{1}}\boxplus\cdots\boxplus\mathcal{D}_{k_{l}}\times\mathcal{D}_{r}$.
By \ref{eq: outer product lct}
\[
\mathrm{lct}(U_{q},\Delta_{N}|_{U_{q}};\mathcal{D}^{(N)}|_{U_{q}})=\min\Bigg\{\min_{1\leq j\leq l}\Bigg\{\mathrm{lct}\left(\left.V_{j},\Delta_{k_{i}}^{p_{i_{j}},w_{i_{j}}}\right|_{V_{j}};\left.\mathcal{D}_{k_{j}}\right|_{V_{j}}\right)\Bigg\},\mathrm{lct}\left((\mathbb{P}^{1})^{r};\mathcal{D}_{r}\right)\Bigg\}.
\]
But for any $1\leq j\leq l$, $\pi^{*}(\Delta_{k_{i}}^{p_{i_{j}},w_{i_{j}}})\leq\Delta_{N}^{p_{i_{j}},w_{i_{j}}}$,
and $\pi^{*}(\mathcal{D}_{k_{j}})\leq\mathcal{D}^{(N)}$ where $\pi\colon(\mathbb{P}^{1})^{N}\rightarrow(\mathbb{P}^{1})^{k_{j}}$
is the projection and similarly $\pi^{*}(\mathcal{D}_{r})\leq\mathcal{D}^{(N)}$.
In view Definition \ref{eq:lct val} and the formula \ref{eq: outer product lct},
$\mathrm{lct}(V_{j},\Delta_{k_{i}}^{p_{i_{j}},w_{i_{j}}};\mathcal{D}_{k_{j}})\geq\mathrm{lct}(U_{q'},\Delta_{N}|_{U_{q'}};\mathcal{D}^{(N)}|_{U_{q'}})$
with $q'=(p_{i_{j}},\dots,p_{i_{j}})$ but the numbers $\mathrm{lct}(U_{q'},\Delta_{N}|_{U_{q'}};\mathcal{D}^{(N)}|_{U_{q'}})$
are contenders in the inf \ref{eq: log p1,  partition}. Moreover,
$\mathrm{lct}(U_{q'},\Delta_{N}|_{U_{q'}};\mathcal{D}^{(N)}|_{U_{q'}})=\mathrm{lct}((\mathbb{P}^{1})^{N},\Delta_{N}^{p_{i_{j}},w_{i_{j}}};\mathcal{D}^{(N)})$
and any such term is lesser than or equal to $\mathrm{lct}(U,\Delta_{N}|_{U};\mathcal{D}^{(N)}|_{U})$
, and also $\mathrm{lct}((\mathbb{P}^{1})^{r};\mathcal{D}_{r})$.
This proves formula  \ref{eq:pf of Prop gamma N log Fano}.

\emph{Step 2: The following formula holds }
\[
\mathrm{lct}((\mathbb{P}^{1})^{N},\Delta_{N}^{p,w};\mathcal{D}^{(N)})=\frac{2(1-w)}{V}.
\]
To prove this, consider the affine chart $\C\subset\P^{1}$ centered
at $p$ and the restrictions to $\mathbb{C}^{N}\subset(\mathbb{P}^{1})^{N}$
of the (multivalued) sections cutting out $\mathcal{D}^{(N)}$ and
$\Delta_{N}^{p,w},$ respectively: 
\begin{equation}
\left(\prod_{i<j}^{N}(z_{i}-z_{j})\right)^{-1/k}\text{and\,\,}\left(\prod_{i}^{N}z_{i}\right)^{w}.\label{eq: support divisor equation in log case}
\end{equation}
Now regard these as functions of $N+1$ variables $z_{1},\dots,z_{N}$
and $z_{N+1}$, which are independent of $z_{N+1}.$ In other words,
we consider their pullbacks to $\C^{N+1}$ under the natural projection
from $\C^{N+1}$ to $\C^{N}.$ By the product property \ref{eq: outer product lct},
this does not affect the computation of the corresponding log canonical
threshold. After making the affine change of variables $z_{i}=u_{i}-u_{N+1}$
for $i=1,\dots,N$ and $z_{N+1}=u_{N+1}$ the two expressions \ref{eq: support divisor equation in log case}
become 
\begin{equation}
\left(\prod_{1\leq i<j\leq N}(u_{i}-u_{j})\right)^{-1/k}\text{and\,\,}\left(\prod_{1\leq i\leq N}(u_{i}-u_{N+1})\right)^{w}.\label{eq: braid arrangement in one more dimension}
\end{equation}
We can thus use the Fulton--MacPherson compactification in the $N+1$
variables $u_{1},\dots,u_{N+1},$ $\pi_{FM}\colon\mathbb{P}^{[N+1]}\rightarrow\mathbb{P}^{N+1},$
as a log resolution of the pullback of $\mathcal{D}^{(N)}+\Delta_{N}^{p,w}.$
However, now there are two types of subvarieties $W$ in $\mathcal{W}$
that we need to treat differently when computing $\mathrm{lct}((\mathbb{P}^{1})^{N},\Delta_{N}^{p,w};\mathcal{D}^{(N)})$:
those not involving $u_{N+1},$ that we denote by $\mathcal{W}',$
and those involving $u_{N+1},$ that we denote by $\mathcal{W}_{p}.$
Denoting by $v_{W}$ the corresponding divisorial valuations over
$\mathbb{P}^{N+1}$ we have
\[
v_{W}(\Delta^{p,w})=0,\,\text{if \ensuremath{W\in\mathcal{W}'}},\,\,\,\,\,\,v_{W}(\Delta^{p,w})=w\mathrm{codim}W\,\text{if \,W\ensuremath{\in\mathcal{W}_{p}},}
\]
as follows exactly as in Lemma \ref{lem: pullback divisors for Fulton-MacPherson},
using Lemma \ref{lem:equiv.  of the two val.  ass.  w W}. From the
definition it follows that

\[
A_{((\mathbb{P}^{1})^{N},\Delta^{p,w})}(v_{W})=A_{(\mathbb{P}^{1})^{N}}(v_{W})-v_{W}(\Delta^{p,w})
\]
and consequently, 
\begin{align*}
\mathrm{lct}((\mathbb{P}^{1})^{N},\Delta_{N}^{p,w};\mathcal{D}^{(N)}) & =\min\left(\min_{W\in\mathcal{W}'}\frac{A_{((\mathbb{P}^{1})^{N+1},\Delta_{N}^{p,w})}(v_{W})}{v_{W}(\mathcal{D}^{(N)})},\min_{W\in\mathcal{W}_{p}}\frac{A_{((\mathbb{P}^{1})^{N+1},\Delta_{N}^{p,w})}(v_{W})}{v_{W}(\mathcal{D}^{(N)})}\right)\\
 & =\min\left(\min_{W\in\mathcal{W}'}\frac{2k}{\mathrm{codim}W+1},\min_{W\in\mathcal{W}^{p}}\frac{\mathrm{codim}W-w\mathrm{codim}W}{\frac{1}{k}\mathrm{codim}W(\mathrm{codim}W-1)/2}\right).
\end{align*}
Since
\[
\min_{W\in\mathcal{W}'}\frac{2k}{\mathrm{codim}W+1}=\frac{2k}{N}=\frac{2(N-1)}{VN},
\]
we thus get 
\[
\min_{W\in\mathcal{W}^{p}}\frac{\mathrm{codim}W-w_{l}\mathrm{codim}W}{\frac{1}{k}\mathrm{codim}W(\mathrm{codim}W-1)/2}=\min_{W\in\mathcal{W}^{p}}\frac{2k(1-w)}{\mathrm{codim}W-1}=\frac{2k(1-w)}{N-1}=\frac{2(1-w)}{V},
\]
which concludes the proof of the formula in Step 2. 

Finally, letting $\pi$ be the blowup along the subvariety $W'\coloneqq\{x_{1}=\cdots=x_{N}\}$
and $\pi_{W'}$ the corresponding divisorial valuation, we find that
\[
\frac{A_{((\mathbb{P}^{1})^{N},\Delta_{N})}(v_{W'})}{v_{W'}(\mathcal{D}^{(N)})}=\frac{2k}{N}=\frac{2(N-1)}{VN}
\]
which proves the claim when $\frac{2(N-1)}{VN}\leq\min_{l}\frac{2(1-w_{l})}{V}.$
The other case follows similarly by blowing up along a subvariety
of the form $\{x_{1}=\dots=x_{N}=p_{l}\}$.
\end{proof}
The formula for $\gamma^{(N)}(X,\Delta)$ in the previous proposition
corrects the formula given in \cite[Thm 4.5]{berm12}, which is only
correct up to an error term of order $O(1/N)$ (that comes from the
discrepancy between the upper and the lower bound on $\gamma^{(N)}(X,\Delta)$
established by analytic means in \cite{berm12}). 

\subsubsection{The case of nontrivial $\G$}

To investigate which log Fano curves are Gibbs polystable, but not
Gibbs stable, we need to restrict to the cases with nontrivial $\G.$
This means the number $m$ of components of $\Delta$ satisfies $m\leq2.$
It is also required that $(X,\Delta)$ be Gibbs semistable. This implies
that there is some point $p_{i}$ in the support of $\Delta$ such
that $(1-w_{i})=V/2,$ by the formula for $\gamma(X,\Delta)$ resulting
from Proposition \ref{prop:gamma N for log Fano curve}. If $m=1,$
then this equality implies that $V=0$ in which case $(X,\Delta)$
is not log Fano. We can thus restrict to $m=2$, since $m=0$ is the
content of Proposition \ref{prop: alg proof of Gibbs polystab of P1}.
In this case $w_{i}=w$ for some $w\in]0,1[,$ by Proposition \ref{prop:gamma N for log Fano curve}.
Such a divisor $\Delta$ on $\P^{1}$ will be denoted by $\Delta_{w}.$ 
\begin{prop}
\label{prop: Gibbs polystability of P1 with two pts}Let $\Delta_{w}$
be the divisor on $\P^{1}$ supported on two distinct points $p_{1}$
and $p_{2}$ with the same weight $w\in]0,1[.$ $(\mathbb{P}^{1},\Delta_{w})$
is uniformly Gibbs polystable. More precisely, for $N=2,3$, $\gamma^{(N)}(X,\Delta_{w})^{\G}=\frac{2(1-1/N)}{V}$
and for any $N\geq3$,
\[
\gamma_{\mathrm{PS}}^{(N)}(X,\Delta_{w})\coloneqq\mathrm{lct}\left((\mathbb{P}^{1})_{\mathrm{ss}}^{N},\Delta_{N};\mathcal{D}^{(N)}\right)=\min\Bigg\{\frac{2(1-1/N)}{V},\frac{N-1}{\lfloor N/2\rfloor-1}\Bigg\}
\]
\begin{itemize}
\item Now assume $N>3$. When $V\geq2(\lfloor N/2\rfloor-1)/N$, $\gamma_{\mathrm{PS}}^{(N)}(X,\Delta_{w})$
is computed by the valuation induced by the diagonal in $X^{N}$ and
when $V\leq2(\lfloor N/2\rfloor-1)/N$, $\gamma_{\mathrm{PS}}^{(N)}(X,\Delta_{w})$
is computed by the divisorial valuations defined by the submanifolds
in $X^{N}$ given by $\{x_{i_{1}}=\cdots=x_{i_{\lfloor N/2\rfloor}}=p_{1}\}$
and $\{x_{i_{1}}=\cdots=x_{i_{\lfloor N/2\rfloor}}=p_{2}\}$ for any
given index $I$ of length $\lfloor N/2\rfloor$.
\item When $V\geq1$, $\gamma_{\mathrm{PS}}^{(N)}(X,\Delta_{w})$ is asymptotically
computed by the decomposable valuation $v_{N}$ induced by $(p,\dots,p)\in(\P^{1})^{N}$
for any given point $p$ not in the support of $\Delta_{w}$ and when
$V\leq1$, $\gamma_{\mathrm{PS}}^{(N)}(X,\Delta_{w})_{\mathrm{}}$
is asymptotically computed by the sequence $v_{N}$ appearing in Prop
\ref{prop: Gibbs polystability of P1 with two pts}, if $p_{1}$ and
$p_{2}$ are the points in the support of $\Delta_{w}$.
\end{itemize}
\end{prop}

\begin{proof}
In this case the group $\G$ is isomorphic to $\C^{*}.$ Denote by
$\mathcal{W}_{p}^{\mathrm{ss}}$ the elements of $\mathcal{W}_{p}$
from Lemma \ref{prop:gamma N for log Fano curve} that have a non-empty
intersection with the corresponding semistable locus on $X^{N}$ (explicitly
described in Section \ref{subsec:Semi-stability-in-GIT}). Following
closely the proof of Proposition \ref{prop:gamma N for log Fano curve},
but taking into account that we restrict to the semistable locus,
as in \ref{prop: alg proof of Gibbs polystab of P1}, we find that
\[
\mathrm{lct}((\left.(\mathbb{P}^{1})_{\mathrm{ss}}^{N})\right|_{U_{p}},\Delta_{N};\mathcal{D}^{(N)})=\min\Bigg\{\min_{W\in\mathcal{W}'}\Bigg\{\frac{2k}{\mathrm{codim}W+1}\Bigg\},\min_{W\in\mathcal{W}_{p}^{\mathrm{ss}}}\Bigg\{\frac{\mathrm{codim}W-w\mathrm{codim}W}{\frac{1}{k}\mathrm{codim}W(\mathrm{codim}W-1)/2}\Bigg\}\Bigg\}
\]
when $N>3$. We have, as in Proposition \ref{prop:gamma N for log Fano curve},
\[
\min_{W\in\mathcal{W}'}\Bigg\{\frac{2k}{\mathrm{codim}W+1}\Bigg\}=\frac{2(N-1)}{VN},
\]
while 
\[
\min_{W\in\mathcal{W}_{p}^{\mathrm{ss}}}\Bigg\{\frac{\mathrm{codim}W-w\mathrm{codim}W}{\frac{1}{k}\mathrm{codim}W(\mathrm{codim}W-1)/2}\Bigg\}=\min_{W\in\mathcal{W}_{p}^{\mathrm{ss}}}\Bigg\{\frac{2k(1-w)}{\mathrm{codim}W-1}\Bigg\}=\frac{2(N-1)(1-w)}{(\lfloor N/2\rfloor-1)V}=\frac{N-1}{(\lfloor N/2\rfloor-1)}
\]
The cases $N=2,3$ are similar except that $\mathcal{W}_{p}^{\mathrm{ss}}$
is empty. The last part, concerning divisorial valuations that compute
the log canonical threshold, follows similarly as in Proposition \ref{prop:gamma N for log Fano curve}. 
\end{proof}

\subsection{Gibbs polystability of log Fano curves}

We recall \cite{li000} that a log Fano curve $(X,\Delta)$ is \emph{K-(semi)stable}
iff $V/2>(\geq)1-w_{l}$ for all $l$ (and K-stable iff it is uniformly
K-stable). Moreover, when $\G$ is nontrivial, $(X,\Delta)$ is a
\emph{K-polystable} iff the Futaki character of $(X,\Delta)$ vanishes
iff $\Delta=0$ or $\Delta=\Delta_{w}$ for $w\in]0,1[$ (as in Prop
\ref{prop: Gibbs polystability of P1 with two pts}). Comparing with
Propositions \ref{prop: alg proof of Gibbs semistab of P1}, \ref{prop: alg proof of Gibbs polystab of P1},
\ref{prop:gamma N for log Fano curve} and \ref{prop: Gibbs polystability of P1 with two pts}
we thus find: 
\begin{thm}
\label{thm:gibbs  poly iff K-poly text}Let $(X,\Delta)$ be a log
Fano curve. Then $(X,\Delta)$ is Gibbs semistable if and only if
$(X,\Delta)$ K-semistable. Moreover, $(X,\Delta)$ is (uniformly)
Gibbs {[}poly{]}stable if and only if $(X,\Delta)$ is (uniformly)
K-{[}poly{]}stable. 
\end{thm}

This proves Conjecture \ref{conj:Gibbs poly intro} for log Fano curves.
It was previously known that Gibbs semistability was equivalent to
K-semistability and that (uniform) Gibbs stability was equivalent
to (uniform) K-stability (see \cite{berm12}).

\subsection{\label{subsec:Higher-dimensional-examples}Higher dimensional examples}

\subsubsection{\label{subsec:Gibbs-stable-orbifolds}Gibbs stable orbifolds that
are not exceptional}

Following \cite{c-s,bi}, a log Fano variety $(X,\Delta)$ is said
to be \emph{exceptional} if $\text{lct}(X,\Delta;D)>1$ for any $\Q$-divisor
$D$ linearly equivalent to $-K_{(X,\Delta)}.$ By \cite[Thm 1.7]{bi},
this is equivalent to $\alpha(X,\Delta)>1,$ where $\alpha(X,\Delta)$
denotes Tian's alpha invariant, whose algebro-geometric definition
appears in formula \ref{eq:lower bound in terms of alpha}. If $(X,\Delta)$
is exceptional, then $(X,\Delta)$ is Gibbs stable at any level $k.$
In fact, $(X,\Delta)$ is even uniformly Gibbs stable, by the lower
bound \ref{eq:lower bound in terms of alpha}.

Consider now the log pair $(X,\Delta)=(\P^{n},\Delta_{w}),$ where
\[
\Delta_{w}\coloneqq w(H_{1}+\cdots+H_{n+2}),\,\,\,w\in[0,1[
\]
 for $(n+2)$ given hyperplanes $H_{i}$ in $\P^{n}$ with normal
crossings. The log pair $(\P^{n},\Delta_{w})$ is Fano iff $w<1-1/(n+2)$
and exceptional iff $w>1-1/(n+1);$ a simple proof is provided in
Remark \ref{rem:except} below. 
\begin{example}
When $w=1-1/m$ for $m\in\Z_{+}$ the corresponding orbifold $(\P^{n},\Delta_{w})$
is the global orbifold quotient of the Fermat hypersurface $X_{m}$
of degree $m$ in $\P^{n+1}$ under the standard action of the group
$(\Z/m\Z)^{n}$ on $\P^{n+1},$ where $X_{m}$ is defined as the zero-locus
of $x_{0}^{m}+\cdots+x_{n+1}^{m}$ in $\P^{n+1}.$ Indeed, $(\Z/m\Z)^{n}$
is the Galois group for the Galois branched cover $X_{m}\rightarrow\P^{n}$
defined by $(x_{0}:\cdots:x_{n+1})\mapsto(x_{0}^{m}:\cdots:x_{n}^{m}).$
\end{example}

The following result provides the first examples of Gibbs stable log
Fano varieties --- beyond complex log Fano curves --- that are not
exceptional:
\begin{prop}
\label{prop:non-excep text}The log Fano variety $(\P^{n},\Delta_{1-1/(n+1)}),$
which may be identified with the Fano orbifold $X_{n+1}/(\Z/(n+1)\Z)^{n},$
is Gibbs stable at any level $k.$
\end{prop}

Although the group $\text{Aut}_{0}(\P^{n},\Delta)$ is trivial, the
proof of the previous proposition exploits that $\text{Aut}_{0}(\P^{n})$
is nontrivial, by applying the following lemma to $X=\P^{n}$:
\begin{lem}
Let $(X,D)$ be a projective log pair and $G$ a connected algebraic
group acting algebraically on $X.$ If $D$ is $G$-invariant, then
so is the center $Z_{v}$ of any divisorial valuation $v$ over $X$
computing $\mathrm{lct}(X,D).$ In particular, given a Fano variety
$X$ and a positive rational number $k$ such that $-kK_{X}$ is defined
as line bundle, any divisorial valuation computing the log canonical
threshold of the $\Q$-divisor $\mathcal{D}^{(N_{k})}$ on $X_{N_{k}}$
is invariant under the diagonal action of $\mathrm{Aut}_{0}(X).$
\end{lem}

\begin{proof}
Consider the action of $G$ on valuations $v$ over $X$ defined by
\[
(g\cdot v)(f)=v((g^{-1})^{*}f)
\]
 for any $g\in G$ and rational function $f$ on $X.$ If $v$ is
a divisorial valuation, then so is $g\cdot v.$ Indeed, if $v$ is
represented by a prime divisor $E$ on a model $\pi\colon Y\rightarrow X,$
then $g\cdot v$ is represented by the same divisor $E$ wrt the model
$g^{*}\pi\colon Y\rightarrow X.$ In particular, $A(v)=A(g\cdot v)$
and $Z_{gv}=g(Z_{v}).$ It follows that if $D$ is $G$-invariant
and $v$ computes $\mathrm{lct}(X,D),$ then $g\cdot v$ also computes
$\mathrm{lct}(X,D)$ and, as a consequence, $g(Z_{v})$ is the center
of a valuation computing $\mathrm{lct}(X,D),$ for any $g\in G.$
But, for a general log pair $(X,D),$ the number of centers of valuations
computing $\mathrm{lct}(X,D)$ is finite. Indeed, after scaling $D$
we may as well assume that $\mathrm{lct}(X,D)=1.$ Then, by definition,
$v$ computes $\mathrm{lct}(X,D)$ iff $v$ is a\emph{ log canonical
place} (in the terminology of the MMP) and it is well-known that there
are only a finite number of centers in $X$ of log canonical places
\cite{am}. Finally, we can identify the center $Z_{v}$ with a point
in the Chow variety $\mathcal{C}$ induced by a fixed ample line bundle
on $X.$ The group $G$ naturally acts algebraically on $\mathcal{C}$
and since $G$ is connected so is the $G$-orbit $GZ_{v},$ viewed
as a subset of $\mathcal{C}.$ But since the orbit $GZ_{v}$ is finite
this implies that it consists of a single point, i.e. $Z_{v}$ is
$G$-invariant, as desired. Finally, the last statement of the lemma
follows from the fact that $\mathcal{D}^{(N_{k})}$ is invariant under
the diagonal action of $\text{Aut}_{0}(X)$ on $X^{N_{k}}$ (see the
proof of Lemma \ref{lem:character}).
\end{proof}
To prove Proposition \ref{prop:non-excep text} we start by expressing
$-K_{(\P^{n},\Delta_{w})}\eqqcolon d\mathcal{O}(1),$ where $d=n+1-w(n+2)$
and $\Delta_{w}=w\Delta_{1}.$ The assumption $w=1-1/(n+1)$ is equivalent
to 
\[
d+w=1,
\]
 since $d+w=(n+1)(1-w).$ Likewise, rescaling
\[
\mathcal{D}^{(N)}\eqqcolon dD_{N},\,\,\Delta_{N}=w\Delta_{1,N},
\]
 where $\Delta_{1,N}$ is the divisor on $(\P^{n})^{N}$ defined by
$\pi_{1}^{*}\Delta_{1}+\cdots+\pi_{N}^{*}\Delta_{1},$ in terms of
the $N$ different projections $\pi_{i}\colon(\P^{n})^{N}\rightarrow\P^{n},$
we note that 
\begin{equation}
(i)\ \mathrm{lct}(\P^{n})^{N},\Delta_{1,N})=1,\,\,\,(ii)\ \mathrm{lct}((\P^{n})^{N},D_{N})\geq1.\label{eq:item one and two}
\end{equation}
The first item follows directly from the fact that $\Delta_{1,N}$
has simple normal crossings. To prove the second, one observe that
by definition, $D_{N}$ is linearly equivalent to $\pi_{1}^{*}\mathcal{O}(1)+\cdots+\pi_{N}^{*}\mathcal{O}(1)$
and since $\mathrm{lct}((\P^{n},D))\geq1$ for any $\Q$-divisor $D$
on $\P^{n}$ \cite[Example 1.3]{c-s} item $(ii)$ thus follows from
the product formula for the global log canonical threshold \cite[Lemma 2.29]{c-s}
Now fix a divisorial valuation $v$ over $(\P^{n})^{N}$ and note
that the inverse of $\mathrm{lct}(X,D)_{v},$ appearing in the valuative
formula \ref{eq:lct val} for $\mathrm{lct}(X,D),$ is linear wrt
$D.$ Hence, 
\begin{equation}
\frac{1}{\mathrm{lct}\left((\P^{n})^{N},\mathcal{D}^{(N)}+\Delta_{N}\right)_{v}}=d\frac{1}{\mathrm{lct}\left((\P^{n})^{N},D_{N}\right)_{v}}+w\frac{1}{\mathrm{lct}((\P^{n})^{N},\Delta_{1,N})_{v}}\leq1,\label{eq:ineq for inverse lct in pf}
\end{equation}
giving $\mathrm{lct}\left((\P^{n})^{N},\Delta;\mathcal{D}^{(N)}\right)_{v}\geq1.$
As a consequence, $\mathrm{lct}\left((\P^{n})^{N},\Delta;\mathcal{D}^{(N)}\right)\geq1.$
Now consider a divisorial valuation $v$ over $(\P^{n})^{N}$ computing
$\mathrm{lct}\left((\P^{n})^{N},\mathcal{D}^{(N)}+\Delta_{N}\right).$
If $\mathrm{lct}\left((\P^{n})^{N},D_{N}\right)_{v}>1,$ then $\mathrm{lct}\left((\P^{n})^{N},\mathcal{D}^{(N)}+\Delta_{N}\right)>1,$
by inequality \ref{eq:ineq for inverse lct in pf}. Consider finally
the case when $\mathrm{lct}\left((\P^{n})^{N},D_{N}\right)_{v}=1.$
This implies, by item $(ii)$ in formula \ref{eq:item one and two},
that $v$ computes $\mathrm{lct}\left((\P^{n})^{N},D_{N}\right).$
To prove that $\mathrm{lct}\left((\P^{n})^{N},\mathcal{D}^{(N)};\Delta\right)>1$
it will, by the inequality \ref{eq:ineq for inverse lct in pf}, be
enough to show that $\mathrm{lct}(\P^{n})^{N},\Delta_{1,N})_{v}>1.$
However, by the previous lemma, the center $Z_{v}$ of $v$ is invariant
under group $\mathrm{PGL}(n+1,\C)$ acting diagonally on $(\P^{n})^{N}.$
This, in fact, forces $v(\Delta_{1,N})=0$ and thus $\mathrm{lct}(\P^{n})^{N},\Delta_{1,N})_{v}=\infty.$
Indeed, assume, to get a contradiction, that $v(\Delta_{1,N})>0.$
Then the center $Z_{v}$ is contained in support of $\Delta_{1,N}.$
In particular, fixing any $(x_{1},...,x_{N})\in Z_{v}$ this implies,
since $Z_{v}$ is $\mathrm{PGL}(n+1,\C)$-invariant, that $(gx_{1},...,gx_{N})$
is contained in $\Delta_{1,N}$ for any $g\in\mathrm{PGL}(n+1,\C).$
This implies, since $\mathrm{PGL}(n+1,\C)$ acts transitively on $\P^{n},$
that the support of the divisor $\Delta_{1}$ is all of $\P^{n},$
which is a contradiction. 
\begin{rem}
\label{rem:except}Taking $N=1$ and replacing the divisor $\mathcal{D}^{(N)}$
in the previous proof with any $\Q$-divisor that is linearly equivalent
to $-K_{(\P^{n},\Delta_{w})}$ shows that $(\P^{n},\Delta_{w})$ is
exceptional iff $d+w<1.$ Indeed, when $d+w<1,$ i.e. when $w>1-1/(n+1),$
the inequality \ref{eq:ineq for inverse lct in pf} is strict and
when $d+w\leq1$ its left hand side is greater than or equal to $1.$
This is seen by taking $v$ to be the divisorial valuation on $\P^{n}$
defined by any of the hyperplanes $H_{i}$ in the support of $\Delta_{w}$.
\end{rem}

\subsection{The minimal case on $\mathbb{P}^{n}$}

As a preparation, we recall the following result \cite[Theorem 11.2]{dol}
describing the semistable locus $(\mathbb{P}^{n})_{\mathrm{ss}}^{N}$
of $(\mathbb{P}^{n})^{N}$ for general $N$ with respect to the diagonal
action of $\G=\mathrm{PGL}(n+1,\mathbb{C})$ and polarization $-K_{(\mathbb{P}^{n})^{N}}$:
A point $(x_{1},\dots,x_{N})\in(\mathbb{P}^{n})^{N}$ is semistable
iff for any proper linear subspace $W\subset\mathbb{P}^{n},$
\begin{equation}
\#\{i\colon x_{i}\in W\}\leq\frac{\mathrm{dim}W+1}{n+1}N.\label{eq: ss Pn}
\end{equation}

With this explicit description we can compute $\gamma^{(N)}(\mathbb{P}^{n})$
for the minimal value of $N$. 
\begin{prop}
\label{prop: Gibbs poly of minimal Pn}$\mathbb{P}^{n}$ is Gibbs
polystable at the minimal level $k_{\mathrm{min}}\coloneqq1/(n+1)$,
more precisely
\[
\gamma_{\mathrm{PS}}^{(N_{\mathrm{min}})}(\mathbb{P}^{n})=\infty.
\]
\end{prop}

\begin{proof}
Consider $X=\mathbb{P}^{n}$ at the minimal level $k_{\min}\coloneqq1/(n+1)$,
i.e. $-kK_{X}=\mathcal{O}(1)$ and $N_{\mathrm{min}}=n+1$. First
recall the well-known fact that the quantized Futaki character of
$\mathbb{P}^{n}$ vanishes at every level. Note that in homogeneous
coordinates $X_{ij}$, $1\le i,j\leq N$ on $(\mathbb{P}^{n})^{n+1}$,
$\mathrm{det}S^{(N_{\mathrm{min}})}=\mathrm{det}[X_{ij}]_{i,j}$.
Note that when $N=n+1$, condition \ref{eq: ss Pn} for $W$ of codimension
1 states that not all $x_{i}$ lie on $W$, or equivalently that $\mathrm{det}S^{(n+1)}(x_{1},\dots,x_{n+1})\neq0$.
Thus $(\mathbb{P}^{n})_{\mathrm{ss}}^{N_{\mathrm{min}}}\subset\{\mathbb{\mathrm{det}}S^{(N_{\mathrm{min}})}\neq0\}$
so we get
\[
\gamma_{\mathrm{PS}}^{(N_{\mathrm{min}})}(\mathbb{P}^{n})=\mathrm{lct}\Big((\mathbb{P}^{n})_{\mathrm{ss}}^{N_{\mathrm{min}}},\frac{1}{k_{\mathrm{min}}}\mathrm{det}S^{N_{\min}}\Big)=\infty.
\]
\end{proof}
In fact, when $N=n+1$ in \ref{eq: ss Pn}, the conditions for arbitrary
$W$ follows from the ones for $W$ of codimension 1 and thus $(\mathbb{P}^{n})_{\mathrm{ss}}^{N_{\mathrm{min}}}=\{\mathbb{\mathrm{det}}S^{(N_{\mathrm{min}})}\neq0\}$
in this special case. 

By Theorem \ref{thm:gibbs poly and lct intro}, $\gamma_{\mathrm{PS}}^{(N)}(X)$
is not increasing, in general. This phenomenon persists for $X=\mathbb{P}^{n}$
when $n\geq2$ by the following observation. Consider $X=\mathbb{P}^{2}$
at level $k=\frac{2}{3}$ so that $-kK_{\mathbb{P}^{2}}=\mathcal{O}(2)$.
Then it does not hold that $(\mathbb{P}^{2})_{\mathrm{ss}}^{6}\subset\{\mathrm{det}S^{(6)}\neq0\}$.
To see this, pick $6$ distinct points $(y_{1},\dots,y_{6})$ on a
smooth quadric in $\mathbb{P}^{2}$, then $\mathrm{det}S^{(6)}(y_{1},\dots,y_{6})=0$.
Since they are distinct, the only condition to be checked for semistability
is whether $5$ or more points are contained in a single line. But
a line meets a quadric in $\mathbb{P}^{2}$ in at most $2$ points
by Bezout, so at most two of the $y_{1},\dots,y_{6}$ are contained
in a common line. By \ref{eq: ss Pn}, $(y_{1},\dots,y_{6})$ is semistable,
giving
\[
\mathrm{lct}_{(y_{1},\dots,y_{6})}((\mathbb{P}^{2})^{6},\frac{3}{2}\mathrm{det}S^{(6)})\leq2/3
\]
and thus $\gamma^{(6)}(\mathbb{P}^{2})\leq2/3$ and $\mathbb{P}^{2}$
is not Gibbs polystable at level $k=2/3$.

\section{\label{sec:Symmetry-breaking,-thermodynamic}Analytic stability thresholds
and thermodynamics}

In this section we will, in particular, compute the analytic polystability
threshold $\delta_{\text{PS}}^{\mathrm{A}}(X,\Delta)$ (defined by
formula \ref{eq:def of anal st thr}) when $(X,\Delta)$ is a log
Fano curve, thus completing the proof of Theorem \ref{thm:gibbs poly and lct intro}.

\subsection{\label{subsec:Setup analytic}Setup }

We continue with the notation introduced in Section \ref{subsec:Setup Gibbs}.
In particular, the complex Lie group $\text{Aut}_{0}(X,\Delta)$ is
denoted by $\G$ and we will assume that it is\emph{ reductive, }i.e.
$\G$ is the complexification of any maximal compact subgroup of $\G.$
We will break the $\G$-symmetry by fixing a maximal compact subgroup
$\K$ of $\G;$ the corresponding Lie algebras will be denoted by
$\mathfrak{k}$ and $\mathfrak{g}.$ The elements of the complex Lie
algebra $\mathfrak{g}$ may be identified with holomorphic vector
fields $V^{1,0}$ on $X$ of type $(1,0)$ (whose flow preserves the
divisor $\Delta).$ Identifying a holomorphic vector field $V^{1,0}$
with the real vector field $\text{Re}V^{1,0}$ defined by its real
part, the reductivity of $\G$ translates into the following decomposition
of the Lie algebra $\mathfrak{g},$ viewed as a real vector space:
\[
\mathfrak{g}=\mathfrak{k}\oplus J\mathfrak{k},
\]
 where $J$ denotes the endomorphism of the real tangent bundle $TX$
defined by the complex structure on $X.$

Furthermore, we fix a $\K$-invariant metric $\phi_{0}$ on $-K_{(X,\Delta)}$
with positive curvature form $\omega_{0}(=dd^{c}\phi_{0})$ (using
additive notation for metrics on line bundles, recalled below). The
fixed metric $\phi_{0}$ on $-K_{(X,\Delta)}$ corresponds to a volume
form on $X$ with divisorial singularities that we denote by $\mu_{0},$
which (up to scaling the metric on $-K_{(X,\Delta)}$$)$ we may assume
is a probability measure; see Section \ref{subsec:Metrics vs vol forms with div}.
We will denote by $\mathcal{P}(X)$ the space of all probability measures
on $X,$ endowed with its weak topology and by $\mathcal{P}(X)_{0}$
the subspace of all measures $\mu$ with vanishing moment: 
\begin{equation}
\int_{X}\boldsymbol{m}\mu=0,\label{eq:moment cond on mu in text}
\end{equation}
where $\boldsymbol{m}$ denotes the moment map induced by the action
of $\K$ on $(X,-K_{(X,\Delta})$ and the symplectic form defined
by $\omega_{0}$ (see Section \ref{subsec:The-moment-map}). Functionals
whose arguments are measures $\mu$ will be denoted by capital letters
(e.g. $E(\mu)),$ while those whose arguments are functions $u$ or
metrics $\phi$ will be denoted by calligraphic letters, such as $\mathcal{E}(u).$ 

\subsection{Preliminaries}

\subsubsection{\label{subsec:Metrics vs vol forms with div}Metrics on $-K_{(X,\Delta)}$
versus volume forms with divisorial singularities}

We will use additive notation for tensor powers of line bundles $L\rightarrow X$
and their Hermitian metrics. This means that we identify a Hermitian
metric $\left\Vert \cdot\right\Vert $ on $L$ with a collection of
smooth local functions $\phi_{U}\coloneqq-\log(\left\Vert e_{U}\right\Vert ^{2})$
associated to a given covering of $X$ by open subsets $U$ and trivializing
holomorphic sections $e_{U}$ of $L\rightarrow U.$ The normalized
curvature form on $X$ of the metric may then, locally, be expressed
as 
\[
dd^{c}\phi_{U}\coloneqq\frac{i}{2\pi}\partial\bar{\partial}\phi_{U}
\]
Accordingly, as is customary, we will symbolically denote by $\phi$
a given Hermitian metric on $L$ and by $dd^{c}\phi$ its curvature
form, which represent the first Chern class $c_{1}(L)$ of $L.$ We
will denote by $\mathcal{H}(L)$ the space of all smooth metrics $\phi$
on $L$ such that $dd^{c}\phi>0$ and by $\text{PSH}(L)$ the space
of all \emph{psh metrics} on $L,$ i.e. $\phi_{U}$ is a plurisubharmonic
(psh) function on $X.$ If $\phi\in\mathrm{PSH}(L)$, then $dd^{c}\phi$
defines a positive current on $L.$ The subspace of $\text{PSH}(L)$
consisting of locally bounded metrics, i.e. $\phi_{U}\in L_{\text{loc}}^{\infty},$
will be denoted by $\text{PSH}(L)_{\mathrm{b}}.$ More generally,
when $L$ is merely defined as $\Q$-line bundle a metric $\phi$
on $L$ may be defined by the collection $\phi_{U}\coloneqq k^{-1}\phi_{U}^{(k)},$
where $\phi^{(k)}$ is a metric on $kL$ for a positive integer such
that $kL$ is a bona fide line bundle. 

Now assume that $(X,\Delta$) is a log Fano manifold and take $L$
to be the $\mathbb{Q}$-line bundle $-K_{(X,\Delta)}$. Then, following
\cite{bbegz}, any metric $\phi$ in $\mathrm{PSH}(L)_{\mathrm{b}}$
induces a measure $\mu_{\phi}$ on $X$, absolutely continuous with
respect to Lebesgue measure. Recall that we have fixed a (multivalued)
section $s_{\Delta}$ cutting out $\Delta$. Then, locally in a holomorphic
coordinate chart $U$,
\begin{equation}
\mu_{\phi}\coloneqq e^{-\phi_{U}}|s_{U}|^{-2}\frac{i^{n^{2}}}{2^{n}}dz\wedge d\bar{z},\,\,\,dz\coloneqq dz_{1}\wedge\cdots\wedge dz_{n}\label{eq: def of mu phi}
\end{equation}

by taking $e_{U}=\partial/\partial z_{1}\wedge\cdots\wedge\partial/\partial z_{n}\otimes e_{U}^{*}$
and $s_{\Delta}=s_{U}e_{U}$. A measure $\mu$ on $X$ will be called
a \emph{volume form with divisorial singularities (along $\Delta$)}
if $\mu=\mu_{\phi}$ for some $\phi$ in $\text{PSH}(L)_{\mathrm{b}}$
with the property that $\phi$ is smooth on $X-\Delta.$ 

\subsubsection{Kähler--Einstein metrics for $(X,\Delta)$ }

Given a log Fano manifold $(X,\Delta)$, a Kähler--Einstein metric
$\omega_{\text{KE }}$ for $(X,\Delta)$ is, by definition, a positive
current $\omega_{\text{KE}}\in c_{1}(-K_{(X,\Delta)})$ such that
$\omega_{\text{KE}}$ is the curvature of a metric $\phi_{\text{KE }}$
in $\text{PSH}(-K_{(X,\Delta)})_{\mathrm{b}}$ and such that the restriction
of $\omega_{\text{KE}}$ to the complement $X-\Delta$ of $\Delta$
in $X$ is a smooth Kähler--Einstein metric, i.e. $\text{Ric }\omega=\omega$
on $X-\Delta.$ In particular, $\omega_{\text{KE}}^{n}$ is a volume
form with divisorial singularities along $\Delta$ (see \cite[Appendix B]{bbegz}
for more general regularity results). When $(X,\Delta)$ is log smooth,
i.e. $\Delta$ has simple normal crossings, $\omega_{\text{KE}}$
has cone singularities along $\Delta$ \cite{g-p,m-r}. Moreover,
when $(X,\Delta)$ is a Fano orbifold $\omega_{\text{KE}}$ is smooth
in the orbifold sense \cite{g-k}. 

\subsubsection{Kähler metrics vs Kähler potentials and metrics on line bundles}

Given a Kähler form $\omega_{0}$ on a compact complex manifold, consider
the space of all\emph{ Kähler potentials:}
\[
\mathcal{H}(X,\omega_{0})\coloneqq\left\{ u\in C^{\infty}(X)\colon\,\,\omega_{u}\coloneqq\omega_{0}+dd^{c}u>0\right\} 
\]
Using the map $u\mapsto\omega_{u}$ we will identify the space of
all Kähler metrics $\omega$ in the cohomology class $[\omega_{0}]\in H^{2}(X,\R)$
with $\mathcal{H}(X,\omega_{0})/\R,$ under the additive action of
$\R$ on $\mathcal{H}(X,\omega_{0})$ (acting as $u\rightarrow u+C).$
Moreover, given a line bundle $L$ over $X$ and a metric $\phi_{0}$
on $L$ with curvature form $\omega_{0}$, expressing 
\begin{equation}
u=\phi-\phi_{0}\label{eq:u in terms of phi}
\end{equation}
 for a metric $\phi$ on $L,$ we may identify the space $\mathcal{H}(X,\omega_{0})$
with the space with $\mathcal{H}(L)$ of smooth metrics on $L$ with
positive curvature form. 

\subsubsection{\label{subsec:Pluripotential-theory-and}Pluripotential theory and
energy}

In this section we recall the notion of pluricomplex energy of a probability
measure $\mu$ introduced in \cite{bbgz}, but using the notation
in \cite{berm8,berm8 comma 5}. Given a Kähler form $\omega_{0}$
on a compact complex manifold $X,$ denote by $\text{PSH}(X,\omega_{0})$
the closure of $\mathcal{H}(X,\omega_{0})$ in $L^{1}(X),$ called
the space of all\emph{ $\omega_{0}$-plurisubharmonic functions. }There
exists an effective constant $C_{\text{sup}}$ only depending on $(X,\omega_{0})$
such that
\begin{equation}
\sup_{X}u\leq\int_{X}u\frac{\omega_{0}^{n}}{V}+C_{\text{sup}},\,\,\,V=\int_{X}\omega_{0}^{n}\label{eq:bound on sup u}
\end{equation}
(by the submean property of local plurisubharmonic functions and the
compactness of $X,$ or by estimating the Green's function of the
Laplacian wrt $\omega).$ Given a line bundle $L$ over $X$ and a
metric $\phi_{0}$ on $L$ with curvature form $\omega_{0}$ we can
identify $\text{PSH}(X,\omega_{0})$ with $\mathcal{H}(L)$ and $\text{PSH}(L),$
using formula \ref{eq:u in terms of phi}. Given a log Fano manifold
$(X,\Delta)$ we will take $L=-K_{(X,\Delta)}$ and $\phi_{0}$ as
in the setup introduced above. Consider the following functional on
$C^{\infty}(X)$: 
\begin{equation}
\mathcal{E}(u)\coloneqq\frac{1}{(n+1)V}\int_{X}\sum_{j=0}^{n}u\omega_{u}^{n-j}\wedge\omega_{0}^{j}\label{eq:def of beaut E}
\end{equation}
 whose differential on $\mathcal{H}(X,\omega_{0})$ is given by 
\begin{equation}
d\mathcal{E}|_{u}=\omega_{u}^{n}/V.\label{eq:diff of beaut E}
\end{equation}
As a consequence, the restriction of $\mathcal{E}$ to $\mathcal{H}(X,\omega_{0})$
is increasing and thus admits a canonical extension to $\text{PSH}(X,\omega_{0}),$
namely the smallest extension which is upper semi-continuous with
respect to the $L^{1}$-topology \cite{bbgz}. The subspace of $\text{PSH}(X,\omega_{0})$
consisting of all $u$ such that $\mathcal{E}(u)>-\infty$ is denoted
by $\mathcal{E}^{1}(X,\omega_{0}).$ Formula \ref{eq:diff of beaut E}
still applies on $\mathcal{E}^{1}(X,\omega_{0}),$ when $\omega_{u}^{n}$
is taken to be the Monge--Ampère measure defined in terms of non-pluripolar
products \cite{bbgz}.

The \emph{pluricomplex energy} $E(\mu)$ of a probability measure
$\mu$ on $X$ (with respect to $\omega_{0})$ is defined by the following
lower semi-continuous (lsc) convex function on $\mathcal{P}(X)$:
\begin{equation}
E(\mu)\coloneqq\sup_{u\in C^{0}(X)\cap\text{PSH}(X,\omega_{0})}\mathcal{E}(u)-\left\langle u,\mu\right\rangle .\label{eq:def of E as sup}
\end{equation}
 By \cite{bbgz}, $E(\mu)<\infty$ iff there exists $u_{\mu}\in\text{PSH}(X,\omega_{0})$
realizing the sup above such that $\mathcal{E}(u_{\mu})>-\infty$
(but $u_{\mu}$ need not be in $C^{0}(X)).$ Moreover, up to an additive
constant, $u_{\mu}$ is uniquely determined by 
\begin{equation}
\omega_{u_{\mu}}^{n}/V=\mu,\label{eq:MA eq for u mu}
\end{equation}
 using the non-pluripolar product introduced in \cite{begz} to define
$\omega_{u_{\mu}}^{n}$. We will fix the additive constant by demanding
that $\int u_{\mu}\omega_{0}^{n}=0$ and call $u_{\mu}$ the \emph{potential}
(function) of $\mu.$ Equivalently, the corresponding potential metric
$\phi_{\mu}$ satisfies
\begin{equation}
(dd^{c}\phi_{\mu})_{\,}^{n}/V=\mu,\,\,\,\,\int(\phi_{\mu}-\phi_{0})(dd^{c}\phi_{0})^{n}=0.\label{eq:def of phi mu}
\end{equation}
By the Calabi--Yau theorem $\mu$ is a smooth volume form iff $\phi_{\mu}\in\mathcal{H}(L)$
and we denote by $\mathcal{V}(X)\subset\mathcal{P}(X)$ the space
of probability volume forms. The functional $E$ on $\mathcal{P}(X)$
satisfies 
\begin{equation}
E(\omega_{0}^{n}/V)=0\label{eq:E is zero on omega zero}
\end{equation}
 and on the convex subspace $\{E<\infty\}\subset\mathcal{P}(X)$ its
differential may be represented as follows \cite{berm6}: 
\begin{equation}
dE\vert_{\mu}=-u_{\mu}\,\,\,\left(=\phi_{0}-\phi_{\mu}\right),\label{eq:diff of E}
\end{equation}
For example, any measure $\mu$ on $X$ with divisorial singularities
(or more, generally with an $L^{p}$-density for some $p>1)$ has
finite energy, $E(\mu)<\infty.$ Thus the properties \ref{eq:E is zero on omega zero}
and \ref{eq:diff of E} determine $E$ on $\mathcal{P}(X)\cap L_{\text{loc}}^{p}$
for any $p>1$ (and the functional $E$ on $\mathcal{P}(X)$ can then
be characterized as the greatest lower semi-continuous extension of
$E$ from $\mathcal{P}(X)\cap L_{\text{loc}}^{p}$ to $\mathcal{P}(X)$). 
\begin{rem}
\label{rem:plurienergy is I minus J}In terms of standard functionals
in Kähler geometry, when $\mu$ is smooth, 
\begin{equation}
E(\mu)=(\mathcal{I}-\mathcal{J})(u_{\mu})=\mathcal{E}(u_{\mu})-\left\langle u_{\mu},\mu\right\rangle ,\label{eq:E as I minus J}
\end{equation}
 where $\mathcal{I}$ and $\mathcal{J}$ are the standard energy functionals
in Kähler geometry, defined by 
\[
\mathcal{J}(u)\coloneqq-\mathcal{E}(u)+\int_{X}u\frac{\omega_{0}^{n}}{V},\,\,\,\,\,\mathcal{I}(u)\coloneqq\int_{X}u\frac{(\omega_{0}^{n}-\omega_{u}^{n})}{V}.
\]
All three functionals $\mathcal{I},(\mathcal{I}-\mathcal{J})$ and
$\mathcal{J}$ are comparable. In particular, by \cite[Lemma 2.2]{bbgz},
\begin{equation}
n^{-1}\mathcal{J}(u)\leq(\mathcal{I}-\mathcal{J})(u)\leq n\mathcal{J}(u)\label{eq:comparison I minus J and J}
\end{equation}
\end{rem}

\begin{example}
\label{exa:energy for n one}When $n=1$ integration by parts gives
\begin{equation}
\mathcal{E}(u)=-\frac{1}{2V}\int_{X}du\wedge d^{c}u+\int_{X}u\frac{\omega_{0}}{V},\,\,E(\mu)=\frac{1}{2V}\int_{X}du_{\mu}\wedge d^{c}u_{\mu}=-\frac{1}{2}\int_{X}u_{\mu}\mu\label{eq:beaut E on Riemann surf in terms of Dirich}
\end{equation}
if $E(\mu)<\infty$, or equivalently, $u_{\mu}$ is in the Sobolev
space $W^{2,1}(X)$ consisting of all $u\in L^{2}(X)$ such that $du\in L^{2}(X).$
Moreover, when $n=1$ the sup in formula \ref{eq:def of E as sup}
can be taken over all of $C^{0}(X)$ or over $W^{2,1}(X).$ Indeed,
by classical Hilbert space theory, the latter sup is attained by $u_{\mu}\in W^{2,1}(X)$
solving the (Laplace) equation \ref{eq:MA eq for u mu} and hence
$u_{\mu}\in\text{PSH}(X,\omega_{0})$.
\end{example}

\subsubsection{\label{subsec:The-actions-of}The actions of $\G$ on $\text{\ensuremath{\mathrm{PSH}}(}L),$
$\mathrm{PSH}(X,\omega_{0})$ and $\mathcal{P}(X)$ }

The group $\G$ naturally acts from the right on $\text{\ensuremath{\text{PSH}}(}L)$
by pullback: $\phi\cdot g\coloneqq g^{*}\phi.$ Fixing a reference
metric $\phi_{0}$ in $L$ and setting $\omega_{0}\coloneqq dd^{c}\phi_{0},$
the right action of $\G$ on $\text{\ensuremath{\text{PSH}}(}L)$
induces a right action on $\text{PSH}(X,\omega_{0}),$ by identifying
$\text{\ensuremath{\text{PSH}}(}L)$ with $\text{PSH}(X,\omega_{0}).$
In other words, $u\cdot g\coloneqq g^{*}(\phi_{0}+u)-\phi_{0}$ (which
is thus different than the action of $\G$ on $\text{PSH}(X,\omega_{0})$
by pullback). The group $\G$ also naturally acts from the right on
$\mathcal{P}(X)$ by setting $\mu\cdot g\coloneqq(g_{-1})_{*}\mu\coloneqq g^{*}\mu.$
This action is compatible with the action of $\G$ on $\text{\ensuremath{\text{PSH}}(}L)$
in the following sense: if $\mu=(dd^{c}\phi)^{n}/V,$ then $g^{*}\mu=(dd^{c}g^{*}\phi)^{n}/V.$
Moreover, if $L=-K_{(X,\Delta)}$ then $g^{*}\mu_{\phi}=\mu_{g^{*}\phi}.$
Finally, there is also a natural \emph{left} action of $\G$ on $\text{\ensuremath{\text{PSH}}(}L),$
defined by $g\cdot\phi\coloneqq(g^{-1})^{*}\phi,$ which is compatible
with the natural left action of $\G$ on $\mathcal{P}(X),$ defined
by $g\cdot\mu\coloneqq g_{*}\mu.$ 

\subsubsection{\label{subsec:The-space-of geod}The space $\mathcal{R}_{1,1}(X,\omega)$
of geodesics rays in $\mathcal{H}_{1,1}(X,\omega)$ and the $L^{1}$-Finsler
completion of $\mathcal{H}(X,\omega)$}

The space $\mathcal{H}(X,\omega_{0})$ of Kähler potentials can be
endowed with a canonical Riemannian metric, first introduced by Mabuchi,
where the length of a tangent vector $\dot{u}$ at $u$ is defined
by 
\begin{equation}
\left\Vert \dot{u}\right\Vert _{u}\coloneqq\left(\int_{X}\dot{u}^{2}\omega_{u}^{n}\right)^{1/2}\label{eq:def of Mab metric}
\end{equation}
 The corresponding notion of geodesics was extended to the space 
\[
\mathcal{H}_{1,1}(X,\omega_{0})\coloneqq\left\{ u\colon\omega_{u}\in L_{\text{loc}}^{\infty},\,\,\omega_{u}\geq0\right\} .
\]
in \cite{ch}. We will denote by $\mathcal{R}(X,\omega_{0})_{1,1}$
the corresponding space of all \emph{geodesic rays} $u_{t}$ in $\mathcal{H}_{1,1}(X,\omega_{0}),$
parametrized by $t\in[0,\infty[.$ More generally, following \cite{da},
fixing $p\in[1,\infty[$ and replacing the $L^{2}$-norm in formula
\ref{eq:def of Mab metric} with the corresponding $L^{p}$-norm yields
a Finsler metric on $\mathcal{H}(X,\omega_{0}),$ which induces a
metric $d_{p}$ on $\mathcal{H}(X,\omega_{0}).$ As shown in \cite{da},
for $p=1,$ the metric space completion $\overline{(\mathcal{H}(X,\omega_{0}),d_{1})}$
of $(\mathcal{H}(X,\omega_{0}),d_{1})$ may be identified with $\mathcal{E}^{1}(X,\omega_{0})$
(defined in the previous section). Thus $\overline{(\mathcal{H}(X,\omega_{0}),d_{1})}$
contains $\mathcal{H}_{1,1}(X,\omega_{0})$ and, as shown in \cite{da},
the geodesics in $\mathcal{H}_{1,1}(X,\omega_{0}),$ discussed above,
yield constant speed geodesics in the metric space $\overline{(\mathcal{H}(X,\omega_{0}),d_{1})}$.
Moreover, on the subset $\{\mathcal{E}(u)=0\}$ of $\mathcal{E}^{1}(X,\omega_{0})$
the following inequalities hold: 
\begin{equation}
A^{-1}E\left(\frac{\omega_{u}^{n}}{V}\right)-B\leq d_{1}(0,u)\leq AE\left(\frac{\omega_{u}^{n}}{V}\right)+B\label{eq:comparison E and d one}
\end{equation}
 for some positive constants $A$ and $B.$ Indeed, by \ref{eq:E as I minus J},
$E(\frac{\omega_{u}^{n}}{V})=(\mathcal{I}-\mathcal{J})(u)$, which
is comparable to $\mathcal{J}(u)$ \cite[Lemma 2.13]{berm6} and by
\cite[Prop 5.5]{d-r} the inequalities \ref{eq:comparison E and d one}
also hold when $E(\frac{\omega_{u}^{n}}{V})$ is replaced by $\mathcal{J}(u)$
(see also \cite[Remark 3.9]{d-r}).

\subsubsection{\label{subsec:The-Mabuchi-f vs free energ}The Mabuchi functional
$\mathcal{M}$ versus the free energy functional $F_{\beta}$ }

Given an ample line bundle $L$, the Mabuchi functional $\mathcal{M}$
was originally defined when $\Delta$ is trivial, by declaring that
its differential on $\mathcal{H}(L)$ at $\phi$ is identified with
$-(R_{\phi}-\overline{R})(dd^{c}\phi)^{n},$ where $R_{\phi}$ denotes
the scalar curvature of the Kähler metric $dd^{c}\phi$ and $\overline{R}$
its average \cite{mab}. Here we will instead follow the thermodynamical
formalism in \cite{berm6,bbegz}, which applies in the general setup
of log Fano varieties $(X,\Delta).$ We thus define the Mabuchi functional
$\mathcal{M}$ on $\text{PSH}(-K_{(X,\Delta)})_{\mathrm{b}}$ with
respect to the log pair $(X,\Delta)$ and the fixed reference metric
$\phi_{0}\in\text{PSH}(-K_{(X,\Delta)})_{\mathrm{b}}$, by 
\begin{equation}
\mathcal{M}(\phi)\coloneqq F((dd^{c}\phi)_{\,}^{n}/V)\coloneqq F_{-1}((dd^{c}\phi)_{\,}^{n}/V)\label{eq:M as free energ}
\end{equation}
 where $F_{\beta}(\mu)$ is the \emph{free energy functional} on $\mathcal{P}(X)$
at \emph{inverse temperature} $\beta\in\R,$ defined as $F_{\beta}(\mu)=\infty$
if $E(\mu)=\infty$ and when $E(\mu)<\infty,$ 
\[
F_{\beta}(\mu)\coloneqq\beta E(\mu)+D(\mu|\mu_{0})\in]-\infty,\infty],
\]
where $\mu_{0}\coloneqq\mu_{\phi}$ and where $D(\mu|\nu)$ denotes
the \emph{entropy of $\mu$ relative to $\nu$ }for given measures
$\mu$ and $\nu$ on a topological space $Y,$ also known as the \emph{Kullback--Leibler
divergence:} 
\begin{equation}
D(\mu|\nu)=\int_{X}\log\left(\frac{\mu}{\nu}\right)\mu\label{eq:def of rel entropy D}
\end{equation}
if $\mu$ is absolutely continuous with respect to $\nu$ and otherwise
$D(\mu|\nu)\coloneqq\infty.$ Set
\begin{equation}
F(\beta)\coloneqq\inf_{\mu\in\mathcal{P}(X)}F_{\beta}(\mu),\,\,\,\,F(\beta)_{0}\coloneqq\inf_{\mu\in\mathcal{P}(X)_{0}}F_{\beta}(\mu)\label{eq:def of F beta}
\end{equation}
We recall that the function $(\mu,\nu)\mapsto D(\mu|\nu)$ is non-negative
and lower semicontinuous (lsc) on $\mathcal{P}(X)^{2}$ and $\mu\mapsto D(\mu|\nu)$
is convex on $\mathcal{P}(X).$ These properties follow directly from
the following well-known Legendre duality formula \cite{d-z}:
\begin{equation}
D(\mu|\nu)=\sup_{u\in C^{0}(X)}\left\{ \int u\mu-\log\int e^{u}\nu\right\} .\label{eq:D as leg tra}
\end{equation}

\subsubsection{\label{subsec:The-Futaki-invariants}The Futaki character of $(X,\Delta)$
and its quantization}

Let $V^{1,0}$ be a holomorphic vector field on $X$ of type $(1,0)$
corresponding to an element in $\mathfrak{g}.$ Its \emph{Futaki invariant}
$\text{Fut}_{(X,\Delta)}(V^{1,0})$ may be defined by 
\begin{equation}
\text{Fut}_{(X,\Delta)}(V^{1,0})\coloneqq-\int_{X}\mathcal{L}_{V^{1,0}}\phi_{0}\frac{(dd^{c}\phi_{0})^{n}}{V}\in\C,\label{eq:def of Fut}
\end{equation}
 where $\mathcal{L}_{V^{1,0}}\phi_{0}$ denotes the Lie derivative
of $\phi_{0}$ along $V^{1,0},$ defined by decomposing $V^{1,0}=\text{Re}V^{1,0}+i\text{Im}V^{1,0}$
and setting 
\begin{equation}
\mathcal{L}_{W}\phi\coloneqq\lim_{t\rightarrow0}\left(\exp(tW)^{*}\phi_{0}-\phi_{0}\right)\label{eq:def of Lie on phi}
\end{equation}
 for any real vector field $W$ on $X$ that lifts to $-K_{(X,\Delta)}$,
and extended $\mathbb{C}$-linearly to any complex vector field on
$X$ that lifts to $-K_{(X,\Delta)}$. The corresponding complex-valued
function $\text{Fut}_{(X,\Delta)}$ on the Lie algebra $\mathfrak{g},$
defining an element in $\mathfrak{g}^{*},$ is called the\emph{ Futaki
character }and is independent of the choice of $\phi_{0}$: 
\begin{lem}
\label{lem:Fut}For any given $\psi,\phi\in\mathrm{PSH}(X,-K_{(X,\Delta)})_{\mathrm{b}}$
and holomorphic vector field $V^{1,0}$ on $X$ corresponding to an
element in $\mathfrak{g}$,
\begin{equation}
\text{\ensuremath{\mathrm{Fut}}}_{(X,\Delta)}(V^{1,0})=-\int_{X}\mathcal{L}_{V^{1,0}}\psi\frac{(dd^{c}\psi)^{n}}{V}=-\int_{X}\mathcal{L}_{V^{1,0}}\phi\frac{(dd^{c}\phi)^{n}}{V}.\label{eq:pf of lemma fut first}
\end{equation}
As a consequence, the definition \ref{eq:def of Fut} of $\text{Fut}_{(X,\Delta)}(V^{1,0})$
is independent of the choice of $\phi_{0}$ and $\text{Fut}_{(X,\Delta)}(\cdot)\equiv0$
iff $\omega_{0}^{n}/V\in\mathcal{P}(X)_{0}.$ Moreover, 
\begin{equation}
\frac{d}{dt}\mathcal{M}\left(\exp(t\mathrm{Re}V^{1,0})^{*}\psi\right)\Big|_{t=0}=\mathrm{Re}\,\ensuremath{\mathrm{Fut}}_{(X,\Delta)}(V^{1,0})\label{eq:pf of lemma fut second}
\end{equation}
More generally, \ref{eq:pf of lemma fut first} holds if $\psi$ and
$\phi$ have finite energy and \ref{eq:pf of lemma fut second} holds
if $(dd^{c}\psi)^{n}/V$ has finite entropy. In particular, if $(X,\Delta)$
admits a KE metric, then $\text{Fut}_{(X,\Delta)}(\cdot)\equiv0.$ 
\end{lem}

\begin{proof}
This is well-known, at least for smooth data, but for completeness
we provide simple proofs. To prove \ref{eq:pf of lemma fut first}
first note that, for any given $g\in\G,$ 
\[
\mathcal{E}\left(g^{*}\psi\right)-\mathcal{E}\left(g^{*}\phi\right)=\mathcal{E}(\psi)-\mathcal{E}(\phi).
\]
This follows directly from the fact that the differential $d\mathcal{E},$
which is the one-form on $\mathcal{H}(X,\omega_{0})$ defined by formula
\ref{eq:diff of beaut E}, satisfies $g^{*}(d\mathcal{E})=d(g^{*}\mathcal{E}).$
Replacing $g$ with a one-parameter subgroup $g_{\tau}$ and differentiating
with respect to $\tau$ then proves formula \ref{eq:pf of lemma fut first}.
Next, observe that $\text{Fut}_{(X,\Delta)}(V^{1,0})=0$ for all $V^{1,0}\in\mathfrak{g}$
iff $\int_{X}\mathcal{L}_{\text{Re}V^{1,0}}\phi_{0}(dd^{c}\phi_{0})^{n}=0$
for all $V^{1,0}\in\mathfrak{g}$ (by applying $\text{Fut}_{(X,\Delta)}$
to both $V^{1,0}$ and $iV^{1,0}).$ Now, using that $\mathfrak{g}$
is reductive, we can decompose $\text{Re}V^{1,0}=W_{1}+JW_{2}$ where
$W_{i}\in\mathfrak{k}$, the Lie algebra of $\K$. Since $\phi_{0}$
is assumed $\K$-invariant $\text{Fut}_{(X,\Delta)}(V^{1,0})=\int h_{W_{2}}(\omega_{0})^{n}/V,$
which shows that $\text{Fut}_{(X,\Delta)}(\cdot)\equiv0$ iff $\int_{X}\boldsymbol{m}(\omega_{0})^{n}/V=0,$
i.e. iff $\omega_{0}^{n}/V\in\mathcal{P}(X)_{0}.$ Finally, applying
formula \ref{eq:M as free energ}, we can decompose
\[
\mathcal{M}(\psi)=-\mathcal{E}(\psi)+D((dd^{c}\psi)^{n}/V|\mu_{\psi}),
\]
expressed in terms of the functional $\mathcal{E}$ on $\text{PSH}(L)_{\mathrm{b}}$
and the entropy of the measure $(dd^{c}\psi)^{n}/V$ relative to the
measure $\mu_{\psi},$ defined by formula \ref{eq: def of mu phi}.
For any $g\in\G$ we have $g^{*}\mu_{\psi}=\mu_{g^{*}\psi}.$ Hence,
setting $g_{t}\coloneqq\exp(tV),$ 
\[
\mathcal{M}(g_{t}^{*}\psi)=-\mathcal{E}(g_{t}^{*}\psi)+D((dd^{c}\psi)^{n}/V|\mu_{\psi}).
\]
Differentiating with respect to $t$ and using \ref{eq:def of beaut E}
thus proves \ref{eq:pf of lemma fut second}. More generally, the
proof of \ref{eq:pf of lemma fut first} only used that $\phi$ and
$\psi$ have finite energy. Moreover, when $(dd^{c}\psi)^{n}$ has
finite entropy it is shown in \cite{bbgz} that $\psi$ has finite
energy, showing that $D((dd^{c}\psi)^{n}|\mu_{\psi})<\infty.$ Finally,
if $(X,\Delta)$ admits a Kähler--Einstein metric $\psi,$ then,
by \cite{bbegz}, $(dd^{c}\psi)^{n}$ has finite entropy and $\psi$
minimizes the Mabuchi functional on the space of metrics with finite
energy. Hence, $\frac{d}{dt}\mathcal{M}\left(\exp(t\text{Re}V^{1,0})^{*}\psi\right)\Big|_{t=0}$,
showing that $\text{Fut}_{(X,\Delta)}(\cdot)\equiv0.$
\end{proof}
The independence of the choice of $\phi_{0}$ is also a consequence
of the fact that $\text{Fut}_{(X,\Delta)}$ arises as the leading
term of the large-$k$ asymptotics of the \emph{quantized Futaki character}
$\text{Fut}_{(X,\Delta),k}(\cdot)$ at level $k,$ defined as follows
when $k$ is a positive integer such that $k\Delta$ is divisor:
\begin{equation}
\text{Fut}_{(X,\Delta),k}(V^{1,0})\coloneqq-\frac{1}{kN}\text{Tr}(\mathcal{L}_{V^{1,0}}),\label{eq:def of quant Fut}
\end{equation}
where $\mathcal{L}_{V^{1,0}}$ now denotes the endomorphism of $H^{0}(X,-kK_{(X,\Delta)})$
induced by the infinitesimal action of $V^{1,0}$ on $H^{0}(X,-kK_{(X,\Delta)})$
(see the proof of \cite[Prop 4.7]{b-n}). More generally, when it
is only assumed that $-kK_{(X,\Delta)}$ is a line bundle we define
$\text{Fut}_{(X,\Delta),k}$ by formula \ref{eq:def of real character}
in Lemma \ref{lem:character}.

\subsubsection{\label{subsec:The-moment-map}The moment maps $\boldsymbol{m}$ and
$\boldsymbol{m}_{N}$}

The compact group $\K$ fixing the symplectic form $\omega_{0}$ naturally
acts on $(X,-K_{(X,\Delta)})$ and thus induces a\emph{ moment map
}$\boldsymbol{m}$ from $X$ into the dual $\mathfrak{k}^{*}$ of
the Lie algebra $\mathfrak{k}$ of $\K.$ This map is $\K$-equivariant
(with respect to the coadjoint action on $\mathfrak{k}^{*}$) \cite{ho}.
Fixing an $\text{ad}\,\K$-invariant $L^{2}$-norm $\left|\cdot\right|$
on $\mathfrak{k}$ and an orthonormal basis $e_{1},\dots,e_{r}$ in
$\mathfrak{k}$ we can identify $\mathfrak{k}^{*}$ with $\R^{r}$
and represent 
\[
\boldsymbol{m}\colon\,\,X\rightarrow\R^{r},\,\,\,x\mapsto(h_{e_{1}}(x),\dots,h_{e_{r}}(x)),
\]
where $h_{W}$ denotes the Hamiltonian of a real vector field $W$
on $X$ corresponding to an element $W$ in $\mathfrak{k},$ defined
by $h_{W}\coloneqq\mathcal{L}_{JW}\phi$ (defined as in formula \ref{eq:def of Lie on phi}).
We will denote by $\boldsymbol{m}_{N}$ the moment map of the diagonal
action of the group $\K$ on $X^{N},$ endowed with the symplectic
form 
\[
\omega_{N}\coloneqq\frac{1}{N}\sum_{i=1}^{N}\omega(x_{i}),
\]
 where $\omega(x_{i})$ denotes the pullback of $\omega$ on $X$
to $X^{N}$ under the projection to the $i^{\mathrm{th}}$ factor
of $X^{N}$. Concretely,
\[
\boldsymbol{m}_{N}\colon X^{N}\rightarrow\R^{r},\,\,\,\boldsymbol{m}_{N}(x_{1},\dots,x_{N})=\frac{1}{N}\sum_{i=1}^{N}\boldsymbol{m}(x_{i})\in\R^{r}
\]
By the Kempf--Ness theorem \cite{ho}, the $\G$-orbit of $\left\{ \boldsymbol{m}_{N}=0\right\} $
coincides with the \emph{polystable} locus $(X^{N})_{\text{ps }}$
of $X^{N}$ (the subset of points $x$ in the semistable locus $(X)_{\text{ss}}^{N}$
such that $\G x$ is closed as a subset of $(X)_{\mathrm{ss}}^{N}$):
\begin{equation}
\G\left\{ \boldsymbol{m}_{N}=0\right\} =(X^{N})_{\text{ps }}\subset(X^{N})_{\mathrm{ss}}.\label{eq:G-orbit of moment zero}
\end{equation}
Moreover, $0$ is a regular value of $\boldsymbol{m}_{N}$ iff all
the stabilizers in $\K$ of the points in $\left\{ \boldsymbol{m}_{N}=0\right\} $
are finite iff the semistable locus of $X^{N}$ coincides with the
polystable locus \cite[Cor 3.10]{ho}. 

\subsection{\label{subsec:The-structure-of}The structure of the space $\mathcal{P}(X)_{0}$
of probability measures with vanishing moment}

With the preliminary results in place we turn to the proof of the
following result, which will play a key role in what follows. 
\begin{prop}
\label{lem:inf of E over G}Assume that the Futaki character of $(X,\Delta)$
vanishes. Then, given $\nu\in\mathcal{P}(X)$ such that $E(\nu)<\infty,$
there exists $\mu\in\G\nu$ minimizing $E$ on $\G\nu.$ Moreover,
$\mu$ minimizes $E$ on $\G\nu$ iff $\mu\in\G\nu\cap\mathcal{P}(X)_{0}$
and such a measure $\mu$ is uniquely determined modulo the action
of $\K.$ 
\end{prop}

\begin{proof}
First note that $E$ is $\K$-invariant (since $\omega_{0}$ is $\K$-invariant).
Given $\nu\in\mathcal{P}(X)$ we thus obtain a function on the quotient
$\K\backslash\G$ of $\G$ by the left action of\emph{ $\K$}: 
\begin{equation}
\K\backslash\G\rightarrow\R,\,\,\,f(g)\coloneqq E(g_{*}\nu),\label{eq:def of as E on G mod K}
\end{equation}
Since $\G=\K\exp(J\mathfrak{k})$ (by polar decomposition of complex
reductive Lie groups) the function $f$ can be identified with a function
on $\exp(J\mathfrak{k}).$ Next, let us show that, given $\nu\in\mathcal{P}(X)$
such that $E(\nu)<\infty,$ $\mu$ is a critical point of $E$ on
the orbit $\G\nu$ iff $\mu\in\mathcal{P}(X)_{0}.$ Since $E$ is
$\K$-invariant it will be enough to show that for any real vector
field $W$ on $X$ induced by an element in $\mathfrak{k}\subset\mathfrak{g}$
the following formula holds for any given $\mu\in\mathcal{P}(X)$
such that $E(\mu)<\infty$: 
\begin{equation}
\left.\frac{df(\exp(tJW))}{dt}\right|_{t=0}\coloneqq\frac{d}{dt}E(\exp(tJW)_{*}\mu)\Big|_{t=0}=\int h_{W}\mu,\label{eq:time deriv of E in pf}
\end{equation}
where $h_{W}$ is the Hamiltonian function corresponding to $W$ and
$\omega_{0}.$ We start by proving the formula when $\mu$ is a smooth
form. Set $\mu_{t}\coloneqq\exp(tJW)_{*}\mu.$ By \ref{eq:diff of E},
denoting by $\phi$ the potential of $\mu,$ 
\[
\frac{d}{dt}E(\mu_{t})|_{t=0}=-\int(\phi-\phi_{0})\frac{d\mu_{t}}{dt}\Big|{}_{t=0}=\int(\phi-\phi_{0})\mathcal{L}_{JW}\mu.
\]
 Integrating by parts and using that $h_{W}\coloneqq\mathcal{L}_{JW}\phi_{0}$
(see Section \ref{subsec:The-moment-map}) gives
\[
\frac{d}{dt}E(\mu_{t})|_{t=0}=-\int\mathcal{L}_{JW}(\phi-\phi_{0})\mu=-\int\mathcal{L}_{JW}\phi\frac{(dd^{c}\phi)^{n}}{V}+\int h_{W}\mu,
\]
which proves \ref{eq:time deriv of E in pf}, using the previous lemma
and the assumption that the Futaki character of $(X,\Delta)$ vanishes.
Now assume only that $E(\mu)<\infty$ and take a sequence $u_{j}$
in $\mathcal{H}(X,\omega_{0})$ decreasing to the potential $u_{\mu}$
of $\mu.$ Since $u_{\mu}\in\mathcal{E}^{1}(X,\omega)$ it follows
that $\mu_{j}\coloneqq\omega_{u_{j}}^{n}/V\rightarrow\mu$ and $E(\mu_{j})\rightarrow E(\mu),$
by \cite[Prop 1.1]{bbgz}. Likewise, given an automorphism $g$ of
$X,$ the sequence $g\cdot u_{j}$ in $\mathcal{H}(X,\omega_{0})$
decreases to $g\cdot u_{\mu}$ and thus $g_{*}\mu_{j}\rightarrow g_{*}\mu$
and $E(g_{*}\mu_{j})\rightarrow E(g_{*}\mu)$. Now, applying the previous
case to $\mu_{j}$ and $tW$ for all $t\in]0,\epsilon]$ gives 
\[
\epsilon^{-1}\left(E(\exp(\epsilon JW)_{*}\mu_{j})-E(\mu_{j})\right)=\epsilon^{-1}\int_{0}^{\epsilon}\left(\int_{X}h_{tW}\mu_{j}\right)dt.
\]
 Letting $j\rightarrow\infty$ and using that $h_{tW}$ is continuous
thus gives 
\[
\epsilon^{-1}\left(E(\exp(\epsilon JW)_{*}\mu)-E(\mu)\right)=\epsilon^{-1}\int_{X}\left(\int_{0}^{\epsilon}h_{tW}dt\right)\mu.
\]
 Finally, letting $\epsilon\rightarrow0$ and using that $h_{tW}(x)$
is continuous wrt $(t,x)$ proves formula \ref{eq:time deriv of E in pf}.

Now assume that $\nu\in\mathcal{P}(X)$ and $E(\nu)<\infty.$ Then
$E(g_{*}\nu)<\infty$ for all $g\in\G$ and we will next show that
the function $f,$ defined by formula \ref{eq:def of as E on G mod K},
is \emph{strictly convex}, when identified with a function on the
real vector space $\exp(J\mathfrak{k}).$ To this end fix a real vector
field $W$ on $X$ induced by an element in $\mathfrak{k}\subset\mathfrak{g}.$
Differentiating formula \ref{eq:time deriv of E in pf} gives 
\[
\frac{d^{2}}{d^{2}t}\Bigg|_{t=0}E(\exp(tJW)_{*}\nu)=\int_{X}dH(JW)\nu=\int_{X}\omega_{0}(W,JW)\nu,
\]
 using that $h$ is a Hamiltonian function for $W.$ Since $\omega_{0}(W,JW)$
is non-negative, applying the previous formula to any element in $\exp(\R JW)_{*}\nu$
shows that $t\mapsto E(\exp(tJW)_{*}\nu)$ is a smooth convex function,
which is, in fact, \emph{strictly} convex. Indeed, $\omega(W,JW)$
only vanishes on the subvariety of $X$ defined by the zero-locus
of $W$ and $\nu$ does not put mass on analytic subvarieties (since
$E(\nu)<\infty,$ $\nu$ does not charge pluripolar subsets \cite{bbgz}).

It follows directly from the strict convexity of $f$ that its critical
points are unique. All that remains is thus to prove that $f$ is
proper on \emph{$\K\backslash\G$} (and thus admits a critical point).
First assume that $\nu=\omega_{0}^{n}/V,$ which minimizes $E$ on
all of $\mathcal{P}(X).$ Then $f$ is a strictly convex function
admitting a minimum (when $g$ is the identity), which implies that
$f$ is proper. Next, for any given $\nu\in\mathcal{P}(X)$ such that
$E(\nu)<\infty$ there exists a unique $u\in\mathcal{E}^{1}(X,\omega)$
such that $\omega_{u}^{n}/V=\nu$ and $\mathcal{E}(u)=0.$ Since we
have assumed that the Futaki character of $(X,\Delta)$ vanishes,
it follows that $\mathcal{E}(u\cdot g)=0$ for all $g\in\G,$ using
Lemma \ref{lem:Fut}. Hence, using the inequality \ref{eq:comparison E and d one},
it is equivalent to show that the function $g\mapsto d_{1}(u\cdot g,0)$
is proper. Having already shown that $g\mapsto E(g^{*}\omega_{0}^{n}/V)$
is proper --- which is equivalent to $d_{1}(0\cdot g,0)$ being proper
--- it will thus be enough to show that $d_{1}(u\cdot g,0)-d_{1}(0\cdot g,0)$
is bounded by a constant independent of $g.$ But, applying the triangle
inequality to $d_{1}$ yields 
\[
\left|d_{1}(u\cdot g,0)-d_{1}(0\cdot g,0)\right|\leq d_{1}(u\cdot g,0\cdot g)=d(u,0),
\]
 using, in the equality, that $\G$ preserves the metric $d_{1}$
\cite[Lemma 5.9]{d-r}.
\end{proof}
By the previous proposition, the following map is well-defined, if
the Futaki character of $(X,\Delta)$ vanishes:
\[
\pi_{0}\colon\,\mathcal{P}(X)\cap\{E<\infty\}\rightarrow\mathcal{P}(X)_{0}\cap\{E<\infty\}/\K,\,\,\,E\left(\pi_{0}(\nu)\right)=\inf_{\G\nu}E.
\]

Additionally, it restricts to a well-defined map from the space of
probability volume forms $\mathcal{V}(X)$ to $\mathcal{V}(X)_{0}/\K$
where $\mathcal{V}(X)_{0}\coloneqq\mathcal{V}(X)\cap\mathcal{P}(X)_{0}$.
\begin{lem}
\label{lem:mon of pi zero}Assume that the Futaki character of $(X,\Delta)$
vanishes. Then $F_{-\gamma}\left(\pi_{0}(\nu)\right)\leq F_{-\gamma}\left(\nu\right)$
when $\gamma\leq1$ and $F_{-\gamma}\left(\pi_{0}(\nu)\right)\geq F_{-\gamma}\left(\nu\right)$
when $\gamma\geq1.$ As a consequence, for $\gamma\leq1,$
\begin{equation}
\inf_{\mathcal{P}(X)_{0}}F_{-\gamma}=\inf_{\mathcal{P}(X)}F_{-\gamma}=\inf_{\mathcal{V}(X)}F_{-\gamma}=\inf_{\mathcal{V}(X)_{0}}F_{-\gamma}.\label{eq:inf F gamma is inf over moment zero}
\end{equation}
\end{lem}

\begin{proof}
First note that $F_{-\gamma}$ is $\K$-invariant, since $\phi_{0}$
is $\K$-invariant. Moreover, by definition, $F_{-\gamma}=(1-\gamma)E+F_{-1}.$
Since $F_{-1}$ is $\G$-invariant, under the assumption about the
vanishing of the Futaki character, the first statement of the lemma
thus follows and, as a consequence, also the first identity in formula
\ref{eq:inf F gamma is inf over moment zero}. By restricting to volume
forms, also the last identity follows. Thus, to prove the remaining
identity in \ref{eq:inf F gamma is inf over moment zero}, it will
be enough to show that for any given $\mu$ in $\mathcal{P}(X)$ such
that $F_{-\gamma}(\mu)<\infty$, there exists $\nu_{j}\in\mathcal{V}(X)$
such that $F_{-\gamma}(\nu_{j})\rightarrow F_{-\gamma}(\mu).$ But
this follows from \cite[Lemma 3.1]{bdl1}, implying that for any given
$\mu$ such that $E(\mu)<\infty$ and $D(\mu|\mu_{0})<\infty,$ there
exists a sequence $\mu_{j}$ converging towards $\mu$ in $\mathcal{P}(X)$
such that $E(\mu_{j})$ and $D(\mu_{j}|\mu_{0})$ converge towards
$E(\mu)$ and $D(\mu|\mu_{0}),$ respectively. 
\end{proof}
The following result relates the analytic threshold $\delta_{\text{PS}}^{\mathrm{A}}(X,\Delta),$
defined by formula \ref{eq:def of anal st thr}, to the free energy
functional (Section \ref{subsec:The-Mabuchi-f vs free energ}) restricted
to the space $\mathcal{P}(X)_{0}$ of probability measures with vanishing
moment:
\begin{prop}
\label{prop:anal stab in terms of F }The invariant $\delta_{\text{PS}}^{\mathrm{A}}(X,\Delta)$
is independent of the choice of reference $\phi_{0}$ and if the Futaki
character of $(X,\Delta)$ vanishes, then 
\[
\delta_{\text{PS}}^{\mathrm{A}}(X,\Delta)=\sup\left\{ \gamma>0\colon\inf_{\mu\in\mathcal{V}(X)_{0}}F_{-\gamma}(\mu)>-\infty\right\} .
\]
 
\end{prop}

\begin{proof}
To prove the first statement of the proposition it is enough to show
that when changing the metric $\phi_{0}$ the corresponding pluricomplex
energy $E_{\phi_{0}}$ only changes by an additive bounded term. To
this end note that, by formula \ref{eq:diff of E}, given $\phi_{0}$
and $\psi_{0}$ in $\mathcal{H}(L)$ the differential $d(E_{\phi_{0}}-E_{\phi_{1}})|_{\mu}$
is represented by the function $u\coloneqq\phi_{0}-\phi_{1},$ for
any $\mu$ in $\mathcal{P}(X)$ satisfying $E(\mu)<\infty.$ Hence,
there exists a constant $C$ such that
\[
E_{\phi_{0}}(\mu)-E_{\phi_{1}}(\mu)=C+\left\langle u,\mu\right\rangle ,
\]
 showing that $\left|E_{\phi_{0}}(\mu)-E_{\phi_{1}}(\mu)\right|\leq C+\left\Vert u\right\Vert _{L^{\infty}},$
as desired. Finally, the second statement follows from combining formulae
\ref{eq:M as free energ} and \ref{eq:E as I minus J} with proposition
\ref{lem:inf of E over G} and using that $F_{-1}$ is $\G$-invariant
(Section \ref{subsec:The-Futaki-invariants}). 
\end{proof}
We are now ready for the proof of the following result, where the
equivalence $1\iff2$ follows directly from \cite{bbegz,d-r}. The
equivalence $1\iff4$ provides partial evidence for the Large Deviation
Principle (LDP) in Conjecture A, stated in Section \ref{subsec:Motivation:-A-conjectural},
as it establishes the lower semicontinuity of the rate functional
--- a property built into the definition of a LDP.
\begin{thm}
\label{thm:KE equiv}On a log Fano manifold $(X,\Delta)$ the following
are equivalent:
\begin{enumerate}
\item $(X,\Delta)$ admits a Kähler--Einstein metric.
\item $\delta_{\text{PS}}^{\mathrm{A}}(X,\Delta)>1$. 
\item The Futaki character of $(X,\Delta)$ vanishes, $\G$ is reductive
and there exists $\epsilon>0$ such that $F_{-(1+\epsilon)}$ is bounded
from below on $\mathcal{P}(X)_{0},$ i.e. $F(-\gamma)_{0}>-\infty$
for some $\gamma>1.$
\item The Futaki character of $(X,\Delta)$ vanishes, $\G$ is reductive
and the functional $F_{-1}$ is lsc on $\mathcal{P}(X)_{0}$.
\end{enumerate}
Additionally, if $(X,\Delta)$ admits a Kähler--Einstein metric then
there is a unique one whose normalized volume form $\mu_{\text{KE}}$
is in $\mathcal{P}(X)_{0}$. Moreover, the measure $\mu_{\text{KE}}$
is the unique minimizer of $F_{-1}$ on $\mathcal{P}(X)_{0}$ and
\begin{equation}
\mu_{\text{KE }}=\lim_{\beta\searrow-1}\mu_{\beta}\,\,\text{in \ensuremath{\mathcal{P}(X)}},\label{eq:mu KE as limit}
\end{equation}
 where $\mu_{\beta}$ is the unique minimizer of $F_{\beta}$ on $\mathcal{P}(X),$
or equivalently, on $\mathcal{P}(X)_{0}.$ Moreover, $\mu_{\mathrm{KE}}$
is $\K$-invariant.
\end{thm}

\begin{proof}
In the proof we will identify a metric $\phi\in\mathcal{H}(-K_{(X,\Delta)})$
with the corresponding function $u(\coloneqq\phi-\phi_{0})$ in $\mathcal{H}(X,\omega_{0})$
and denote by $u\cdot g$ the function corresponding to $g^{*}\phi.$
As shown in \cite{d-r}, $(X,\Delta)$ admits a Kähler--Einstein
metric iff there exist $\epsilon,C>0$ such that 
\begin{equation}
\mathcal{M}(u)\geq\epsilon\inf_{g\in\G}d_{1}(0,u\cdot g)-C\,\,\,\text{on \ensuremath{\mathcal{H}(X,\omega_{0}}})\cap\{\mathcal{E}(u)=0\},\label{eq:Mab geq d one in pf}
\end{equation}
where $d_{1}$ is the $L^{1}$--Finsler metric on $\mathcal{H}(X,\omega_{0})$
(see Section \ref{subsec:The-space-of geod}). Using the inequality
\ref{eq:comparison E and d one}, this proves $1\iff2$. Next, note
that, $3\implies2$ follows from the previous proposition. To prove
that $1\implies3$ we recall that it was pointed out in \cite[Remark 3.9]{d-r}
that, by density, the inequality \ref{eq:Mab geq d one in pf} is
equivalent to the corresponding inequality on the metric completion
of the metric space $\left(\mathcal{H}(X,\omega_{0}),d_{1}\right),$
which, by \cite{da}, may be identified with $\mathcal{E}^{1}(X,\omega_{0}).$
Combined with the inequalities \ref{eq:comparison E and d one}, this
shows that, assuming $1,$ there exist $\epsilon,C>0$ such that
\[
F(\mu)\geq\epsilon\inf_{\G}E(\mu)-C\text{\,\,on \ensuremath{\mathcal{P}(X)\cap\{E<\infty\}.} }
\]
 In particular, restricting this inequality to $\mathcal{P}(X)_{0},$
and using Prop \ref{lem:inf of E over G}, yields $F(\mu)\geq\epsilon E(\mu)-C$
on $\mathcal{P}(X)_{0}\cap\{E<\infty\},$ i.e. $F_{-(1+\epsilon)}\geq-C$
on $\mathcal{P}(X)_{0}\cap\{E<\infty\}$ and hence, by definition,
on $\mathcal{P}(X)_{0},$ as desired. To prove the remaining equivalence
$1\iff4,$ first assume that $F_{-1}$ is lsc on $\mathcal{P}(X)_{0}.$
Since $\mathcal{P}(X)_{0}$ is compact it follows that $F_{-1}$ admits
a minimizer in $\mathcal{P}(X)_{0}.$ Such a minimizer is KE, since
the Futaki character vanishes (see below). Conversely, assume that
$(X,\Delta)$ admits a KE metric. Take a sequence $\mu_{j}$ converging
towards $\mu_{\infty}$ in $\mathcal{P}(X)_{0}.$ To prove that $F_{-1}$
is lsc we may as well assume that $F_{-1}(\mu_{j})\leq C.$ Using
$1\iff3,$ it follows that $D(\mu_{j}|\mu_{\phi_{0}})\leq C.$ Thus,
by \cite[Thm 2.17]{bbegz}, $F_{-1}(\mu_{\infty})\leq\liminf F_{-1}(\mu_{j}),$
as desired. 

To prove the last statement of the theorem, recall that any KE volume
form $\mu$ in $\mathcal{P}(X)$ minimizes $F_{-1}$ on $\mathcal{P}(X)$
\cite{bbegz}. Thus, by Lemma \ref{lem:mon of pi zero}, there exists
a KE volume form $\pi_{0}(\mu)$ minimizing $F_{-1}$ on $\mathcal{P}(X)_{0}$
(but, a priori, $\pi_{0}(\mu)$ is only determined modulo $\K).$
To prove that there exists a unique KE volume form $\mu_{\text{KE}}$
in $\mathcal{P}(X)_{0},$ we recall that, as shown in \cite[Lemma 6.3]{berm15},
for any given $\phi_{0}\in\mathcal{H}(-K_{(X,\Delta)})$ there exists
a unique KE volume form $\mu_{-1}$ in $\mathcal{P}(X)$ minimizing
$E$ on the space of KE volume forms in $\mathcal{P}(X),$ i.e. on
$\G\mu_{-1}.$ Moreover, $\mu_{-1}$ arises as the limit \ref{eq:mu KE as limit}.
In the present case, where $\phi_{0}$ is assumed $\K$-invariant,
it thus follows from Prop \ref{lem:inf of E over G} that $\mu_{-1}\in\mathcal{P}(X)_{0}$
and conversely, any KE volume form in $\mathcal{P}(X)_{0}$, coincides
with $\mu_{-1}.$ It also follows from the aforementioned uniqueness
that $\mu_{\text{KE}}$ is $\K$-invariant.
\end{proof}
\begin{rem}
When $X$ admits a KE metric and $\G$ is nontrivial it is well known
that $F(\beta)>-\infty$ iff $\beta\geq-1.$ But, $F(\beta)_{0}$
provides a \emph{finite} extension of $F(\beta)$ from $[-1,\infty[$
to some neighbourhood of $[-1,\infty[$ in $\R,$ as follows from
the third point in the previous theorem, combined with Lemma \ref{lem:mon of pi zero}.
\end{rem}

\subsubsection{Comparison with Tian's orthogonality condition}

Consider the space $\mathcal{H}(X,\omega_{0})/\R$ of all Kähler potentials,
modulo additive constants. Assume that $\omega_{0}$ is Kähler--Einstein
and consider the following orthogonality condition on $u\in\mathcal{H}(X,\omega_{0})/\R$:
\begin{equation}
\int_{X}u\psi\omega_{0}^{n}=0\,\,\,\forall\psi\in\text{\ensuremath{\ker}(}\Delta-\lambda_{1})\label{eq:Tians OG cond}
\end{equation}
 where $\text{\ensuremath{\Delta} }$ denotes the Laplacian with respect
to $\omega$ and $\lambda_{1}$ denotes its smallest strictly positive
eigenvalue. This condition appears in Tian's conjecture \cite[Conj 5.5]{ti},
settled in \cite{d-r}, saying that if $\omega_{0}$ is Kähler--Einstein,
then there exists positive constants $\epsilon$ and $C$ such that
\begin{equation}
\mathcal{M}(u)\geq\epsilon(\mathcal{I}-\mathcal{J})(u)-C\label{eq:Tians conj for M}
\end{equation}
 for all $u\in\mathcal{H}(X,\omega_{0})/\R$ satisfying the orthogonality
condition \ref{eq:Tians OG cond}.\footnote{Strictly speaking, in the original conjecture \cite[Conj 5.5]{ti}
the role of $\mathcal{M}$ is played by the Ding functional $\mathcal{D}$
(see Section \ref{subsec:The-constrainted-Ding}), still, both versions
of the conjecture were established in \cite{d-r}. Alternatively,
by Prop \ref{prop:M G-coerc iff D G-coerc}, the two coercivity properties
are equivalent. } It should be emphasized that the functional $(\mathcal{I}-\mathcal{J})$
appearing in the right hand side above may be replaced by the functionals
$\mathcal{I}$ or $\mathcal{J}$ (up to changing the constant $\epsilon),$
since all three functionals $\mathcal{I},(\mathcal{I}-\mathcal{J})$
and $\mathcal{J}$ are comparable. The inequality \ref{eq:Tians conj for M}
should be compared to item 3 in Theorem \ref{thm:KE equiv}, saying
that the inequality \ref{eq:Tians conj for M} holds for all $u\in\mathcal{H}(X,\omega_{0})/\R$
such that $\omega_{u}^{n}\boldsymbol{m}$ integrates to zero on $X.$
Note that, while Tian's orthogonality condition \ref{eq:Tians OG cond}
is linear in $u$, the moment condition appearing in Theorem \ref{thm:KE equiv}
is linear with respect to the Monge--Ampère measure $\omega_{u}^{n}.$
The relation between these two conditions is clarified in the following
lemma, which gives a refinement of \cite[Lemma A.2]{z-z}:
\begin{lem}
\label{lem:inf J}Let $X$ be a Fano manifold and assume that the
Futaki character of $X$ vanishes and that the group $\G$ is reductive.
Fix a Kähler form $\omega_{0}$ that is invariant under a maximal
compact subgroup $\K$ of $\G.$ Given $u\in\mathcal{E}^{1}(X,\omega_{0})/\R$
the following conditions are equivalent for $v\in u\cdot\G\subset\mathcal{E}^{1}(X,\omega_{0})/\R$:
\[
(i)\,\inf_{g\in\G}\mathcal{J}(u\cdot g)=\mathcal{J}(v),\,\,\,(ii)\int_{X}v\Delta\boldsymbol{m}\omega_{0}^{n}=0,\,\,\,(iii)\int_{X}\boldsymbol{m}\omega_{v}\wedge\omega_{0}^{n-1}=0
\]
 where $\Delta$ denotes the Laplacian with respect to the Kähler
metric $\omega_{0}.$ Moreover, modulo the action of $\K,$ there
exists a unique $v\in u\cdot\G\subset\mathcal{E}^{1}(X,\omega_{0})/\R$
satisfying any of these three conditions. When $\omega_{0}$ is Kähler--Einstein
the condition $(ii)$ is equivalent to Tian's orthogonality condition
\ref{eq:Tians OG cond} for $v.$
\end{lem}

\begin{proof}
First observe that $\mathcal{J}(u)$ is $\K$-invariant, since $\omega_{0}$
is assumed $\K$-invariant. It thus follows, precisely as in the proof
of Prop \ref{lem:inf of E over G}, that the inf over $\G$ in $(i)$
may be replaced by the inf over $J\mathfrak{k},$ under the identification
of $\K\backslash\G$ with $J\mathfrak{k}.$ Now fix a nontrivial real
vector field $W$ on $X$ induced by an element in $\mathfrak{k}\subset\mathfrak{g}.$
Let us first show that 
\begin{equation}
V\frac{d\mathcal{J}(v\cdot\exp(-tJW))}{dt}\Bigg|_{t=0}=c_{n}\int_{X}v\Delta h_{W}\omega_{0}^{n}\label{eq:formul for deriv av J along W}
\end{equation}
 for a positive constant $c_{n}$ only depending on $n.$ It is enough
to consider the case when $v\in\mathcal{H}(X,-K_{X}),$ using an approximation
argument, precisely as in the proof of Prop \ref{lem:inf of E over G}.
Setting $v_{t}\coloneqq v\cdot\exp(-tJW),$ 
\begin{equation}
V\frac{d\mathcal{J}(v_{t})}{dt}_{|t=0}=\int_{X}\left.\frac{dv_{t}}{dt}\right|_{t=0}\omega_{0}^{n}=-\int_{X}\left\langle dv,JW\right\rangle \omega_{0}^{n}=-\int_{X}\nabla v\cdot JW\omega_{0}^{n}\label{eq:pf of Lemma Tian}
\end{equation}
using, in the first equality, that the Futaki character vanishes (formula
\ref{eq:def of Fut}). Expressing $\left\langle dv,JW\right\rangle $
as the scalar product $\nabla v\cdot JW$ between the vector field
$JW$ and the gradient $\nabla v$ of $v$ with respect to the Kähler
metric $\omega_{0}$ and using that $JW=c_{n}\nabla h_{W}$ for a
positive constant $c_{n}$ (where $h_{W}$ denotes the Hamiltonian
function corresponding to $W)$ we deduce that 
\[
V\frac{d\mathcal{J}(v_{t})}{dt}_{|t=0}=-c_{n}\int_{X}\nabla v\cdot\nabla h_{W}\omega_{0}^{n}=c_{n}\int_{X}v\Delta h_{W}\omega_{0}^{n},
\]
In particular, $\frac{d}{dt}\mathcal{J}(v_{t})=0$ for all $W\in\mathfrak{k}$
iff condition $(ii)$ in the present lemma holds. To prove that $(i)$
is equivalent to $(ii),$ together with the uniqueness statement in
the lemma, it is thus enough to show that $\mathcal{J}(v_{t})$ is
strictly convex for any given $W\in\mathfrak{k}.$ To this end we
first note that formula \ref{eq:formul for deriv av J along W} can
be reformulated as
\[
V\frac{d\mathcal{J}(v\cdot\exp(-tJW))}{dt}_{|t=0}=c_{n}\int_{X}h_{W}(\Delta v+n)\omega_{0}^{n-1}=\frac{c_{n}}{(n-1)!}\int_{X}h_{W}\omega_{v}\wedge\omega_{0}^{n-1},
\]
 showing that condition $(ii)$ is also equivalent to condition $(iii).$
Indeed, integrating by parts, we can replace the integrand $v\Delta h_{W}$
in formula \ref{eq:formul for deriv av J along W} with $\Delta vh_{W}$
and then use that the integral of $h_{W}\omega_{0}^{n}$ vanishes,
by Lemma \ref{lem:Fut}, since the Futaki character of $X$ is assumed
to vanish. It now follows that the function $t\mapsto\mathcal{J}(v\cdot\exp(-tJW))$
is strictly convex, precisely as in the proof of Prop \ref{lem:inf of E over G},
using that the measure $\omega_{v}\wedge\omega_{0}^{n-1}$ does not
charge pluripolar subsets. Indeed, since $v\in\mathcal{E}^{1}(X,\omega)$
the measure defined as the product $\omega_{v}\wedge\omega_{0}^{n-1}$
of the current $\omega_{v}$ with the smooth form $\omega_{0}^{n-1}$
is identical to the non-pluripolar product of $\omega_{v}$ with $\omega_{0}^{n-1},$
which, by definition, does not charge pluripolar subsets (this identity
follows, for example, from writing $v$ as a decreasing limit of elements
in $\mathcal{E}^{1}(X,\omega)$ and using that both types of products
are continuous wrt decreasing limits, since $v$ has full Monge--Ampère
mass \cite[Prop 2.1]{bbgz}). Next note that the pullback of $\mathcal{J}$
to $J\mathfrak{k}$ is proper and thus admits a critical point. Indeed,
since $\mathcal{J}$ is comparable to $\mathcal{I}-\mathcal{J},$
which can be identified with the functional $E$ on $\mathcal{P}(X),$
the properness of $\mathcal{J}$ follows directly from the properness
of $E,$ established in the course of the proof of Prop \ref{lem:inf of E over G}.
Finally, when $\omega_{0}$ is Kähler--Einstein, $\Delta h_{W}=\lambda_{1}h_{W}$
for any $W\in\mathfrak{k},$ where $\lambda_{1}$ denotes the smallest
strictly positive eigenvalue of $\Delta$ \cite[Thm 2.2]{b-m}. Since
$\boldsymbol{m}=(h_{W_{1}},\dots,h_{W_{r}})$ for a basis $W_{i}$
in $\mathfrak{k}$, it follows that $(ii)$ is equivalent to Tian's
condition \ref{eq:Tians OG cond}, when $\omega_{0}$ is Kähler--Einstein.
\end{proof}

\subsubsection{Comparison with the $\G$-reduced weighted delta invariant of Delcroix--Jubert}

In \cite{de-ju}, the following analytic coercivity threshold for
a Fano manifold $X$ was introduced as a special case of a general
weighted and twisted setting:
\begin{equation}
\delta^{\G}(X)\coloneqq\sup\{\delta\in\mathbb{R}|\exists C_{\delta}\in\mathbb{R}\text{ s.t. }\forall u\in\mathcal{H}(X,\omega_{0}),D(\omega_{u}^{n}/V|\mu_{0})-\delta\inf_{g\in\G}E(g_{*}\omega_{u}^{n}/V)\geq C_{\delta}\}.\label{eq: def of reduced delta}
\end{equation}
The definition is reminiscent of the analytic polystability threshold
and, indeed, they are closely related. 
\begin{prop}
Let $X$ be a Fano manifold, then $\delta^{\G}(X)\geq\delta_{\mathrm{PS}}^{\mathrm{A}}(X)$.
If furthermore the Futaki character of $X$ vanishes, then $\delta_{\mathrm{PS}}^{\mathrm{A}}(X)=\mathrm{\delta^{\G}(X)}$.\label{prop: relation with reduced delta}
\end{prop}

\begin{proof}
The first statement follows from applying the trivial inequality $E(\omega_{u}^{n}/V)\geq\inf_{g\in\G}E(g_{*}\omega_{u}^{n}/V)$
in the definition of $\delta_{\mathrm{PS}}^{\mathrm{A}}(X)$. For
the second statement, in the definition of $\delta^{\G}(X)$, restrict
to $u\in\mathcal{H}(X,\omega_{0})$ such that $\omega_{u}^{n}/V\in\mathcal{P}(X)_{0}$
and apply Propositions \ref{lem:inf of E over G} and \ref{prop:anal stab in terms of F },
which assumes that the Futaki character of $X$ vanishes, to find
that $\delta^{\G}(X)\leq\delta_{\mathrm{PS}}^{\mathrm{A}}(X)$. 
\end{proof}

\subsection{\label{subsec:The-constrainted-Ding}The constrained Ding functional
vs the penalized Ding functional }

Given a log Fano manifold $(X,\Delta),$ a reference metric $\phi_{0}\in\mathcal{H}(-K_{(X,\Delta)})$
and $\gamma\in\R$ the \emph{twisted Ding functional} on $\mathcal{H}(X,\omega_{0})$,
is defined by
\begin{equation}
\mathcal{D}_{\gamma}(u)\coloneqq-\mathcal{E}(u)-\frac{1}{\gamma}\log\int_{X}e^{-\gamma u}\mu_{0}.\label{eq:def of D gamma}
\end{equation}
 It was introduced by Ding in \cite{din} (when $\Delta=0)$\footnote{In the notation of \cite{din}, $\mathcal{D}_{\gamma}(u)=F_{-\gamma}(u),$
but here the notation $F_{-\gamma}$ is reserved for the free energy
functional.}. When $\gamma=1$ we will simply refer to $\mathcal{D}_{\gamma}$
as the \emph{Ding functional}, denoted by $\mathcal{D},$ whose critical
points are the potentials of KE metrics. Moreover, 
\begin{equation}
\gamma^{-1}F_{-\gamma}(\mu)\geq\mathcal{D}_{\gamma}(u_{\mu}),\,\,\,\,\,\inf_{\mu\in\mathcal{P}(X)}\gamma^{-1}F_{-\gamma}(\mu)=\inf_{u\in\mathcal{E}^{1}(X)}\mathcal{D}_{\gamma}(u).\label{eq:F gamma geq D gamma}
\end{equation}
The first inequality follows directly from the definitions (see the
proof of Lemma \ref{lem:constrained F geq penalized Ding} below)
and the second equality is shown in \cite[Thm 1.1]{berm6} using a
duality argument involving the Legendre--Fenchel transform wrt the
standard duality pairing between measures $\mu$ and functions $u$
on $X.$ It is also shown in \cite{berm6} that $\gamma^{-1}F_{-\gamma}$
is coercive iff $\mathcal{D}_{\gamma}$ is coercive. 

Now, in order to take the action of the group $\G$ into account we
assume, as before, that $\phi_{0}$ is $\K$-invariant for a maximal
compact subgroup $\K.$ A key role in the present work (and in the
companion paper \cite{a-b-s3}) will be played by a penalized version
of the twisted Ding functional, defined by $\mathcal{D}_{\gamma}(u)+\delta_{\gamma}(u),$
where $\delta_{\gamma}(u)$ is the following ``entropic distance'':
\begin{equation}
\delta_{\gamma}(u)\coloneqq\frac{1}{\gamma}D(\mathcal{P}(X)_{0}|\mu_{\gamma u}),\,\,\,\,D(\mathcal{P}(X)_{0}|\nu)\coloneqq\inf_{\mu\in\mathcal{P}(X)_{0}}D(\mu|\nu),\label{eq:def of entropic distance}
\end{equation}
between the measure $\mu_{u\gamma}\coloneqq e^{-\gamma u}\mu_{0}/\int_{X}e^{-\gamma u}\mu_{0}\in\mathcal{P}(X)$
and the convex subspace $\mathcal{P}(X)_{0}.$ Adding $\delta_{\gamma}(u)$
to $\mathcal{D}_{\gamma}(u)$ can be viewed as a softening of the
constraint that $\mu_{\gamma u}\in\mathcal{P}(X)_{0},$ i.e. 
\begin{equation}
\int_{X}\boldsymbol{m}\mu_{\gamma u}=0.\label{eq:Aubin constraint with gamma}
\end{equation}

The definition of $\delta_{\gamma}(u)$ is motivated by the following
\begin{lem}
\label{lem:constrained F geq penalized Ding}For any $\gamma>0$ and
$\mu\in\mathcal{P}(X)_{0}$
\[
\gamma^{-1}F_{-\gamma}(\mu)\geq\mathcal{D}_{\gamma}(u_{\mu})+\delta_{\gamma}(u_{\mu})
\]
\end{lem}

\begin{proof}
It follows directly from the definitions that $\gamma^{-1}F_{-\gamma}(\mu)=\mathcal{D}_{\gamma}(u_{\mu})+\gamma^{-1}D(\mu|\mu_{u\gamma})$
and that $D(\mu|\nu)\geq D(\mathcal{P}(X)_{0}|\nu)$ when $\mu\in\mathcal{P}(X)_{0}.$
\end{proof}
Furthermore, will show in the companion paper \cite{a-b-s3} that
the infimum over $\mathcal{P}(X)_{0}$ of the lhs in the inequality
in the previous lemma coincides with the infimum of the rhs over all
$u\in\mathcal{E}^{1}(X),$ by generalizing the duality argument in
\cite{berm6} (but this result will not be needed in the present paper).
In the case when $X=\P^{1}$ the constraint \ref{eq:Aubin constraint with gamma}
on $u$ appears in Aubin's strengthening of the Moser--Trudinger
inequality on the two-sphere \cite[Cor 2]{A}, which shows that ---
under this constraint --- $\mathcal{D}_{\gamma}(u)$ is bounded from
below for any $\gamma\in]0,2[.$ Here we obtain a new statistical
mechanical proof of Aubin's result by combining the effective large
deviation bound in Theorem \ref{thm:effective LD bound with Ding}
with Theorem \ref{thm:gibbs poly and lct intro} (and another analytic
proof is given in the companion paper \cite{a-b-s3} that applies
to any log Fano curve). Moreover, from the stronger large deviation
bound in Theorem \ref{thm:effective LD bound with F} it follows that
$F_{-\gamma}(\mu)$ is bounded from below on $\mathcal{P}(\P^{1})_{0}$
when $\gamma\in]0,2[,$ which implies Aubin's result, by the following
\begin{lem}
\label{lem:F gamma bounded below implies D}Given $\gamma>0$, assume
that $F_{-\gamma}(\mu)$ is bounded from below on $\mathcal{P}(X)_{0}.$
Then $\mathcal{D}_{\gamma}(u)$ is bounded from below for all $u\in\mathcal{E}^{1}(X,\omega_{0})$
satisfying the constraint \ref{eq:Aubin constraint with gamma}, i.e.
$\mu_{\gamma u}\in\mathcal{P}(X)_{0}.$ Moreover, when $n=1$, the
space $\mathcal{E}^{1}(X,\omega_{0})$ may be replaced by the Sobolev
space $W^{2,1}(X).$
\end{lem}

\begin{proof}
The first statement follows directly from the inequality $\frac{1}{\gamma}F_{\gamma}\left(\frac{e^{-\gamma u}\mu_{0}}{\int_{X}e^{-\gamma u}\mu_{0}}\right)\leq\mathcal{D}_{\gamma}(u)$
appearing in \cite[formula 3.6]{berm6}. The refinement when $n=1$
is shown in a similar way, using that, in this case, the sup in the
defining formula \ref{eq:def of E as sup} for $E(\mu)$ can be taken
over $W^{2,1}(X)$ (a detailed argument appears in the companion paper
\cite{a-b-s3}).
\end{proof}
It is natural to ask whether the converse of the preceeding lemma
holds. However, as discussed in the companion paper \cite{a-b-s3},
such a converse appears unlikely in general. A weaker converse result
is a corollary of the following
\begin{prop}
\label{prop:constrained Ding vs properness}Given $\gamma\geq1,$
assume that there exist constants $\sigma$ and $C$ such that 
\[
\mathcal{D}_{\gamma}(u)\geq-\sigma\mathcal{J}(u)+C
\]
 for any $u\in\mathcal{H}(X,\omega_{0})$ satisfying the constraint
\ref{eq:Aubin constraint with gamma}. Then 
\begin{equation}
\mathcal{D}(v)\geq\left(1-\frac{(1+\sigma)}{\gamma^{1/n}}\right)\mathcal{J}(v)+\gamma C\label{eq:ineq on D v in pf}
\end{equation}
 for any $v\in\mathcal{H}(X,\omega_{0})$ satisfying the constraint
$\int_{X}\boldsymbol{m}\mu_{v}=0$ and, as a consequence, for any
$v\in\mathcal{H}(X,\omega_{0}),$ 
\begin{equation}
\mathcal{D}(v)\geq\left(1-\frac{(1+\sigma)}{\gamma^{1/n}}\right)\inf_{g\in\G}\mathcal{J}(g\cdot v)+\gamma C\label{eq:ineq on D v with inf in pf}
\end{equation}
if the Futaki character of $(X,\Delta)$ vanishes. The same statements
hold when $\mathcal{H}(X,\omega_{0})$ is replaced by $\mathcal{E}^{1}(X,\omega_{0}).$ 
\end{prop}

\begin{proof}
Given $v\in\mathcal{H}(X,\omega_{0})$ satisfying $\int_{X}\boldsymbol{m}\mu_{v}=0$,
set $u=v/\gamma$, which is in $\mathcal{H}(X,\omega_{0})$ since
$\gamma\geq1,$ and satisifies $\int_{X}\boldsymbol{m}\mu_{\gamma u}=0.$
Since the inequality \ref{eq:ineq on D v in pf} is invariant when
$v$ is shifted by a constant we may as well assume that $\int v\omega_{0}^{n}=0.$
Then, $J=-\mathcal{E}$ and the assumed inequality for $u$ thus gives
\[
\gamma J(\gamma^{-1}v)(1+\sigma)-\log\int e^{-v}\mu_{0}\geq\gamma C.
\]
 Finally, since $J(tv)\leq t^{1+1/n}J(v)$ for any $v\in\mathcal{H}(X,\omega_{0})$
and $t\in[0,1]$ (see \cite[Remark 2]{din}) we have $\gamma J(\gamma^{-1}v)\leq J(v)/\gamma^{1/n},$
which shows that the inequality \ref{eq:ineq on D v in pf} indeed
holds for any $v\in\mathcal{H}(X,\omega_{0})$ satisfying the constraint
$\int_{X}\boldsymbol{m}\mu_{v}=0.$ All that remains it to verify
the claim that, given $v\in\mathcal{E}^{1}(X,\omega_{0})$ there exists
$g\in\G$ such that $\int_{X}\boldsymbol{m}\mu_{v\cdot g}=0.$ Indeed,
since the Futaki character of $(X,\Delta)$ is assumed to vanish,
$\mathcal{D}$ is $\G$-invariant and we can thus apply the inequality
\ref{eq:ineq on D v in pf} to $v\cdot g$ and take the inf over all
$g\in\G$ to deduce the inequality \ref{eq:ineq on D v with inf in pf}.
Finally, to prove the claim we identify $v\in\mathcal{E}^{1}(X,\omega_{0})$
with the corresponding metric $\phi\coloneqq\phi_{0}+v\in\text{PSH}(-K_{(X,\Delta)}).$
Then $\mu_{v\cdot g}=g^{*}\mu_{\phi}$ and the existence of the desired
$g\in\G$ thus follows from Prop \ref{lem:inf of E over G}.
\end{proof}
\begin{cor}
Assume that the Futaki character of $(X,\Delta)$ vanishes, $\gamma\geq1$
and that $\mathcal{D}_{\gamma}(u)$ is bounded from below on all $u\in\mathcal{E}^{1}(X,\omega_{0})$
satisfying the constraint \ref{eq:Aubin constraint with gamma}, i.e.
$\mu_{\gamma u}\in\mathcal{P}(X)_{0}.$ Then $F_{-(1+n^{-1}(1-\gamma^{-1/n}))}(\mu)$
is bounded from below on $\mathcal{P}(X)_{0}.$ 
\end{cor}

\begin{proof}
This follows directly from combining the previous proposition with
the inequality $\mathcal{J}\geq n^{-1}(\mathcal{I}-\mathcal{J})$
and Prop \ref{lem:inf of E over G}. 
\end{proof}
Combining the previous results yields the following proposition, where
$(\mathcal{I}-\mathcal{J})^{\G}\coloneqq\inf_{g\in\G}(\mathcal{I}-\mathcal{J}).$ 
\begin{prop}
\label{prop:M G-coerc iff D G-coerc}Assume that the Futaki character
of $(X,\Delta)$ vanishes and fix $\epsilon>0.$ If $\mathcal{D}\geq\epsilon(\mathcal{I}-\mathcal{J})^{\G}+C$
on $\mathcal{E}^{1}(X,\omega_{0}),$ then $\mathcal{M}\geq\epsilon(\mathcal{I}-\mathcal{J})^{\G}+C$
on $\mathcal{E}^{1}(X,\omega_{0})$. Conversely, if the latter inequality
holds, then $\mathcal{D}\geq n^{-1}\left(1-(1+\epsilon)^{-1/n}\right)(\mathcal{I}-\mathcal{J})^{\G}+C.$
The same statements hold when $\mathcal{E}^{1}(X,\omega_{0})$ is
replaced by $\mathcal{H}(X,\omega_{0})$. 
\end{prop}

\begin{proof}
The first statement follows directly from applying the inequality
in formula \ref{eq:F gamma geq D gamma} with $\gamma=1.$ To prove
the converse statement, note that the assumed inequality for $\mathcal{M}$
is equivalent, by Prop \ref{lem:inf of E over G}, to the inequality
$F_{-\gamma}\geq C$ on $\mathcal{P}(X)_{0}$ for $\gamma=1+\epsilon.$
This implies, by the proof of Lemma \ref{lem:F gamma bounded below implies D},
that $\mathcal{D}_{\gamma}(u)\geq C/\gamma$ for any $u\in\mathcal{E}^{1}(X,\omega_{0})$
satisfying the constraint \ref{eq:Aubin constraint with gamma}. The
proof is thus concluded by invoking Prop \ref{prop:constrained Ding vs properness}
and using that $\mathcal{J}\geq n^{-1}(\mathcal{I}-\mathcal{J}).$ 
\end{proof}

\section{\label{sec:The-reduced-probability}The $\K$-reduced probability
measures on $X^{N}$}

We continue with the setup introduced in Section \ref{sec:Symmetry-breaking,-thermodynamic}.
Given a log Fano manifold $(X,\Delta),$ we thus fix a compact subgroup
$\K$ of $\G$ and a $\K$-invariant positively curved metric on $-K_{(X,\Delta)}.$
Recall that $\boldsymbol{m}_{N}$ denotes the corresponding moment
map for the diagonal action of $\K$ on $X^{N}$, taking values in
the dual $\mathfrak{k}^{*}$ of the Lie algebra $\mathfrak{k}$ of
$\K$ (Section \ref{subsec:The-moment-map}). 

Given a positive rational number $k$ such that $-kK_{(X,\Delta)}$
defines a line bundle, we denote by $N$ the dimension of $H^{0}(X,-kK_{(X,\Delta)})$
and by $\nu^{(N)}$ the measure on $X^{N}$ defined by 
\[
\nu^{(N)}=\frac{i^{(Nn)^{2}}}{2^{Nn}}\left(\det S^{(N)}\otimes S_{\Delta_{N}}\right)^{-1/k}\wedge\overline{\left(\det S^{(N)}\otimes S_{\Delta_{N}}\right)}^{-1/k},
\]
 where $\det S^{(N)}$ denotes the section of $-kK_{(X^{N},\Delta_{N})}\rightarrow X^{N}$
obtained by taking $s_{1},\dots,s_{N}$ in formula \ref{eq:slater determinant intro}
to be a basis of $H^{0}(X,-kK_{(X,\Delta)})$ which is orthonormal
with respect to the Hermitian product on $H^{0}(X,-kK_{(X,\Delta)})$
induced by the metric $\phi_{0}$ on $-K_{(X,\Delta)}$ and the measure
$\mu_{\phi_{0}}$ on $X.$ 
\begin{lem}
\label{lem:character}Given a log Fano manifold $(X,\Delta),$ there
exists a real character $\chi_{N}$ on $\G$ such that $g^{*}\nu_{N}=\chi_{N}(g)\nu_{N}$
for any $g\in\G.$ More precisely, given $W$ in the Lie algebra of
$\G$, 
\begin{equation}
\chi_{N}(e^{tW})=e^{tN2\mathrm{Re}\left(\mathrm{Fut}_{(X,\Delta),k}(W)\right)},\label{eq:def of real character}
\end{equation}
 for $t\in\R,$ where $\mathrm{Fut}_{(X,\Delta),k}(W)$ denotes the
quantized Futaki invariant at level $k.$ Moreover, if $k$ is a positive
integer and $k\Delta$ is a divisor, then $\mathrm{Fut}_{(X,\Delta),k}=0$
iff $\det S^{(N)}$ is $\G$-invariant.
\end{lem}

\begin{proof}
First assume that $k$ is a positive integer and $k\Delta$ is a divisor.
Then $\G$ acts naturally on $(X,-kK_{(X,\Delta)})$ and, as a consequence,
also on $H^{0}(X,-kK_{(X,\Delta)}).$ By basic linear algebra, $g^{*}\nu_{N}=\left|\mathrm{det}A(g)\right|^{-2/k}\nu_{N}$
where $A(\cdot)$ denotes the representation of $\G$ on $H^{0}(X,-kK_{(X,\Delta)})$
and $\mathrm{det}A(g)$ defines a $\C^{*}$-valued algebraic character
on $\G.$ Hence, $\chi_{N}(g)\coloneqq\left|\mathrm{det}A(g)\right|^{-2/k}$
is a real-valued character on $\G$ and $\chi_{N}(e^{tW})=e^{tN2\text{Re}\left(\text{Fut}_{(X,\Delta),k}(W)\right)}$
by the very definition of $\text{Fut}_{(X,\Delta),k}(W).$ In particular,
if $\text{Fut}_{(X,\Delta),k}=0$ then $\mathrm{det}A(g)$ is a holomorphic
function on $\G$ taking values in $S^{1}.$ But this implies that
$\mathrm{det}A(g)$ is identically equal to $1,$ showing that $\det S^{(N)}$
is $\G$-invariant. Finally, consider the case when it is only assumed
that $-kK_{(X,\Delta)}$ defines a line bundle. Then $\G$ still acts
naturally on the real line bundle over $X^{N}$ defined as the top
exterior power of the real cotangent bundle of $X^{N},$ having $\left|\det S^{(N)}\right|^{-2/k}$
as a smooth section. We can thus proceed essentially as before. 
\end{proof}
Define the\emph{ energy per $N$ particles} as the following lsc function
on $X^{N}$
\[
E^{(N)}(x_{1},\dots,x_{N})\coloneqq-(kN)^{-1}\log\left(\left\Vert \det S^{(N)}(x_{1},\dots,x_{N})\right\Vert _{\phi_{0}}^{2}\right).
\]
By \ref{eq:support of anti-canon divisor}, 
\begin{equation}
E^{(N)}(x_{1},\dots,x_{N})=\infty\iff\exists s\in H^{0}(X,-kK_{(X,\Delta)})-\{0\}\text{ such that }s(x_{1})=\cdots=s(x_{N})=0.\label{eq:E N infinity}
\end{equation}

As shown in \cite{berm8}, $E^{(N)}$ converges towards the pluricomplex
energy $E$ on $\mathcal{P}(X),$ as $N\rightarrow\infty,$ in the
sense of $\Gamma$-convergence. The first part of the following result
can thus be viewed as an analog of Prop \ref{lem:inf of E over G}
for finite $N$ and the second part is an analog of the fact that
the unique minimizer $\omega_{0}^{n}/V$ of the functional $E$ on
$\mathcal{P}(X)$ is in $\mathcal{P}(X)_{0},$ if the Futaki character
of $(X,\Delta)$ vanishes (by Lemma \ref{lem:Fut}).
\begin{prop}
Let $(X,\Delta)$ be a log Fano manifold. Assume that the quantized
Futaki character of $(X,\Delta)$ at level $k$ vanishes. Then, under
the diagonal action of $\G$ of $X^{N},$ 
\[
E^{(N)}(x_{1},\dots,x_{N})<\infty\implies(x_{1},\dots,x_{N})\,\text{is semi-stable. }
\]
As a consequence, the closure of the $\G$-orbit of $(x_{1},\dots,x_{N})$
intersects $\left\{ \boldsymbol{m}_{N}=0\right\} .$ Moreover, any
configuration $(x_{1},\dots,x_{N})_{*}$ minimizing $E^{(N)}$ satisfies
$\boldsymbol{m}_{N}((x_{1},\dots,x_{N})_{*})=0.$ Additionally, if
$H^{0}(X,-kK_{(X,\Delta)}\otimes T^{*}X)\neq\{0\}$, then $E^{(N)}(x_{1},\dots,x_{N})<\infty$
implies that $(x_{1},\dots,x_{N})$ has a finite stabilizer. In particular,
then $\boldsymbol{m}_{N}$ is a submersion at $(x_{1},\dots,x_{N})_{*}.$
\end{prop}

\begin{proof}
First, assume that $k$ is a positive integer such that $k\Delta$
is a divisor. Then $\G$ naturally acts on $-kK_{(X,\Delta)}$ and
the section $\det S^{(N)}$ is $\mathbb{G}$-invariant, by Lemma \ref{lem:character}.
Hence, if $\det S^{(N)}(x_{1},\dots,x_{N})\neq0,$ then $(x_{1},\dots,x_{N})$
is semi-stable, by the very definition of semi-stability. Next, if
$\boldsymbol{x}\coloneqq(x_{1},\dots,x_{N})_{*}$ minimizes $E^{(N)}$
then it follows, in particular, that $\boldsymbol{x}$ minimizes the
restriction of $E^{(N)}$ to the orbit $\G\boldsymbol{x}.$ Note that,
since $\mathrm{det}S^{(N)}$ is $\mathbb{G}$-invariant, the restriction
$f$ of $E^{(N)}$ to $\G\boldsymbol{x}$ defines a Kempf--Ness functional,
i.e. up to scaling,
\[
f(g\boldsymbol{x})=-\log\left\Vert gS_{\boldsymbol{x}}\right\Vert (g\boldsymbol{x})
\]
 for a non-zero vector $S_{\boldsymbol{x}}$ in the complex line over
$\boldsymbol{x},$ which in this case is $(-kK_{X^{N_{k}}})|_{\boldsymbol{x}}.$
But in general, a point $\boldsymbol{x}$ is a minimum of the Kempf--Ness
function on $\G\boldsymbol{x}$ iff the corresponding moment map vanishes
at $\boldsymbol{x}$ \cite{ho}. In the general case when it is only
assumed that $-kK_{(X,\Delta)}$ is defined as a line bundle, one
can replace $\det S^{(N)}$ with a tensor power $(\det S^{(N)})^{\otimes m}$
for $m$ a sufficiently divisible positive integer and proceed essentially
as before. Finally, assume that $H^{0}(X,-kK_{(X,\Delta)}\otimes T^{*}X)\neq\{0\}.$
Assume that $E^{(N)}(x_{1},\dots,x_{N})\neq\infty$ and --- in order
to get a contradiction --- that the stabilizer of $(x_{1},\dots,x_{N})$
is infinite. Then there exists a nontrivial holomorphic vector field
$V$ on $X$ vanishing at all points $x_{1},\dots,x_{N}.$ Take an
element $S_{k}$ in $H^{0}(X,-kK_{(X,\Delta)}\otimes T^{*}X)-\{0\}.$
Contracting $S_{k}$ with $V$ then gives a nontrivial holomorphic
section of $-kK_{(X,\Delta)}$ vanishing at all points $x_{1},\dots,x_{N},$
which contradicts $E^{(N)}(x_{1},\dots,x_{N})\neq\infty.$ Finally,
we recall that, in general, the moment map is a submersion at $\boldsymbol{x}$
iff the stabilizer is finite iff $\boldsymbol{x}$ is stable.
\end{proof}
Note that the assumption $H^{0}(X,-kK_{(X,\Delta)}\otimes T^{*}X)\neq\{0\}$
is satisfied when $k$ is sufficiently large, since $-K_{(X,\Delta)}$
is ample. In general, this assumption is needed to get finite stabilizers,
as illustrated in the following example. 
\begin{example}
Let $(X,\Delta)=(\P^{1},0).$ In this case $k\in\Z_{+}/2,$ $N=2k+1=2,3,\dots$
and $H^{0}(X,-kK_{(X,\Delta)}\otimes T^{*}X)\neq\{0\}$ holds iff
$k\geq1,$ i.e. $N\geq3.$ For any $N\geq2$, $E^{(N)}(x_{1},\dots,x_{N})\neq0$
iff all $x_{i}$ are distinct. Hence, $\left\{ E^{(N)}(x_{1},\dots,x_{N})\neq\infty\right\} $
is a subset of $(X^{N})_{\mathrm{ss}}$ , which is a proper subset
when $N>3$ (by the condition \ref{eq:SS cond}). Consider now the
case when $k=1/2$ (i.e. $N=2)$ and take $\phi_{0}$ to be the Fubini--Study
metric i.e. $\phi_{0}$ is an $\mathrm{SO}(3)$-invariant metric.
Identifying $\P^{1}$ with the unit-sphere in $\R^{3},$ $\left\Vert \det S^{(2)}(x_{1},x_{2})\right\Vert _{\phi_{0}}$
is proportional to $|x_{1}-x_{2}|$ and $\boldsymbol{m}_{2}=(x_{1}+x_{2})/2.$
Hence, $(x_{1},x_{2})$ minimizes $E^{(2)}(x_{1},x_{2})=0$ iff $x_{1}$
and $x_{2}$ are anti-podal points, which is equivalent to $\boldsymbol{m}_{2}(x_{1},x_{2})=0.$
This is also a consequence of the previous proposition. However, in
this case the stabilizer is infinite (and, as a consequence, $\boldsymbol{m}_{2}$
is not a submersion at $(x_{1},x_{2})$). Indeed, any two distinct
points on $\P^{1}$ are fixed by a $\C^{*}-$action. Note also that
$N=1$ on $\P^{1}$ is not covered by the previous proposition, which
is consistent with the fact that there are no semi-stable points in
this case (since the condition \ref{eq:SS cond} is never satisfied).
In particular, the moment map has no zeroes. 
\end{example}

\subsection{Gelfand--Leray measures and contracted measures}

Before introducing the probabilistic setup we start with some general
definitions. Assume that $0$ is a regular value of the moment map
$\boldsymbol{m}_{N}.$ Given a measure $\nu$ on $X^{N}$ with a smooth
density and a basis $e_{1},\dots,e_{r}$ of $\mathfrak{k}$, we will
denote by $\nu_{\boldsymbol{m}_{N}}$ the corresponding Gelfand--Leray
form on the submanifold $\left\{ \boldsymbol{m}_{N}=0\right\} $ of
$X^{N},$ uniquely determined by the following property: 
\begin{equation}
\nu_{\boldsymbol{m}_{N}}\wedge dh_{1}\wedge\cdots\wedge dh_{r}=\nu\label{eq:gelf-ler prop}
\end{equation}
 when evaluated at any given point in $\left\{ \boldsymbol{m}_{N}=0\right\} $.
Here $h_{1},\dots,h_{r}$ are Hamiltonians for the vector fields on
$X$ corresponding to $e_{1},\dots,e_{r}$ with respect to the action
of $\K$. Note that $\mathbf{m}_{N}=\sum_{i=1}^{r}h_{i}e_{i}^{*}$.
In particular, $\nu_{\boldsymbol{m}_{N}}$ defines a measure on $X^{N},$
supported on $\left\{ \boldsymbol{m}_{N}=0\right\} $, that we shall
call the\emph{ Gelfand--Leray measure.} 

\subsection{\label{subsec:Probabilistic-setup}Probabilistic and statistical
mechanical setup}

\subsubsection{\label{subsec:The-Gelfand-Leray-measure}The Gelfand--Leray measure
$\nu_{0}^{(N)}$ and the $\K$-reduced canonical probability measure
$\mu_{0}^{(N)}$ and partition function $Z_{N,0}.$ }

Assume that $0$ is a \emph{regular value} of $\boldsymbol{m}_{N}.$
We then denote by $\nu_{0}^{(N)}$ the Gelfand--Leray measure supported
on $\left\{ \boldsymbol{m}_{N}=0\right\} \subset X^{N},$ corresponding
to the measure $\nu^{(N)}$ on $X^{N}.$ While this measure depends
on the choice of basis in $\mathfrak{k},$ changing the basis only
has the effect of multiplying $\nu_{0}^{(N)}$ with a positive constant
$C$ (the determinant of the corresponding change of bases matrix).
Since $C$ is independent of $N$ this does not effect the large-$N$
limits in Conjecture \ref{conj:large N limtis intro} and Theorem
\ref{thm:conv towards inf for n one intro}. Set
\begin{equation}
\mu_{0}^{(N)}\coloneqq\nu_{0}^{(N)}/Z_{N,0},\,\,\,Z_{N,0}\coloneqq\int\nu_{0}^{(N)}.\label{eq:def of K-red can prob measure and partition}
\end{equation}
 Assuming that $Z_{N,0}<\infty$ the measure $\mu_{0}^{(N)}$ defines
a probability measure on $X^{N},$ that we shall call the\emph{ $\K$-reduced
canonical probability measure }and the corresponding normalizing constant
$Z_{N,0}$ will be called the\emph{ $\K$-reduced partition function.}

\subsubsection{The thickened setup}

When $0$ is not a regular value of $\boldsymbol{m}_{N}$ we replace
the measure $\nu_{0}^{(N)}$ with 
\begin{equation}
\nu_{\epsilon_{N}}^{(N)}\coloneqq1_{\left\{ ||\boldsymbol{m}_{N}||_{\infty}<\epsilon_{N}\right\} }\nu^{(N)}\left/\epsilon_{N}^{\mathrm{dim}\mathfrak{k}}\right.\label{eq:def of reg G-L when sing}
\end{equation}
 for an appropriate sequence of positive numbers $\epsilon_{N},$
called the \emph{thickening sequence}. The corresponding normalization
constants $Z_{N,\epsilon_{N}}$ and probability measures $\mu_{\epsilon_{N}}^{(N)}$
will be called the \emph{thickened }partition functions and\emph{
}probability measures.
\begin{rem}
When $0$ is a regular value of $\boldsymbol{m}_{N}$, $\nu_{0}^{(N)}$
coincides with the limit of $\nu_{\epsilon_{N}}^{(N)}$ as $\epsilon_{N}\rightarrow0$
(for $N$ fixed). It would be interesting to know if, when $\{\boldsymbol{m}_{N}=0\}$
is non-empty, this limit exists and if it coincides with the measure
supported on the regular locus of $\{\boldsymbol{m}_{N}=0\},$ defined
by the Gelfand--Leray measure on the regular locus.
\end{rem}

\subsubsection{The $\K$-reduced Gibbs measures at inverse temperature $\beta$ }

Given a number $\beta\in\R,$ set 
\[
\nu_{\beta}^{(N)}\coloneqq e^{-\beta NE^{(N)}}\mu_{0}^{\otimes N}.
\]
When $\beta=-1$, $\nu_{\beta}^{(N)}$ coincides with the measure
$\nu^{(N)}$ on $X^{N}$ (but $\nu_{\beta}^{(N)}$ depends on the
fixed metric $\phi_{0}$ when $\beta\neq-1).$ The measures $\nu_{\beta}^{(N)}$
on $X^{N}$ induce a family of Gelfand--Leray measures $\left(\nu_{\beta}^{(N)}\right)_{\boldsymbol{m}_{N}}.$
Assuming that the \emph{$\K$-reduced partition function $Z_{N,0}(\beta)$}
\emph{at inverse temperature} $\beta$ is finite,
\[
Z_{N,0}(\beta)\coloneqq\int\left(\nu_{\beta}^{(N)}\right)_{\boldsymbol{m}_{N}}<\infty,
\]
 we obtain a measure on $X^{N}$ that we shall call the \emph{$\K$-reduced
Gibbs measures at inverse temperature $\beta$ }by setting 
\begin{equation}
\mu_{0,\beta}^{(N)}\coloneqq\frac{\left(\nu_{\beta}^{(N)}\right)_{\boldsymbol{m}_{N}}}{Z_{N,0}(\beta)}=\frac{e^{-\beta NE^{(N)}}\left(\mu_{0}^{\otimes N}\right)_{\boldsymbol{m}_{N}}}{Z_{N,0}(\beta)}.\label{eq:Gibbs measure general text}
\end{equation}
When $\beta=-1$ this measure coincides with the $\K$-reduced canonical
probability measure $\mu_{0}^{(N)}$, defined in the previous section,
and when $\beta=0$ it coincides with the classical microcanonical
measure associated to the Hamiltonians $h_{1},\dots,h_{r},$ that
is ubiquitous in the statistical mechanics literature (especially
when $r=1)$. In general, $Z_{N,0}(\beta)$ is finite when $\beta>-\alpha(X,\Delta)$
(by the inequality \ref{eq:lower bound in terms of alpha}). But if
Conjecture \ref{conj:lct is an stab intro} is correct, then this
bound can be improved considerably. 
\begin{rem}
\label{rem:vortex}From a statistical mechanics perspective, the function
$E^{(N)}$ plays the role of the $N$-particle energy (or Hamiltonian)
and the corresponding Gibbs measure $\mu_{0,\beta}^{(N)}$ describes
the corresponding equilibrium state of $N$ interacting particles
confined to $\left\{ \boldsymbol{m}_{N}=0\right\} ,$ at inverse temperature
$\beta.$ The corresponding probability space $(X^{N},\mu_{0,\beta}^{(N)})$
may be viewed as a singular incarnation of the notion of a \emph{mixed
ensemble} introduced in \cite{e-h-t}. In particular, when $X=\P^{1},$
identified with the two-sphere, and $\omega_{0}$ is the standard
$\mathrm{SO}(3)$-invariant symplectic form on $X,$ $E^{(N)}$ coincides
with the energy particle in Onsager's statistical mechanical description
of the 2D incompressible Euler equations, adapted to the two-sphere
\cite{ki2,k-w}. The corresponding Gibbs measures $\mu_{0,\beta}^{(N)}$
are invariant under the Hamiltonian flow of $E^{(N)}$ on $X^{N}$
and, building on Theorem \ref{thm:conv towards inf for n one intro},
we show in the companion paper \cite{a-b-s1} that the empirical measure
$\delta_{N},$ viewed as a random variable on the probability space
$(X^{N},\mu_{0,\beta}^{(N)})$ converges in law towards the $\mathrm{SO}(3)$-invariant
probability measure on the two-sphere, for any $\beta>-1.$ 
\end{rem}

\begin{prop}
For any $\beta\in\R$, the probability measure $\mu_{0,\beta}^{(N)}$
is $\K$-invariant, when it is well-defined.
\end{prop}

\begin{proof}
First observe that for any $g\in\G,$ $g^{*}\det S^{(N)}=\chi_{N}(g)\det S^{(N)}$
for a homomorphism $\chi_{N}\colon\G\rightarrow\C^{*}$ (see Lemma
\ref{lem:character} below). It follows that the restriction of $\chi_{N}$
to $\K$ takes values in $S^{1}\subset\C^{*}.$ Indeed, since $\K$
is compact, the image $\chi_{N}(\G)$ is a compact subgroup of $\C^{*}$
and thus contained in the unique maximal compact subgroup $S^{1}$
of $\C^{*}.$ Since $\phi_{0}$ is $\K$-invariant, it follows that
$\left\Vert \det S^{(N)}(x_{1},\dots,x_{N})\right\Vert _{\phi_{0}}^{2}$
is $\K$-invariant on $X^{N}.$ Finally, since the measure $\mu_{0}$
and $|\boldsymbol{m}|$ are also $\K$-invariant, we deduce that $\mu_{0,\beta}^{(N)}$
is $\K$-invariant. 
\end{proof}

\subsection{\label{subsec:Motivation:-A-conjectural}Motivation: A conjectural
Large Deviation Principle}

We next formulate a conjecture, which serves as motivation for the
conjectures stated in the introduction of the paper. Recall that a
sequence of measures $\Gamma_{N}$ on a topological space $\mathcal{X}$
satisfies a \emph{Large Deviation Principle (LDP)} with \emph{speed}
$N$ and \emph{rate function} $F\colon\mathcal{X}\rightarrow]-\infty,\infty]$
if $F$ is lsc and for any (Borel) measurable subset $\mathcal{B}$
of $\mathcal{X}$,
\[
-\inf_{\mu\in\mathcal{\mathring{B}}}F\leq\liminf_{N\rightarrow\infty}\frac{1}{N}\log\Gamma_{N}(\mathcal{B})\leq\limsup_{N\rightarrow\infty}\frac{1}{N}\log\Gamma_{N}(\mathcal{B})\leq-\inf_{\mu\in\mathcal{\mathcal{\overline{B}}}}F,
\]
where $\mathcal{\mathring{B}}$ and $\mathcal{\mathcal{\overline{B}}}$
denote the interior and closure of $\mathcal{B},$ respectively (see
\cite{d-z} for background on large deviation theory). An LDP is symbolically
expressed as $\Gamma_{N}\sim e^{-NF}.$ 
\begin{conjecture}
\label{conj:LDP} The following conjectures hold when $N$ is a sequence
of regular values:
\begin{itemize}
\item \textbf{A(i)} If $F(\beta-\epsilon)_{0}>-\infty,$ for some $\epsilon>0,$
then the measures $(\delta_{N})_{*}\left((\nu_{\beta}^{(N)})_{\boldsymbol{m}_{N}}\right)$
on $\mathcal{P}(X)_{0}$ satisfy an LDP with speed $N$ and rate functional
$F_{\beta}(\mu)$. \textbf{(ii)} If $\mu_{\beta}^{(N)}$ are well-defined
probability measures for $N\gg1,$ then the probability measures $(\delta_{N})_{*}\text{\ensuremath{\mu_{\beta}^{(N)}}}$
on $\mathcal{P}(X)_{0}$ satisfy an LDP with speed $N$ and rate functional
$F_{\beta}(\mu)-C$ for some $C\in\R.$ 
\item \textbf{B(i)} If $F(\beta-\epsilon)>-\infty,$ for some $\epsilon>0,$
then the measures $(\delta_{N})_{*}\left(\nu_{\beta}^{(N)}\right)$
on $\mathcal{P}(X)$ satisfy an LDP with speed $N$ and rate functional
$F_{\beta}(\mu)$. \textbf{(ii}) if $\nu_{\beta}^{(N)}/\int\nu_{\beta}^{(N)}$
are well-defined probability measures for $N\gg1,$ then the probability
measures $(\delta_{N})_{*}\text{\ensuremath{\left(\nu_{\beta}^{(N)}/\int\nu_{\beta}^{(N)}\right)}}$
on $\mathcal{P}(X)$ satisfy an LDP with speed $N$ and rate functional
$F_{\beta}(\mu)-C$ for some $C\in\R.$ 
\end{itemize}
Without assuming that $N$ is a sequence of regular values, the corresponding
LDPs also hold on $\mathcal{P}(X)$ when $\left\{ \boldsymbol{m}_{N}=0\right\} $
is replaced by $\left\{ |\boldsymbol{m}_{N}|\leq\epsilon_{N}\right\} $
for any appropriate sequence $\epsilon_{N},$ with rate functional
$F_{\beta}+\chi_{0},$ where $\chi_{0}=0$ on $\mathcal{P}(X)_{0}$
and $\chi_{0}=\infty$ away from $\mathcal{P}(X)_{0}.$ 
\end{conjecture}

Conjecture B follows from results in \cite{berm12}, when $n=1$,
and when $\beta>0$ it is established in \cite{berm8,berm8 comma 5}
for any $n.$ In the general case, Conjecture A is motivated by the
following heuristic argument. First, as shown in \cite{berm8}, the
following mean field type approximation holds:

\[
E^{(N)}(x_{1},\dots,x_{N})\rightarrow E(\mu),\,\,\text{as \ensuremath{\delta_{N}(x_{1},\dots,x_{N})\rightarrow\mu}}
\]
 in the sense of $\Gamma$-convergence. This suggests the formal approximation
\[
(\delta_{N})_{*}\text{\ensuremath{\nu_{\beta}^{(N)}}}\sim e^{-N\left(\beta E+\chi_{0}\right)}(\delta_{N})_{*}\left(\mu_{0}^{\otimes N}\right),
\]
as probability measures on $\mathcal{P}(X),$ where $\chi_{0}=0$
on $\mathcal{P}(X)_{0}$ and $\chi_{0}=\infty$ away from $\mathcal{P}(X)_{0}.$
Next, according to Sanov's classical LDP \cite{d-z},
\[
(\delta_{N})_{*}\left(\mu_{0}^{\otimes N}\right)\sim e^{-ND(\mu|\mu_{0})}.
\]
Hence, formally, $(\delta_{N})_{*}\text{\ensuremath{\nu_{\beta}^{(N)}}}$
satisfies an LDP with rate functional $F_{\beta}+\chi_{0}.$ Finally,
since $(\delta_{N})_{*}\text{\ensuremath{\nu_{\beta}^{(N)}}}$ is
supported on $\mathcal{P}(X)_{0}$ such an LDP is equivalent to the
LDP for $(\delta_{N})_{*}\text{\ensuremath{\nu_{\beta}^{(N)}}}$ in
Conjecture \ref{conj:LDP}. The argument for the probability measures
$(\delta_{N})_{*}\text{\ensuremath{\mu_{\beta}^{(N)}}}$ is similar. 

We next show that Conjecture \ref{conj:LDP} implies the conjectures
stated in the introduction of the paper.
\begin{prop}
If Conjecture \ref{conj:LDP} A holds for all $\beta<0,$ then 
\[
\gamma_{\mathrm{PS}}(X,\Delta)=\sup_{\gamma>0}\left\{ F(-\gamma)_{0}>-\infty\right\} .
\]
 As a consequence, if Conjecture \ref{conj:LDP} A holds, then $\gamma_{\mathrm{PS}}(X,\Delta)>1$
and the Futaki character of $(X,\Delta)$ vanishes iff $(X,\Delta)$
is K-polystable. Moreover, Conjecture \ref{conj:large N limtis intro}
also follows. Assuming also the validity of Conjecture \ref{conj:LDP}
B(i), it follows that $(X,\Delta)$ is uniformly Gibbs polystable
iff $(X,\Delta)$ is K-polystable. Finally, assuming only the validity
of Conjecture B, 
\[
\gamma(X,\Delta)=\sup_{\gamma>0}\left\{ F(-\gamma)>-\infty\right\} .
\]
\end{prop}

\begin{proof}
Assume that $F(\beta-\epsilon)_{0}>-\infty,$ for some $\epsilon>0.$
Then it follows directly from the definition of the assumed LDP for
$(\delta_{N})_{*}\left((\nu_{\beta}^{(N)})_{\boldsymbol{m}_{N}}\right)$
that 
\begin{equation}
-\lim_{N\rightarrow\infty}N^{-1}\log\int(\nu_{\beta}^{(N)})_{\boldsymbol{m}_{N}}=\inf_{\mathcal{P}(X)_{0}}F_{\beta}(\mu)\eqqcolon F(\beta)_{0}.\label{eq:conv towards F beta in proof}
\end{equation}
 As a consequence, $\int(\nu_{\beta}^{(N)})_{\boldsymbol{m}_{N}}<\infty,$
for $N\gg1,$ showing that 
\begin{equation}
\gamma_{\text{PS}}(X,\Delta)\geq\sup_{\gamma>0}\left\{ F(-\gamma)_{0}>-\infty\right\} .\label{eq:ineq in conj implies conj}
\end{equation}
 Conversely, take $\gamma>\gamma_{\text{PS}}(X,\Delta).$ Then $\mu_{-\gamma}^{(N)}$
is a well-defined probability measure for $N\gg1,$ as follows from
the argument in the proof of Prop \ref{prop:Gibbs poly iff integrable}.
Hence, the assumed LDP for the probability measures $(\delta_{N})_{*}\text{\ensuremath{\mu_{\beta}^{(N)}}}$
yields
\[
0=-\lim_{N\rightarrow\infty}N^{-1}\log\int_{X^{N}}\mu_{\beta}^{(N)}=\inf_{\mathcal{P}(X)_{0}}F_{\beta}(\mu)-C,
\]
 showing that $F(\beta)_{0}=C>-\infty.$ Thus, the inequality \ref{eq:ineq in conj implies conj}
is an equality. As a consequence, $\gamma_{\text{PS}}(X,\Delta)>1$
iff there exists $\gamma>1$ such that $F(-\gamma)_{0}>-\infty.$
By Theorem \ref{thm:KE equiv}, assuming that the Futaki character
of $(X,\Delta)$ vanishes, this means that $\gamma_{\text{PS}}(X,\Delta)>1$
iff $\delta_{\text{PS}}^{\mathrm{A}}(X,\Delta)>1$ iff $(X,\Delta)$
admits a KE metric, which by the solution of the YTD conjecture equivalently
means that $(X,\Delta)$ is K-polystable (see \cite{c-d-s} for the
case $\Delta=0$ with $X$ smooth and \cite{li,l-x-z} for the case
of a general log Fano variety $(X,\Delta)).$ Finally, the last statement
of the proposition is shown as above. 

Next, to prove that Conjecture \ref{conj:LDP} A implies Conjecture
\ref{conj:large N limtis intro}, assume that $(X,\Delta)$ admits
a KE metric. It follows, from Theorem \ref{thm:KE equiv} , that $F(\beta+\epsilon)_{0}>-\infty,$
for some $\epsilon>0,$ Hence, the convergence \ref{eq:conv towards F beta in proof}
for $\beta=-1$ implies, by Lemma \ref{lem:mon of pi zero}, the convergence
\ref{eq:conv towards Mab in conj intro} towards $\inf\mathcal{M}$
which is finite (since the inf is attained at any KE metric). As a
consequence, $\int(\nu^{(N)})_{\boldsymbol{m}_{N}}$ is finite for
$N$ sufficiently large. This means that $\mu_{-1}^{(N)}(=\mu_{0}^{(N)})$
is a well-defined probability measure for $N\gg1.$ Thus, by assumption,
$(\delta_{N})_{*}\text{\ensuremath{\mu_{-1}^{(N)}} }$ satisfies an
LDP with rate functional $F_{-1}(\mu)-C$ on $\mathcal{P}(X)_{0}.$
By standard results for LDPs, to prove the convergence towards $\mu_{\text{KE}}$
in \ref{conj:Gibbs poly intro} it is thus enough to show that $\mu_{\text{KE}}$
is the unique minimizer of the rate functional. But this is shown
in Theorem \ref{thm:KE equiv}.

Now assume that $(X,\Delta)$ is uniformly Gibbs polystable, i.e.
$(X,\Delta)$ is Gibbs semistable and $\gamma_{\text{PS}}(X,\Delta)>1$.
Then, by \cite{f-o}, $(X,\Delta)$ is K-semistable, which implies
that the Futaki character of $(X,\Delta)$ vanishes. Hence, as shown
above, $(X,\Delta)$ is K-polystable. Conversely, if $(X,\Delta)$
is K-polystable, then it admits a KE metric and, as a consequence,
the inf of $F_{-1}(\mu)$ over $\mathcal{P}(X)$ is finite. Assuming
that Conjecture \ref{conj:LDP} B holds it follows, as before, that,
for $\beta>-1,$ 
\begin{equation}
-\lim_{N\rightarrow\infty}N^{-1}\log\int_{X^{N}}\nu_{\beta}^{(N)}=\inf_{\mathcal{P}(X)}F_{\beta}(\mu),\,\,\,\gamma(X,\Delta)\geq\sup_{\gamma>0}\left\{ F(-\gamma)>-\infty\right\} .\label{eq:conv towars inf F on all of P in pf}
\end{equation}
 Hence, for any $\beta>-1,$ $\int_{X^{N}}\nu_{\beta}^{(N)}<\infty$
for $N$ sufficiently large, which means that $(X,\Delta)$ is Gibbs
semistable. Since, as shown above, Conjecture \ref{conj:LDP} A implies
that $\gamma_{\text{PS}}(X,\Delta)>1,$ it follows that $(X,\Delta)$
is uniformly Gibbs polystable. 
\end{proof}
\begin{rem}
\label{rem:ldp implies Gibbs semi iff Fut zero}If Conjecture \ref{conj:LDP}
holds, then $(X,\Delta)$ is Gibbs polystable iff its Futaki character
vanish and $\gamma_{\text{PS}}(X,\Delta)>1,$ by the previous proposition.
\end{rem}

\subsubsection{A digression on the need for $\epsilon>0$ in Conjecture \ref{conj:LDP}}

It should be stressed that, in general, it is not enough to assume
that Conjecture B holds with $\epsilon=0$ (and, most likely, not
in Conjecture A either). Indeed, if Conjecture B holds with $\epsilon=0,$
then it would imply the implication 
\[
F_{-1}\,\,\text{bounded below on \ensuremath{\mathcal{P}(X)\implies F_{-1}\,\,\text{is lsc on \ensuremath{\mathcal{P}(X)}}}}.
\]
 But this implication fails in general. Indeed, by \cite[Thm 5.11]{berm11},
$F_{-1}$ is lsc on $\mathcal{P}(X)$ iff $X$ admits a KE metric
and $\G$ is trivial (which, in turn, is equivalent to $X$ being
K-stable \cite{c-d-s}). Hence, already $X=\P^{1}$ is a counterexample.
But there are also counterexamples when $\G$ is trivial. Indeed,
$F_{-1}$ is bounded from below on $\mathcal{P}(X)$ iff $(X,\Delta)$
is K-semistable (see \cite{li00} for the case when $X$ is smooth
and $\Delta$ is trivial and the proof of \cite[Thm 2.5]{a-b1} for
the general case). However, there are examples of K-semistable log
Fano manifolds $(X,\Delta)$ with $\G$ trivial, such that $(X,\Delta)$
is not K-stable (for example, this can be achieved when $X=\P^{1}$
and $\Delta$ consists of three points). 

\subsection{Gibbs polystability vs finiteness of the restricted partition function
$Z_{N,0}$}

In this section we show that Gibbs polystability arises naturally
from the construction of the $\K$-reduced canonical probability measure. 
\begin{lem}
\label{lem:integrab in neigh of m zero}$\gamma_{\mathrm{PS}}^{(N)}(X,\Delta)>1$
iff $\nu^{(N)}\in L_{\text{}}^{1}\left(\left\{ \left\Vert \boldsymbol{m}_{N}\right\Vert <\epsilon_{N}\right\} \right)$
for any $\epsilon_{N}$ sufficiently small to ensure that $\left\{ \left\Vert \boldsymbol{m}_{N}\right\Vert <\epsilon_{N}\right\} \subset(X^{N})_{\mathrm{ss}}.$
\end{lem}

\begin{proof}
First, by the analytic definition of the lct (formula \ref{eq:lct anal}),
$\gamma_{\text{PS}}^{(N)}(X,\Delta)>1$ iff $\nu^{(N)}\in L_{\text{loc}}^{1}((X^{N})_{\text{ss }}).$
Next, note that there exists $\epsilon_{N}\in]0,\infty[$ such that
the closure of $\left\{ \left\Vert \boldsymbol{m}_{N}\right\Vert <\epsilon_{N}\right\} $
is contained in $(X^{N})_{\text{ss}}$ (since $\left\{ \boldsymbol{m}_{N}=0\right\} $
is compact and contained in the open set $(X^{N})_{\text{ss}}).$
Hence, $\gamma_{\text{PS}}^{(N)}(X,\Delta)>1$ implies that $\nu^{(N)}\in L_{\text{}}^{1}\left(\left\{ \left\Vert \boldsymbol{m}_{N}\right\Vert <\epsilon_{N}\right\} \right).$
To prove the converse, assume that $\nu^{(N)}\in L_{\text{}}^{1}\left(\left\{ \left\Vert \boldsymbol{m}_{N}\right\Vert <\epsilon_{N}\right\} \right)$
and consider the integral of $\nu^{(N)}$ over a neighbourhood $U$
of a given point $(x_{1},\dots,x_{N})$ in $(X^{N})_{\text{ss }}.$
By the Kempf--Ness theorem \cite[Thm 3.3]{ho}, a point $(x_{1},\dots,x_{N})\in(X^{N})_{\text{ss}}$
iff $\overline{\G(x_{1},\dots,x_{N})}$ intersects $\left\{ \left\Vert \boldsymbol{m}_{N}\right\Vert =0\right\} .$
This means that there exists $g\in\G$ such that $g(x_{1},\dots,x_{N})\in\left\{ \left\Vert \boldsymbol{m}_{N}\right\Vert <\epsilon_{N}\right\} $
and, as a consequence, there exists a neighborhood $U$ of $(x_{1},\dots,x_{N})$
such that $gU\subset\left\{ \left\Vert \boldsymbol{m}_{N}\right\Vert <\epsilon_{N}\right\} .$
Thus, by assumption, $\int_{gU}\nu^{(N)}<\infty$ and since the integrability
of $\nu^{(N)}$ is $\G$-invariant (by the previous lemma) we deduce
that $\int_{U}\nu^{(N)}<\infty,$ showing that $\nu^{(N)}\in L_{\text{loc}}^{1}((X^{N})_{\text{ss }}),$
as desired.
\end{proof}
\begin{prop}
\label{prop:Gibbs poly iff integrable} First assume that $0$ is
a regular value of $\boldsymbol{m}_{N}.$ Then $\gamma_{\mathrm{PS}}^{(N)}(X,\Delta)>1$
iff $Z_{N,0}<\infty.$ In particular, $(X,\Delta)$ is Gibbs polystable
iff it is Gibbs semistable and $Z_{N,0}<\infty$ for all sufficiently
large $N.$ In general, $\gamma_{\mathrm{PS}}^{(N)}(X,\Delta)>1$
iff $Z_{N,\epsilon_{N}}<\infty$ when $\epsilon_{N}$ is sufficiently
small. In particular, $(X,\Delta)$ is Gibbs polystable iff it is
Gibbs semistable and $Z_{N,\epsilon_{N}}<\infty$ when $N$ is sufficiently
large and $\epsilon_{N}$ is sufficiently small. 
\end{prop}

\begin{proof}
When $0$ is not a regular value of $\boldsymbol{m}_{N}$ this follows
directly from the previous lemma. The case when $0$ is a regular
value of $\boldsymbol{m}_{N}$ is established in \cite{a-b-s1} (the
proof appeared in a previous ArXiv version of the present paper). 
\end{proof}

\section{\label{sec:Uniform-Gibbs-polystability}Strong uniform Gibbs polystability
and large deviation bounds}

In this section we will deduce the effective bounds in Theorem \ref{thm:uniform polyGibbs implies KE}
from effective large deviation bounds on thickened partition functions
of independent interest (Theorems \ref{thm:effective LD bound with Ding},
\ref{thm:effective LD bound with F}). The constants in the bounds
will come from upper bounds on Bergman measures established in Section
\ref{subsec:Effective-upper-bounds Berg}. 

We continue with the setup in Section \ref{subsec:Setup analytic},
where $(X,\Delta)$ is a given log Fano manifold and we have fixed
a maximal compact subgroup $\K$ of $\G(\coloneqq\text{Aut}_{0}(X,\Delta))$
and a $\K$-invariant metric $\phi_{0}$ on $-K_{(X,\Delta)}$ with
positive curvature form $\omega_{0}.$ We will denote by $R_{\text{out}}$
the \emph{outer radius} of the image $\boldsymbol{m}(X)$ of the moment
map induced by $(\phi_{0},\K)$ and by $R_{\text{in}}$ the \emph{inner
radius }of the convex hull of $\boldsymbol{m}(X)$ (which is strictly
positive, by Lemma \ref{lem:A N is finite} below). 

We start by recalling some general notions from GIT and defining the
notion of strong uniform Gibbs polystability.

\subsection{Recap of the normalized Hilbert--Mumford weight }

Let $G$ be a complex reductive group acting on a polarized projective
manifold $(X,L)$ and fix a maximal compact subgroup $K$ of $G$
and an $\text{ad}K$-invariant scalar product on the Lie algebra $\mathfrak{k}$
of $K.$ The corresponding norms on $\mathfrak{k}^{*}$ and $\mathfrak{k}$
will both be denoted by $\left|\cdot\right|.$ Denote by $w_{\lambda}(x)$
the \emph{weight} of a point $x\in X$ wrt a nontrivial subgroup $\lambda$
of $G$ isomorphic to $\C^{*},$ which takes values in $\Z$ (also
known as the Hilbert--Mumford weight) and denote by $\left|\lambda\right|$
the norm of a generator of the compact subgroup of $\lambda,$ identified
with an element in the Lie algebra of $K.$ Given a ``thickening
parameter'' $\epsilon\in[0,\infty[$, set 
\[
X_{\text{ss,\ensuremath{\epsilon} }}\coloneqq\left\{ x\in X\colon M(x)\geq-\epsilon\right\} ,\,\,\,M(x)\coloneqq\inf_{\lambda}\frac{w_{\lambda}(x)}{\left|\lambda\right|}.
\]
Thus, $X_{\text{ss},\epsilon}$ coincides with the semi-stable locus
$X_{\text{ss }}$ when $\epsilon=0$ and contains $X_{\text{ss }}$
for any $\epsilon.$ The sets $X_{\text{ss,\ensuremath{\epsilon} }}$
are Zariski open for generic $\epsilon$ (avoiding the finite set
$-M(X)\cap]0,\infty[$). This follows from properties of the \emph{Hesselink
stratification} of $X-X_{\text{ss }}$ (where the function $-M(x)$
defines the order of the different of strata \cite[Section 5.6]{ho2}).

\subsection{Strong uniform Gibbs polystability}

Coming back to our setup where $(X,\Delta)$ is a given log Fano manifold,
consider the action of the group $\G$ on $X^{N}$ and set, for $\epsilon>0$
and $N\epsilon$ generic,
\[
(X^{N})_{\text{ss},\epsilon}\coloneqq\left\{ X^{N}\colon\frac{1}{N}M(x_{1},...,x_{N})\geq-\epsilon\right\} ,\,\,\,\,\gamma_{\epsilon}^{(N)}(X,\Delta)\coloneqq\text{lct}\left((X^{N})_{\text{ss},\epsilon},\Delta_{N};\mathcal{D}^{(N)}\right)
\]
A sequence $\epsilon_{N}$ of such parameters will be called a \emph{thickening}
\emph{sequence }(since $(X^{N})_{\text{ss}}$ is contained in $(X^{N})_{\text{ss},\epsilon}$).
\begin{defn}
\label{def:strong}A log Fano manifold $(X,\Delta)$ is said to be
\emph{strongly uniformly Gibbs polystable} if $(X,\Delta)$ is Gibbs
semi-stable and there exists a sequence $\epsilon_{N}\in[0,\infty[,$
called the\emph{ thickening} \emph{sequence}, such that $\epsilon_{N}>(r+1)R_{\mathrm{out}}/N$
and
\[
\liminf_{N\rightarrow\infty}\gamma_{\epsilon_{N}}^{(N)}(X,\Delta)>1.
\]
In particular, if one instead would take that $\epsilon_{N}=0,$ then
this means that $X$ is uniformly Gibbs polystable. We expect that
a log Fano manifold is uniformly Gibbs polystable iff it is \emph{strongly}
uniformly Gibbs polystable, as illustrated by the following
\end{defn}

\begin{example}
\label{exa:strong curve}A log Fano curve $(\P^{1},\Delta)$ is uniformly
Gibbs polystable iff it is strongly uniformly Gibbs polystable. Indeed,
for \emph{any} sequence $\epsilon_{N}$ tending to zero 
\[
\lim_{N\rightarrow\infty}\gamma_{\epsilon_{N}}^{(N)}(X,\Delta)=\lim_{N\rightarrow\infty}\gamma_{\text{PS}}^{(N)}(X,\Delta)
\]
as follows from a slight generalization of the proof of Theorem \ref{thm:gibbs poly and lct intro}.
Indeed, when, for example $\Delta=0,$ so that $\G=\mathrm{SL}(2,\C),$
we get, for an appropriate choice of norm, $M(x_{1},...,x_{N})/N=1-2I/N$
where $I$ denotes the numbers of $x_{i}$ coinciding (as follows
from the formula for the weights $w_{\lambda}(x_{1},...,x_{N})$ in
\cite[Example 6.6]{ho2}). Hence, $M(x_{1},...,x_{N})/N\geq-\epsilon_{N}$
iff $I/N\leq1/2+\epsilon_{N}/2$ and one can then proceed precisely
as in the proof of Theorem \ref{thm:gibbs poly and lct intro}.
\end{example}

\subsection{\label{subsec:Effective-upper-bounds Berg}Effective upper bounds
on the Bergman measure $\beta_{k\phi}$ and the invariant $\sigma_{k}(X)$}

Given an ample line bundle $L\rightarrow X$ and a psh metric $\phi$
on $L$ such that $e^{-\phi}\in L_{\text{loc}}^{1},$ consider the
\emph{Bergman measure} $\beta_{\phi}$ defined by the following probability
measure on $X$:
\[
\beta_{\phi}\coloneqq\frac{K_{\phi}}{N(L)},\,\,K_{\phi}\coloneqq\sup_{s\in H^{0}(X,L+K_{X})-\{0\}}\frac{s\wedge\bar{s}e^{-\phi}}{\int_{X}s\wedge\bar{s}e^{-\phi}},\,\,\,N(L)\coloneqq\dim H^{0}(X,L+K_{X}),
\]
 where $s$ has been identified with a holomorphic top form on $X$
with values in $L$ so that $i^{n^{2}}s\wedge\bar{s}e^{-\phi}$ defines
a measure on $X.$ In the case when $X$ is a Fano manifold we will,
for a given positive integer $k,$ take $L=(k+1)(-K_{X})$ and consider
the following sequence of invariants of $X:$ 
\[
\sigma_{k}(X)\coloneqq\inf_{\phi\in\mathcal{H}(-K_{X})^{\K}}\sigma_{k}(\phi),\,\,\,\,\sigma_{k}(\phi)\coloneqq\sup_{X}\left(\frac{\beta_{(k+1)\phi}}{(dd^{c}\phi)^{n}/V}-1\right)
\]
The constants in the effective bounds in Theorems \ref{thm:uniform polyGibbs implies KE}
depend on effective bounds on $\sigma_{k}(X)$ (see Cor \ref{cor:crit for KE text}).
The following non-effective asymptotic bound on $\sigma_{k}(X)$ follows
from standard results: 
\begin{lem}
Denote by $S_{\omega}$ the \emph{scalar curvature} of a Kähler metric
$\omega$ on $X,$ normalized so that $S_{\omega}\omega^{n}=\mathrm{Ric}\,\omega\wedge\omega^{n-1}/(n-1)!$
and set $\overline{S_{\omega}}\coloneqq\int_{X}S_{\omega}\omega^{n}$.
As $k\rightarrow\infty$,
\[
\sigma_{k}(X)\leq k^{-1}\frac{1}{2}\inf_{\omega\in c_{1}(X)}\left\Vert S_{\omega}-\overline{S_{\omega}}\right\Vert _{L^{\infty}(X)}+O(k^{-2}).
\]
Moreover, if $X$ is a homogeneous Fano manifold, then $\sigma_{k}(X)=0.$ 
\end{lem}

\begin{proof}
It follows from the well-known asymptotic expansion of the Bergman
kernel on the diagonal for $pL+E,$ that for any given positive curved
line bundle $L$ over $X$ and any holomorphic line bundle $E$ over
$X,$ that, 
\begin{equation}
\frac{K_{p\phi}}{p^{n}(dd^{c}\phi)^{n}/V}=1-p^{-1}\frac{1}{2}S_{\omega}+O(p^{-2}),\,\,\,\omega\coloneqq dd^{c}\phi,\,\,V\coloneqq L^{n},\label{eq:exp of Bermgan kernel}
\end{equation}
 as $p\rightarrow\infty$, in the $C^{0}$-topology on $X$ (e.g.
\cite[Section 2.5]{bbs}). Integrating this formula over $X$ (or
using the Riemann--Roch theorem) gives 
\[
N(pL)/L^{n}p^{n}=1-p^{-1}\frac{1}{2}\overline{S_{\omega}}+O(p^{-2})\,\,\Big(=1+p^{-1}\frac{1}{2}\frac{K_{X}\cdot L^{n-1}}{(n-1)!}+O(p^{-2})\Big).
\]
As a consequence, for any fixed $\phi\in\mathcal{H}(-K_{X})^{\K},$
\[
\sigma_{k}(X)\leq k^{-1}\sup_{X}\left(\frac{\beta_{(k+1)\phi}}{(dd^{c}\phi)^{n}/V}-1\right)\leq k^{-1}\frac{1}{2}\left\Vert S-\bar{S}\right\Vert _{L^{\infty}(X)}+O(k^{-2}).
\]
Finally, when $X$ is a homogeneous Fano manifold any maximal compact
subgroup $\K$ acts transitively on $X.$ As a consequence, for any
$\K$-invariant metric $\phi$ on $-K_{X}$ the two probability measures
$\beta_{(k+1)\phi}$ and $(dd^{c}\phi)^{n}/V$ coincide (since they
are $\K$-invariant). Thus $\sigma_{k}(X)=0.$
\end{proof}
In particular, if $X$ admits a KE metric $\omega,$ then $\sigma_{k}(X)\leq Ak^{-2}$
for some constant $A,$ depending on $\omega.$ However such a bound
cannot be used to obtain an effective criterion for the existence
of a KE metric. Instead we will use the following proposition which,
for any fixed metric $\phi$ in $\mathcal{H}(-K_{X}),$ yields an
effective upper bound of the form 
\[
\frac{K_{k\phi}}{(dd^{c}\phi)^{n}/n!}\leq1+C/k,
\]
 when $k\geq k_{0}$ for explicit constants $k_{0}$ and $C$ depending
on $\phi.$ In fact, $C$ can be taken arbitrarily close to the constant
$C_{\phi_{0}}$ below, by taking $k_{0}$ sufficiently large.
\begin{prop}
\label{prop:bound on Bergman}Given a metric $\phi$ in $\mathcal{H}(L)$
and $q\in]0,1/2[$ there exists a positive constant $C_{\text{\ensuremath{\phi}}}$
such that the following bounds hold on $X,$
\[
\frac{K_{k\phi}}{(dd^{c}\phi)^{n}/n!}\leq\frac{1}{1-C_{\phi}/k-\epsilon_{k}},\,\,\,\epsilon_{k}\coloneqq\frac{1}{n!}\int_{k^{q}}^{\infty}e^{-s}s^{n-1}ds\left(=O(k^{-\infty})\right),
\]
where $C_{\text{\ensuremath{\phi}}}$ can be estimated as follows.
Fix a finite holomorphic atlas $(U_{\alpha},z_{\alpha})$ and relatively
compact subcharts $V_{\alpha}\Subset U_{\alpha}$ covering $X$ and
let $\phi_{U_{\alpha}}$ be the local functions representing $\phi$
on $U_{\alpha}.$ Set 
\[
\lambda_{\alpha}\coloneqq\inf_{V_{\alpha}}\lambda_{\min}\left((\phi_{U_{\alpha}})_{i\bar{j}}\right),\,\,\,M_{\alpha}\coloneqq\max_{\alpha}\left(\sup_{U_{\alpha}}\left|\Delta_{\C}^{2}\phi_{U_{\alpha}}\right|\right),
\]
 where $\lambda_{\min}\left((\phi_{\alpha})_{i\bar{j}}\right)$ denotes
the minimal eigenvalue of the complex Hessian matrix $(\phi_{U_{\alpha}})_{i\bar{j}}$
and $\Delta_{\C}^{2}$ denotes the square of $\sum_{i}\frac{\partial}{\partial z_{i}\partial\overline{z}_{i}}(=\Delta/4).$
Then 
\begin{equation}
C_{\phi}\coloneqq\frac{1}{2}\max_{\alpha}\frac{M_{\alpha}}{\lambda_{\alpha}^{2}},\label{eq:def of C phi K est}
\end{equation}
if $k$ is sufficiently large so that 
\[
\frac{1}{k^{(1/2-q)}}\leq\lambda_{\alpha}^{1/2}\mathrm{dist}(V_{\alpha,}\partial U_{\alpha}).
\]
\end{prop}

\begin{proof}
Given a point $p$ in $X,$ fix local holomorphic coordinates $z$
centered at $x$ that are\emph{ normalized} at $x,$ in the sense
that $\partial\bar{\partial}\phi=\sum_{j}dz_{i}\wedge d\bar{z_{i}}$
at $p.$ We will first establish the following local estimate at the
point $p,$
\begin{equation}
\left(\frac{K_{k\phi}}{k^{n}(dd^{c}\phi)^{n}/n!}\right)(p)\leq\frac{1}{1-C_{p}/k},\,\,\,\,\,C_{p}\coloneqq\frac{1}{2}\sup_{|z|\leq1/k^{(1/2-q)}}\left|\Delta_{\C}^{2}\phi_{U_{\alpha}}\right|.\label{eq:bound in terms of C p}
\end{equation}
 By the normality of the holomorphic coordinates, it is enough to
show that for any local holomorphic function $f(z)$ (not identically
$0$)
\begin{equation}
\frac{\left|f\right|_{k\phi}^{2}(0)}{k^{n}\int_{\left|z\right|\leq R_{k}/k^{1/2}}\left|f\right|_{k\phi}^{2}\frac{i}{2}dz_{1}\wedge d\bar{z}_{1}\cdots}\leq\frac{1}{\pi^{n}}\frac{1}{1-C_{p}/k},\label{eq:upper bound by C p}
\end{equation}
where $R_{k}$ is a sequence of positive numbers such that $R_{k}/k^{1/2}\rightarrow0$
(that we will later take to be $k^{q}).$ We first note that the following
bounds hold\footnote{This formula fixes some typos appearing in the proof of \cite[Thm 2.1]{ber-ber2},
where the constant $|S^{2n-1}|$ and the factor $k$ in front of the
bracket were missing}, by the argument in the beginning of the proof of \cite[Thm 2.1]{ber-ber2}
(using that $\log(\left|f\right|^{2})$ is subharmonic and then applying
Jensen's inequality and the change of variables $r\rightarrow k^{1/2}r)$:
\[
\left|f\right|_{k\phi}^{2}(0)\leq k^{n}\int_{\left|z\right|=r}\left|f\right|_{k\phi}^{2}d\sigma_{r}\left(|S^{2n-1}|\int_{0}^{R_{k}}e^{-\int_{\left|z\right|=r}k\left(\phi(k^{-1/2}z)-\phi(0)\right)d\sigma_{r}}r^{2n-1}dr\right)^{-1},
\]
 where $d\sigma_{r}$ denotes the rotationally invariant probability
measure on the unit-sphere $S^{2n-1}$ in $\C^{n}$ and $|S^{2n-1}|$
denotes the volume of $S^{2n-1}$ wrt the metric on $S^{2n-1}$ induced
by the Euclidean metric on $\C^{n}$. Since the integral of any homogenous
polynomial on $\R^{2n}$ of degree $3$ vanishes (by symmetry), 
\[
\int_{\left|z\right|=r}k\left(\phi(k^{-1/2}z)-\phi(0)\right)d\sigma_{r}=r^{2}+k\int_{\left|z\right|=r}\left(\phi-T_{3}\phi)(k^{-1/2}z)\right)d\sigma_{r},
\]
 where $T_{3}\phi$ denotes the third order jet of $\phi$ around
$z=0.$ In general, for any smooth function $\psi$ in $\R^{2n}$
the following sharp fourth order remainder estimate holds, where $\Delta$
is the Laplacian wrt the Euclidean metric:
\[
\left|\frac{1}{r^{4}}\int_{|x|=r}\bigl(\psi-T_{3}\psi\bigr)\,d\sigma_{r}\right|\le A_{n}\sup_{B_{r}}|\Delta^{2}\psi|,\,\,A_{n}\coloneqq\frac{1}{32(n+1)n}
\]
(which is an equality when $\psi=r^{4}).$ Indeed, applying \cite[formula 4]{ol}
(with $l=1$ in $\R^{2n})$ gives 
\[
A_{n}=c_{2,2n}=\frac{\Gamma(n)}{4^{2}2!\Gamma(2+n)}=\frac{1}{32(n+1)n}.
\]
Thus, by scaling,
\[
\int_{\left|z\right|=r}k\left(\phi(k^{-1/2}z)-\phi(0)\right)d\sigma_{r}\leq r^{2}+k^{-1}C^{(4)}r^{4},\,\,\,C^{(4)}\coloneqq\frac{1}{32n(n+1)}\sup_{|z|\leq R_{k}k^{-1/2}}|\Delta^{2}\phi|
\]
 and, as a consequence, 
\[
\frac{\left|f\right|_{k\phi}^{2}(0)}{k^{n}\int_{\left|z\right|\leq R_{k}/k^{1/2}}\left|f\right|_{k\phi}^{2}\frac{i}{2}dz_{1}\wedge d\bar{z}_{1}\cdots}\leq\left(|S^{2n-1}|\int_{0}^{R_{k}}gdr\right)^{-1},\,\,g\coloneqq e^{-\left(r^{2}+k^{-1}C^{(4)}r^{4}\right)}r^{2n-1}.
\]
Next, expressing $\int_{0}^{R_{k}}gdr=\int_{0}^{\infty}gdr-\int_{R_{k}}^{\infty}gdr$
and using that $C^{(4)}\geq0$ yields 
\[
\int_{0}^{R_{k}}gdr\geq\int_{0}^{\infty}e^{-\left(r^{2}+k^{-1}C^{(4)}r^{4}\right)}r^{2n-1}dr-\epsilon_{k},\,\,\epsilon_{k}\coloneqq\int_{R_{k}}^{\infty}e^{-r^{2}}r^{2n-1}dr.
\]
 Rewriting $e^{-\left(r^{2}+k^{-1}C^{(4)}r^{4}\right)}=e^{-r^{2}}e^{-k^{-1}C^{(4)}r^{4}}$
and using that $e^{t}\geq1+t$ gives 
\[
\int_{0}^{\infty}e^{-\left(r^{2}+k^{-1}C^{(4)}r^{4}\right)}r^{2n-1}dr\geq\int_{0}^{\infty}e^{-r^{2}}r^{2n-1}dr-k^{-1}C^{(4)}B_{n},
\]
where, 

\[
B_{n}\coloneqq\int_{0}^{\infty}r^{4}e^{-r^{2}}r^{2n-1}dr=\frac{(n+1)!}{2}.
\]
Since $|S^{2n-1}|\int_{0}^{\infty}e^{-r^{2}}r^{2n-1}dr=\left(\int_{0}^{\infty}e^{-x^{2}}dx\right)^{2n}=\pi^{n}$
we deduce that
\[
\left(|S^{2n-1}|\int_{0}^{R_{k}}gdr\right)^{-1}\leq\frac{1}{\pi^{n}}\frac{1}{1-\pi^{-n}|S^{2n-1}|(k^{-1}C^{(4)}B_{n}+\epsilon_{k})},\,\,(\frac{|S^{2n-1}|}{\pi^{n}}=\frac{2}{(n-1)!}).
\]
The bound \ref{eq:bound in terms of C p} now follows by taking $R_{k}=k^{q}$
and using that $\pi^{-n}|S^{2n-1}|B_{n}=n(n+1).$ Finally, we note
that $C_{p}\leq C_{\text{\ensuremath{\phi}}},$ where $C_{\text{\ensuremath{\phi}}}$
is defined by formula \ref{eq:def of C phi K est}. This follows directly
from translational invariance of $\Delta^{2}\phi_{\alpha}$ and the
observation that the effect of passing to normalized coordinates at
a given point $p$ gives rise (by scaling) to the denominator $\lambda_{\text{\ensuremath{\alpha}}}^{2}.$
\end{proof}
\begin{rem}
The coefficent $1/2$ appearing in the effective bound in the previous
proposition is essentially sharp; more precisely, the coefficent $1/2$
in the local estimate \ref{eq:bound in terms of C p} is sharp. Indeed,
the value of $\Delta_{\C}^{2}\phi_{U_{\alpha}}$ at $z=0$ coincides
with the scalar curvature $S_{\omega}$ for $\omega\coloneqq\frac{i}{2}dd^{c}\phi,$
if the normalized coordinates are chosen so that $\partial(\phi_{i\bar{j}})(0)=0.$
The sharpness in question thus follows from the asymptotics \ref{eq:exp of Bermgan kernel}.
\end{rem}

\subsection{Effective large deviation bounds }

Given $\epsilon>0$, consider the \emph{thickened partition function}
$Z_{N,\epsilon}(-\gamma),$ defined by
\[
Z_{N,\epsilon}(-\gamma)\coloneqq\int_{\left\{ \left|\boldsymbol{m}_{N}\right|<\epsilon\right\} }\left\Vert \det S^{(N)}\right\Vert _{\phi_{0}}^{-2\gamma/k}\mu_{0}^{\otimes N}.
\]

\begin{lem}
Given two sequences of positive numbers $\epsilon'_{N}$ and $\epsilon_{N}$
such that $\epsilon'_{N}<\epsilon_{N},$
\[
\gamma_{\epsilon_{N}}^{(N)}(X,\Delta)\leq\inf\left\{ \gamma\in]0,\infty[\colon Z_{N,\epsilon'_{N}}(-\gamma)<\infty\right\} .
\]
\end{lem}

\begin{proof}
By the moment-weight inequality in GIT \cite{g-r-s}
\[
-\frac{1}{N}\inf_{\lambda}\frac{w_{\lambda}(x_{1},\dots,x_{N})}{\left|\lambda\right|}\leq\inf_{g\in G}\left|\boldsymbol{m}_{N}(g\cdot(x_{1},\dots,x_{N}))\right|.
\]
In particular, $\left\{ \vert\boldsymbol{m}_{N}((x_{1},\dots,x_{N}))\vert<\epsilon_{N}\right\} \subset(X^{N})_{\text{ss},\epsilon_{N}}.$
This shows that the closure of $\left\{ \boldsymbol{m}_{N}((x_{1},\dots,x_{N}))<\epsilon_{N}'\right\} $
is contained in $(X^{N})_{\text{ss},\epsilon_{N}}.$ Hence, the lemma
follows directly from the analytic definition of the log canonical
threshold.
\end{proof}
We are now ready to state the following result, where
\[
\gamma_{k}\coloneqq\gamma\frac{1+k^{-1}}{1+\gamma k^{-1}}\,\,\,\Bigg(=\frac{1+k^{-1}}{\gamma^{-1}+k^{-1}}\Bigg).
\]

\begin{thm}
\label{thm:effective LD bound with Ding}Let $X$ be a Fano manifold
and fix $\gamma>0.$ If $\epsilon_{N}>(r+1)R_{\mathrm{out}}/N,$ then
there exists a constant $A_{N}$ such that
\[
-\frac{\log Z_{N,\epsilon_{N}}(-\gamma)}{\gamma N(1+k^{-1})}\leq\inf_{u}\left\{ \mathcal{D}_{\gamma_{k}}(u)+\sigma_{k}(\phi_{0})\mathcal{J}(u)\right\} +\sigma_{k}(\phi_{0})C_{\text{sup}}+\frac{\log N!}{N\gamma(1+k^{-1})}+\frac{A_{N}}{\gamma_{k}},
\]
 where the inf ranges over all $u\in\mathcal{E}^{1}(X,\omega_{0})$
satisfying the constraint $\int_{X}\boldsymbol{m}\mu_{\gamma_{k}u}=0.$
Moreover, $\left\{ \left|\boldsymbol{m}_{N}\right|<\epsilon_{N}\right\} $
is non-empty for any $N$ and, as a consequence, so is $(X^{N})_{\text{ss},\epsilon_{N}}.$
\end{thm}

The constant $A_{N}$ depends on $N,R_{\mathrm{out}}$ and $\epsilon_{N},$
but is not explicit (see Lemma \ref{lem:A N is finite} below). However,
by Lemma \ref{lem:Hoeffding}, if $\epsilon_{N}>(2\log2)^{1/2}R_{\text{out}}/N^{1/2},$
\[
A_{N}\leq\frac{1}{N}\log\left(1-2e^{-\frac{1}{2R^{2}}\epsilon_{N}^{2}N}\right).
\]
 Combining the previous theorem with the effective bound on $\sigma_{k}(\phi_{0})$
resulting from Prop \ref{prop:bound on Bergman} we deduce the following
result, containing the first statement in Theorem \ref{thm:uniform polyGibbs implies KE},
when $X$ is a Fano manifold.
\begin{cor}
\label{cor:crit for KE text}Let $X$ be a Fano manifold. Then $X$
admits a unique Kähler--Einstein metric if the following condition
holds for some $k$:
\[
\gamma^{(N)}(X)>\frac{1+n\sigma_{k}(X)}{1-k^{-1}n\sigma_{k}(X)}.
\]
This condition is satisfied if 
\[
\gamma^{(N)}(X)\geq1+C_{1}k^{-1},
\]
holds for some $k\geq k_{1}$ for effective constants $k_{1}$ and
$C_{1}$ (independent of $k).$ More generally, if $\G$ is nontrivial
and the Futaki character of $X$ vanishes, then $X$ admits a KE metric
if the latter lower bounds hold for some $k$ (with a slightly larger
constant $C_{1}$) when $\gamma^{(N)}(X)$ is replaced by $\gamma_{\epsilon_{N}}^{(N)}(X),$
if $\epsilon_{N}>(r+1)R_{\mathrm{out}}/N.$ 
\end{cor}

\begin{rem}
\label{rem:non-reductive groups strong Gibbs poly}Theorem \ref{thm:effective LD bound with Ding}
holds when $\G$ is merely a reductive subgroup of $\mathrm{Aut}_{0}(X)$.
Consequently, strong uniform $\G$-Gibbs polystability (defined as
in Definition \ref{def:strong} but with $\G$ an arbitrary reductive
subgroup) implies the existence of a Kähler--Einstein metric, and
thus further implies that in fact $\mathrm{Aut}_{0}(X)$ is reductive. 
\end{rem}

We will also establish the following improvement of the bound in Theorem
\ref{thm:effective LD bound with Ding}, which yields an effective
large deviation bound, motivated by the conjectural LDP in Conjecture
\ref{conj:LDP}: 
\begin{thm}
\label{thm:effective LD bound with F}Let $X$ be a Fano manifold
and fix $\gamma>0.$ If $\epsilon_{N}>(r+1)R_{\mathrm{out}}/N$ and
$2\epsilon_{N}/R_{\mathrm{in}}<1,$ then there exists a constant $C_{\boldsymbol{m}}$
only depending on $\boldsymbol{m}$ such that
\[
-\frac{\log Z_{N,\epsilon_{N}}(-\gamma)}{N\gamma(1+k^{-1})}\leq\left(1-\frac{2\epsilon_{N}}{R_{\mathrm{in}}}\right)\inf_{\mathcal{P}(X)_{0}}F_{-\tilde{\gamma}}+\sigma_{k}(\phi_{0})C_{\mathrm{sup}}+\frac{2\epsilon_{N}C_{\boldsymbol{m}}}{R_{\mathrm{in}}\gamma}+\frac{A_{N}}{\gamma_{k}}+\frac{\log N!}{N\gamma(1+k^{-1})},
\]
 where $\tilde{\gamma}\coloneqq\gamma_{k}-\left(n\sigma_{k}(\phi_{0})+n^{-1}\epsilon_{N}\right)\left(1-\frac{2\epsilon_{N}}{R_{\mathrm{in}}}\right)^{-1}$
and $A_{N}$ is the constant appearing in Theorem \ref{thm:effective LD bound with Ding}. 
\end{thm}

In fact, the proof shows that $2\epsilon_{N}$ may be replaced by
$\rho\epsilon_{N}$ for any given constant $\rho\in]1,\infty]$ if
the constant $C_{\boldsymbol{m}}$ is taken to depend on $\rho.$
\begin{cor}
\label{cor:lower bound on delta A ps text}Let $X$ be a Fano manifold
and assume that $\epsilon_{N}>(r+1)R_{\mathrm{out}}/N.$ Then there
exist effective constants $k_{2}$ and $C_{2}$ such that for any
$k\geq k_{2}$:
\[
-C_{2}k^{-1}+\frac{1+k^{-1}}{1/\gamma_{\epsilon_{N}}^{(N)}(X)+k^{-1}}\leq\delta_{\text{PS}}^{\mathrm{A}}(X,\Delta).
\]
\end{cor}

Letting $N\rightarrow\infty$, we also deduce the following result,
which yields the upper bound in Conjecture \ref{conj:large N limtis intro},
for the ``thickened'' partition function $Z_{N,\epsilon_{N}}$:
\begin{cor}
\label{cor:asymptotic upper bound on minus log Z text}Let $X$ be
a Fano manifold and $\epsilon_{N}$ be a sequence of positive numbers
such that $\epsilon_{N}\geq CN^{-1/2}$ for a sufficiently large constant
$C$ (more precisely, $C>(2\log2)^{1/2}R_{\mathrm{out}}).$ Then 
\[
\limsup_{N\rightarrow\infty}-\frac{1}{N\gamma}\log Z_{N,\epsilon_{N}}(-\gamma)\leq\inf_{\mathcal{P}(X)_{0}}F_{-\gamma}
\]
and if the inf in the rhs above is attained then 
\[
-\frac{1}{N\gamma}\log Z_{N,\epsilon_{N}}(-\gamma)\leq\inf_{\mathcal{P}(X)_{0}}F_{-\gamma}+O\Big(\frac{\log N}{k}\Big)+O(N^{-1/2}).
\]
 In particular, if $X$ admits a KE metric, then 
\[
-\frac{1}{N}\log Z_{N,\epsilon_{N}}(-1)\leq\inf_{\mathcal{H}(X,\omega_{0})}\mathcal{M}+O\Big(\frac{\log N}{k}\Big)+O(N^{-1/2}).
\]
\end{cor}

\begin{rem}
Since $O(k^{-1})=O(N^{-1/n})$ the error term in the previous corollary
may be replaced by $O(N^{-1/2})$ when $n=1$ and by $O(\frac{\log N}{k})$
when $n>1.$ 
\end{rem}

\subsection{The proof of Theorem \ref{thm:effective LD bound with Ding}}

Given a function $u$ on $X$ such that $e^{-\gamma u}\in L^{1}(X)$
we will denote by $\mu_{\gamma u}$ the probability measure $\mu_{\gamma u}\coloneqq e^{-\gamma u}\mu_{0}/\int_{X}e^{-\gamma u}\mu_{0}$
on $X.$ 

\subsubsection{Inequality in terms of a quantized Ding functional and a mean entropy
penalty for a finite $N$}

Given $k,$ denote by $\mathcal{D}_{k,-\gamma}$ the ``quantized''
twisted Ding functional on $\mathcal{E}^{1}(X)$ defined by
\begin{equation}
\mathcal{D}_{k,\gamma}(u)\coloneqq-\mathcal{E}_{k}(u)-\frac{1}{\gamma}\log\int_{X}e^{-\gamma u}\mu_{0},\label{eq:def of D k gamma on H}
\end{equation}
 where $\mathcal{E}_{k}(u)$ is the following variant of Donaldson's
$\mathcal{L}_{k}-$functional:
\[
\mathcal{E}_{k}(u)\coloneqq-\frac{1}{N(k+\gamma)}\log\det\left(\int_{X}\left\langle s_{i},s_{j}\right\rangle _{\phi_{0}+u}e^{-\gamma u}\mu_{0}\right).
\]
Here $\left\langle \cdot,\cdot\right\rangle _{\phi_{0}+u}$ denotes
the Hermitian metric on $-kK_{(X,\Delta)}$ induced by the metric
$\phi_{0}+u$ on $-K_{(X,\Delta)}$, and $s_{1},\dots,s_{N}$ denotes
any fixed basis of $H^{0}(X,-kK_{(X,\Delta)})$ which is orthonormal
with respect to the Hermitian product on $H^{0}(X,-kK_{(X,\Delta)})$
induced by the metric $\phi_{0}$ on $-K_{(X,\Delta)}$ and the measure
$\mu_{\phi_{0}}$ on $X.$ The normalization have been chosen to ensure
that $\mathcal{D}_{k,\gamma}(u)$ is invariant under $u\mapsto u+c,$
for any $c\in\R.$ The following inequality generalizes \cite[Prop 2.3]{berm13}
(where $\boldsymbol{m}_{N}\equiv0)$:
\begin{prop}
\label{prop:bounding Z thick by q Ding etc}For any $u$ such that
$e^{-\gamma u}\in L^{1}(X)$:
\begin{equation}
-\frac{1}{\gamma N(1-\gamma/k)}\log Z_{N,\epsilon}(-\gamma)\leq\mathcal{D}_{k,\gamma}(u)-\frac{1}{\gamma}\frac{1}{N}\log\mu_{\gamma u}^{\otimes N}\left(\left\{ \left|\boldsymbol{m}_{N}\right|<\epsilon\right\} \right)+\frac{\log N!}{N(k-\gamma)}.\label{eq:ineq for quant Ding with mean entropy}
\end{equation}
\end{prop}

\begin{proof}
Setting $U_{\epsilon}\coloneqq\left\{ \left|\boldsymbol{m}_{N}\right|<\epsilon\right\} \subset X^{N}$,
we first rewrite
\[
\left\Vert \det S^{(N)}\right\Vert _{\phi_{0}}^{-2\gamma/k}\mu_{0}^{\otimes N}=\left\Vert \det S^{(N)}\right\Vert _{\phi_{0}+u}^{-2\gamma/k}\left(e^{-\gamma u}\mu_{0}\right)^{\otimes N}=
\]
\[
=\left\Vert \det S^{(N)}\right\Vert _{\phi_{0}+u}^{-2\gamma/k}\left(\frac{\left(e^{-\gamma u}\mu_{0}\right)^{\otimes N}}{\int_{U_{\epsilon}}\left(e^{-\gamma u}\mu_{0}\right)^{\otimes N}}\right)\int_{U_{\epsilon}}\left(e^{-\gamma u}\mu_{0}\right)^{\otimes N}.
\]
Applying Hölder's inequality with negative exponent $-\gamma/k$ (or
applying Jensen's inequality to the convex function $t\mapsto t^{-\gamma/k}$
on $]0,\infty[$) thus yields
\[
Z_{N,\epsilon}(-\gamma)^{-1}\leq\left(\int_{U_{\epsilon}}\left\Vert \det S^{(N)}\right\Vert _{\phi_{0}+u}^{2}\frac{\left(e^{-\gamma u}\mu_{0}\right)^{\otimes N}}{\int_{U_{\epsilon}}e^{-\gamma u}\mu_{0}}\right)^{\gamma/k}\left(\int_{U_{\epsilon}}\left(e^{-\gamma u}\mu_{0}\right)^{\otimes N}\right)^{-1}=
\]
\[
\left(\int_{U_{\epsilon}}\left\Vert \det S^{(N)}\right\Vert _{\phi_{0}+u}^{2}\left(e^{-\gamma u}\mu_{0}\right)^{\otimes N}\right)^{\gamma/k}\left(\int_{U_{\epsilon}}\left(e^{-\gamma u}\mu_{0}\right)^{\otimes N}\right)^{-(1-\gamma/k)}
\]
Bounding the integral over $U_{\epsilon}$ in the first bracket above
by an integral over all of $X^{N},$ taking logarithms and dividing
by $\gamma N(1-\gamma/k),$ we deduce that 
\[
-\frac{1}{\gamma N(1-\gamma/k)}\log Z_{N,\epsilon}(-\gamma)\leq\frac{\log\left(\int_{X^{N}}\left\Vert \det S^{(N)}\right\Vert _{\phi_{0}+u}^{2}\left(e^{-\gamma u}\mu_{0}\right)^{\otimes N}\right)}{N(k-\gamma)}-\frac{\log\int_{U_{\epsilon}}\left(e^{-\gamma u}\mu_{0}\right)^{\otimes N}}{\gamma N}
\]
Next, applying the previous inequality to $u-\log\int_{X}e^{-\gamma u}\mu_{0}$,
the right hand side becomes
\[
\frac{1}{N(k-\gamma)}\log\left(\int_{X^{N}}\left\Vert \det S^{(N)}\right\Vert _{\phi_{0}+u}^{2}\mu_{\gamma u}^{\otimes N}\right)-\frac{1}{\gamma}\log\int_{X}e^{-\gamma u}\mu_{0}-\frac{1}{\gamma}\frac{1}{N}\log\int_{U_{\epsilon}}\mu_{\gamma u}^{\otimes N},
\]
using that the latter right hand side is invariant under $u\mapsto u-c$
for any $c\in\R.$ Finally, the desired inequality \ref{eq:ineq for quant Ding with mean entropy}
follows from the following basic identity \cite[Lemma 5.3]{b-b}:
\begin{equation}
\int_{X^{N}}\left\Vert \det S^{(N)}\right\Vert _{\phi_{0}+u}^{2}\mu_{\gamma u}^{\otimes N}=N!\det\left(\int_{X}\left\langle s_{i},s_{j}\right\rangle _{\phi_{0}+u}e^{-\gamma u}\mu_{0}\right).\label{eq:integral det gives Don func}
\end{equation}
\end{proof}

\subsubsection{An effective upper bound on the quantized Ding functional}
\begin{prop}
\label{prop:effect upper bound q Ding}Given $u\in\mathcal{E}^{1}(X,\omega)$
and $\gamma>0,$ define $\epsilon$ by $1-\epsilon\coloneqq(1+k^{-1})/(1+\gamma k^{-1}).$
Then 
\[
\frac{1}{(1-\epsilon)}\mathcal{D}_{k,\gamma}((1-\epsilon)u)\leq\mathcal{D}_{(1-\epsilon)\gamma}(u)+\sigma_{k}(\phi_{0})\left(\mathcal{J}(u)+C_{\mathrm{sup}}\right).
\]
\end{prop}

\begin{proof}
Since this bound follows from the proof of \cite[Lemma 2.4]{berm13},
together with the scaling argument in the beginning of \cite[Section 2.3]{berm13},
we will only outline the argument. We will identify a metric $\phi\in\mathcal{H}(-K_{X})$
with the function $u\coloneqq\phi-\phi_{0}\in\mathcal{H}(X,\omega_{0}).$
Moreover, setting $L\coloneqq-K_{X}$ and rewriting $-kK_{X}=(k+1)L+K_{X}$,
a section $s_{k}\in H^{0}(X,-kK_{X})$ will be identified with a holomorphic
top form on $X$ with values in $(k+1)L$ satisfying
\[
\int_{X}\left\Vert s_{k}\right\Vert _{\phi_{0}+u}^{2}e^{-\gamma u}\mu_{0}=\int_{X}i^{n^{2}}s_{k}\wedge\overline{s_{k}}e^{-(k\phi+\gamma\phi+(1-\gamma)\phi_{0})}.
\]
Then formula \ref{eq:integral det gives Don func} reveals that $\mathcal{E}_{k}\left(\phi\right)$
coincides with a negative multiple of the logarithm of the $L^{2}-$norm
on the top exterior power of $H^{0}(X,(k+1)L+K_{X}),$ when $L$ is
endowed with the metric $k\phi+\gamma\phi+(1-\gamma)\phi_{0}.$ Note
that replacing $\phi$ with 
\[
\phi^{(\epsilon)}\coloneqq\phi(1-\epsilon)+\epsilon\phi_{0},\,\,\,\epsilon=\frac{\gamma-1}{k+\gamma},
\]
 the corresponding metric on $(k+1)L$ has non-negative curvature.
Let us first show that
\begin{equation}
-\frac{1}{(1-\epsilon)}\mathcal{E}_{k}\left(\phi^{(\epsilon)}\right)+\mathcal{E}(\phi)\leq\sigma_{k}(\phi_{0})\left(-\mathcal{E}(\phi)+\sup_{X}(\phi-\phi_{0})\right).\label{eq:ineq in pf of bound on q Ding}
\end{equation}
To this end let $\phi_{t}$ be the weak geodesic in\emph{ $\mathcal{H}^{1,1}(-K_{X})$}
connecting $\phi$ with $\phi_{0}.$ Then 
\begin{equation}
(i)\,t\mapsto\mathcal{E}(\phi_{t})\,\text{is affine,\text{\,\,\,\,\,}\ensuremath{(ii)\,t}\ensuremath{\mapsto}\ensuremath{\mathcal{E}\left(\phi_{t}^{(\epsilon)}\right)}}\,\text{is\,convex}\label{eq:i and ii}
\end{equation}
(this corrects a typographical error in \cite[formula 2.8]{berm13},
where \textquotedblleft convex\textquotedblright{} is mistakenly written
as \textquotedblleft concave\textquotedblright ). The first statement
follows from \cite[Thm 1.7]{bbj} and the second one from \cite{bern0},
using that $\phi_{t}^{(\epsilon)}$ is a subgeodesic (i.e. $\phi_{t}^{(\epsilon)}(z,t)$
is psh in $(z,t),$ when $\phi$ is extended to the strip $[0,1]\times i\R\subset\C$
so that $\phi(t)\coloneqq\phi(\text{Re}\,t).$ As a consequence, the
function $t\mapsto-\frac{1}{(1-\epsilon)}\mathcal{E}_{k}\left(\phi_{t}^{(\epsilon)}\right)+\mathcal{E}(\phi_{t})$
is concave, giving, 
\[
-\frac{1}{(1-\epsilon)}\mathcal{E}_{k}\left(\phi^{(\epsilon)}\right)+\mathcal{E}(\phi)\leq-\frac{1}{1-\epsilon}\mathcal{E}_{k}\left(\phi_{0}^{(\epsilon)}\right)+\mathcal{E}(\phi_{0})+\int_{X}\left(-\frac{d\phi_{t}}{dt}\Big|_{t=0}\right)\left(\frac{1}{(1-\epsilon)}d\mathcal{E}_{k}\left(\phi^{(\epsilon)}\right)-d\mathcal{E}\right)\Bigg|_{\phi_{0}}
\]
At $\phi=\phi_{0}$, the differential $d\mathcal{E}$ equals $(dd^{c}\phi_{0})^{n}/V$
(by \ref{eq:diff of beaut E}) and the differential $\frac{1}{(1-\epsilon)}d\mathcal{E}_{k}\left(\phi^{(\epsilon)}\right)$
equals $\beta_{(k+1)\phi_{0}}$ (by \cite[Lemma 2.1]{b-b}). Since
the lhs in the inequality \ref{eq:ineq in pf of bound on q Ding}
is invariant under additive scalings of $\phi$ we may as well assume
that $\phi$ is sup-normalized, i.e. $\phi-\phi_{0}\leq0.$ Then $\frac{d\phi_{t}}{dt}|_{t=0}\leq0$
(since $t\mapsto\phi_{t}(z)$ is convex) and thus 
\[
\int_{X}\left(-\frac{d\phi_{t}}{dt}\Big|_{t=0}\right)\left(\frac{1}{(1-\epsilon)}d\mathcal{E}_{k}\left(\phi^{(\epsilon)}\right)-d\mathcal{E}\right)\Bigg|_{\phi_{0}}\leq\sigma_{k}(\phi_{0})\int_{X}\left(-\frac{d\phi_{t}}{dt}\Big|_{t=0}\right)(dd^{c}\phi_{0})^{n}/V.
\]
Since $t\mapsto\mathcal{E}(\phi_{t})$ is affine, the integral in
the rhs above equals $-\mathcal{E}(\phi)$ and thus the inequality
\ref{eq:ineq in pf of bound on q Ding} follows. Finally, combining
this inequality with the sup-bound \ref{eq:bound on sup u}, the proof
of the lemma for $u\in\mathcal{H}(X,\omega_{0})$ follows directly
from the definitions of $\mathcal{D}_{k,\gamma}$ and $\mathcal{D}.$
Since any $u\in\text{PSH }(X,\omega_{0})$ is a decreasing limit of
functions $u_{j}\in\mathcal{H}(X,\omega_{0})$ (by Demailly's approximation
result), it follows that the lemma in fact holds for any $u\in\mathcal{E}^{1}(X,\omega),$
using also that $e^{-\gamma u}\in L^{1}(X)$ for any $\gamma>0$ if
$u\in\mathcal{E}^{1}(X,\omega)$ \cite{bbegz}.
\end{proof}

\subsubsection{Upper bounds on the mean entropy for measures with mean zero}

Denote by $\mathcal{S}$ the space of all probability measure $\mu$
on $\R^{r}$ such that $\mu$ is supported on the ball of radius $R$
and $\int\boldsymbol{x}\mu=0.$ The following lemma will be used to
prove the first inequality in Theorem \ref{thm:effective LD bound with Ding}:
\begin{lem}
\label{lem:A N is finite}For any positive integer $N,$ 
\[
e^{-NA_{N}}\coloneqq\inf_{\mu\in\mathcal{S}}\mu^{\otimes N}\left(\left\{ (\boldsymbol{x}_{1},\dots,\boldsymbol{x}_{N})\in(\R^{r})^{N}\colon\left|\frac{1}{N}\sum_{i=1}^{N}\boldsymbol{x}_{i}\right|<\frac{(r+1)R}{N}\right\} \right)>0.
\]
\end{lem}

\begin{proof}
Denote by $K$ the support of $\mu.$ Let us first show that
\begin{equation}
U_{N}\coloneqq\left\{ (\boldsymbol{x}_{1},\dots,\boldsymbol{x}_{N})\in(\R^{r})^{N}\colon\left|\frac{1}{N}\sum_{i=1}^{N}\boldsymbol{x}_{i}\right|<\frac{(r+1)R}{N}\right\} \cap K^{N}\neq\emptyset.\label{eq:of of lemma A N finite}
\end{equation}
 Since $\int_{K}\boldsymbol{x}\mu=0,$ the convex hull $K$ contains
$0.$ Hence, by Caratheodory's theorem, there exist $r+1$ points
$\boldsymbol{a}_{1},\dots,\boldsymbol{a}_{r+1}$ in $K$ and $\lambda_{i}\in[0,1]$
such that 
\[
\lambda_{1}\boldsymbol{a}_{1}+\cdots+\lambda_{r}\boldsymbol{a}_{r+1}=0,\,\,\,\,\sum_{i=1}^{r+1}\lambda_{i}=1.
\]
Expressing $N\lambda_{i}=n_{i}+\delta_{i}\forall i,$ where $n_{i}$
are non-negative integers, $|\delta_{i}|<1$ and $\sum_{i=1}^{r+1}n_{i}=N$,
we deduce that d
\[
|n_{1}\boldsymbol{a}_{1}+\cdots+n_{r}\boldsymbol{a}_{r+1}|\leq|\boldsymbol{a}_{1}|+\cdots+|\boldsymbol{a}_{r+1}|\leq(r+1)R
\]
Hence $(\boldsymbol{x}_{1},\dots,\boldsymbol{x}_{N})\in U_{N}$ if
$\boldsymbol{x}_{1},\dots,\boldsymbol{x}_{n_{1}}\coloneqq\boldsymbol{a}_{1}$
and so on, proving \ref{eq:of of lemma A N finite}. Finally, endow
the space $\mathcal{S}$ with the weak topology and consider the functional
on $\mathcal{S}$ defined by $\mu\mapsto\mu^{\otimes N}(U_{N}).$
Since $U_{N}$ is open this functional is lower semi-continuous and
$\mu^{\otimes N}(U_{N})>0$ for any $\mu\in\mathcal{S},$ by \ref{eq:of of lemma A N finite}.
The space $\mathcal{S}$ being compact, it thus follows that the infimum
of the functional in question is strictly positive, as desired. 
\end{proof}
The following lemma, which is standard in the large deviation theory
for independent random vectors, shows that $A_{N}\leq-N^{-1}\log(1-2e^{-\frac{1}{2R^{2}}\epsilon^{2}N})$
when $2e^{-\frac{1}{2R^{2}}\epsilon^{2}N}<1$:
\begin{lem}
\label{lem:Hoeffding}For any $\mu\in\mathcal{S}$,
\[
\mu^{\otimes N}\left(\left\{ \left|\frac{1}{N}\sum_{i=1}^{N}\boldsymbol{x}_{i}\right|<\epsilon\right\} \right)\geq1-2e^{-\frac{1}{2R^{2}}\epsilon^{2}N}.
\]
\end{lem}

\begin{proof}
It is equivalent to show the bound
\begin{equation}
\mu^{\otimes N}\left(\left\{ \left|\frac{1}{N}\sum_{i=1}^{N}\boldsymbol{x}_{i}\right|\geq\epsilon\right\} \right)\leq2e^{-\frac{1}{2R^{2}}\epsilon^{2}N}.\label{eq:Hoeffding}
\end{equation}
When $m=1$, i.e. $\boldsymbol{x}_{i}$ is real-valued, the latter
bound follows directly from the classical \emph{Hoeffding's inequality
}for the empirical sum of independent real-valued random variables
of mean zero\emph{ }\cite{hoe}. The general case follows from the
large deviation bound \cite[Thm 3.5]{pi} for Banach space-valued
martingales (using that $D=1$ for a Hilbert space), applied to the
Euclidean Hilbert space $\R^{m}$ with the martingale $f_{j}$ defined
by $f_{0}=0,$ $f_{j}=\sum_{i=1}^{j}\boldsymbol{x}_{i}$ when $j\in[1,N]$
and $f_{j}=f_{N}$ when $j>N.$ 
\end{proof}

\subsubsection{Conclusion of the proof of Theorem \ref{thm:effective LD bound with Ding}}

Fix $u\in\mathcal{E}^{1}(X,\omega_{0}),\gamma>0$ and $\epsilon\in\R.$
Applying the upper bound on $-\frac{1}{N\gamma}\log Z_{N,\epsilon}(-\gamma)$
in Prop \ref{prop:bounding Z thick by q Ding etc} to $(1-\epsilon)u\in L^{1}(X)$
for a given $u\in\mathcal{E}^{1}(X)$ gives
\begin{equation}
-\frac{1}{\gamma N(1-\gamma/k)}\log Z_{N,\epsilon}(-\gamma)\leq\mathcal{D}_{k,\gamma}((1-\epsilon)u)-\frac{1}{\gamma}\frac{1}{N}\log\mu_{\gamma(1-\epsilon)u}^{\otimes N}\left(\left\{ \left|\boldsymbol{m}_{N}\right|<\epsilon\right\} \right)+\frac{\log N!}{N(k-\gamma)}.\label{eq:ineq in concl of pf LDP Ding}
\end{equation}
Now take $\epsilon$ to be as in Prop \ref{prop:effect upper bound q Ding},
so that $\gamma(1-\epsilon)\coloneqq\gamma_{k}.$ Then Prop \ref{prop:effect upper bound q Ding}
gives 

\begin{align}
-\frac{1}{\gamma N(1-\gamma/k)(1-\epsilon)} & \log Z_{N,\epsilon}(-\gamma)\leq\label{eq:ineq for quant Ding with mean entropy-1-1-1}\\
\mathcal{D}_{\gamma_{k}}(u)+\sigma_{k}(\phi_{0})\left(\mathcal{J}(u)+C_{\text{sup}}\right)-\frac{1}{\gamma(1-\epsilon)}\frac{1}{N} & \log\mu_{\gamma(1-\epsilon)u}^{\otimes N}\left(\left\{ \left|\boldsymbol{m}_{N}\right|<\epsilon\right\} \right)+\frac{\log N!}{N\gamma(1+k^{-1})}.
\end{align}

The proof of the inequality in Theorem \ref{thm:effective LD bound with Ding}
is thus concluded by invoking Lemma \ref{lem:A N is finite} to bound
the term above involving $\mu_{\gamma(1-\epsilon)u}^{\otimes N},$
using the assumption that $u$ satisfies the constraint $\int_{X}\boldsymbol{m}\mu_{\gamma_{k}u}=0,$
also using that $\gamma_{k}=\gamma(1-\epsilon)$ and $(1-\gamma/k)(1-\epsilon)=1+k^{-1}.$

Finally, to prove the last statement in Theorem \ref{thm:effective LD bound with Ding}
it is, by Lemma \ref{lem:A N is finite}, enough to show that there
exists $\nu\in\mathcal{P}(X)_{0}$ (since we can then apply Lemma
\ref{lem:A N is finite} to $\mu\coloneqq\boldsymbol{m}_{*}\nu$).
The proof is thus concluded by noting that $\mu_{\phi_{0}}\in\mathcal{P}(X)_{0}.$
To see this fix $W$ in the Lie algebra of $\K$ and observe that
the function $t\mapsto\int_{X}\mu_{(e^{tJW})^{*}\phi_{0}}$ is constant,
since the integral coincides with $\int_{X}(e^{tJW})^{*}\mu_{\phi_{0}}.$
Differentiating at $t=0$ thus reveals that 
\begin{equation}
\int h_{W}\mu_{\phi_{0}}=0,\label{eq:average of Ham vanishes}
\end{equation}
where $h_{W}$ denotes the Hamiltonian function corresponding to $W.$
Hence, $\mu_{\phi_{0}}\in\mathcal{P}(X)_{0}.$

\subsection{Proof of Cor \ref{cor:crit for KE text}}

Take $\gamma\in]0,\gamma^{(N)}(X)[$ and $\phi_{0}\in\mathcal{H}(X,\omega_{0}).$
Then, by replacing $\G$ with the trivial group (so that $\boldsymbol{m\equiv}0)$,
Theorem \ref{thm:effective LD bound with Ding} implies that 
\begin{equation}
F_{-(\gamma_{k}-n\sigma_{k}(\phi_{0}))}\geq-C_{0}\label{eq:lower boundon F in pf sufficent crit}
\end{equation}
 on $\mathcal{P}(X)$, for some constant $C_{0}$ depending on $\phi_{0}$
and $\gamma,$ using that $\mathcal{D}_{\gamma}(u)\leq\gamma^{-1}F_{-\gamma}(\omega_{u}/V)$
and $\mathcal{J}(u)\leq n(\mathcal{I}-\mathcal{J})(u)=E(\omega_{u}/V).$
In particular, if 
\begin{equation}
\gamma_{k}\coloneqq\frac{1+k^{-1}}{\gamma^{-1}+k^{-1}}>1+n\sigma_{k}(\phi_{0})\iff\gamma>\frac{1+n\sigma_{k}(\phi_{0})}{1-k^{-1}n\sigma_{k}(\phi_{0})}\label{eq:ineq in pf of Cor crit for KE}
\end{equation}
it then follows from the lower bound \ref{eq:lower boundon F in pf sufficent crit}
that $X$ admits a unique KE metric (since it implies that the Mabuchi
functional $\mathcal{M}$ is coercive). The proof of the first statement
in Cor \ref{cor:crit for KE text} is thus concluded by noting that
it follows directly from the definitions that the inequality in the
rhs of formula \ref{eq:ineq in pf of Cor crit for KE} holds for some
$\gamma\in]0,\gamma^{(N)}(X)[$ and $\phi_{0}\in\mathcal{H}(X,\omega_{0})$
if $\gamma^{(N)}(X)>\frac{1+n\sigma_{k}(X)}{1-k^{-1}n\sigma_{k}(X)}.$
Finally, the existence of the effective constants $k_{1}$ and $C_{1}$
now follows from combining the effective bound in Prop \ref{prop:bound on Bergman}
with the following well-known lower bound that holds for any $k$:
\begin{equation}
\frac{N_{k}}{k^{n}c_{1}(-K_{X})^{n}/n!}\geq1-C_{\text{RR}}k^{-1},\label{eq:effective RR bound}
\end{equation}
 for an effective constant $C_{\text{RR}},$ resulting from the Riemann--Roch
theorem (using that, by Kodaira vanishing, the Dolbeault cohomology
groups $H^{q}(X,-kK_{X})$ are trivial for any $k$ and $q>0).$ 

Now consider the case when $\G$ is nontrivial and assume that the
Futaki character of $X$ vanishes. Take $\gamma\in]0,\gamma_{\epsilon_{N}}^{(N)}(X)[$
and $\phi_{0}\in\mathcal{H}(X,\omega_{0})^{\K}.$ Then Theorem \ref{thm:effective LD bound with Ding}
implies that $\mathcal{D}_{\gamma_{k}}(u)+\sigma_{k}(\phi_{0})\mathcal{J}(u)$
is bounded from below on the subspace of $\mathcal{E}^{1}(X,\omega_{0})$
satisfying the constraint $\mu_{\gamma_{k}u}\in\mathcal{P}(X)_{0}.$
By Prop \ref{prop:constrained Ding vs properness}, this implies that
$\mathcal{D}(v)\geq\left(1-\frac{(1+\sigma_{k}(\phi_{0}))}{\gamma_{k}^{1/n}}\right)\mathcal{J}^{\G}(v)+\gamma_{k}C$,
as long as $\gamma_{k}\geq1$. If 
\[
1-\frac{(1+\sigma_{k}(\phi_{0})}{\gamma_{k}^{1/n}}>0\iff\gamma_{k}>\left(1+\sigma_{k}(\phi_{0})\right)^{n},
\]
it thus follows that $X$ admits a KE, since $\mathcal{D}$ is then
$\G$-coercive. Finally, using that $(1+t)^{n}=1+nt+O(t^{2}),$ as
$t\rightarrow0,$ we can thus conclude the proof as before, by invoking
the effective bounds in Prop \ref{prop:bound on Bergman} and formula
\ref{eq:effective RR bound}. 

\subsection{Proof of Theorem \ref{thm:effective LD bound with F}}

To prove Theorem \ref{thm:effective LD bound with F} it will be enough
to establish the following upper bound on the mean entropy, where
$\delta_{\gamma}(u)$ denotes the entropic distance defined by formula
\ref{eq:def of entropic distance}:
\begin{prop}
\label{prop:upper bound mean entr}Given $\gamma>0$ and $u\in\mathcal{E}^{1}(X)$,
the following inequality holds for any $\epsilon>0$: 
\[
-\frac{1}{N}\log\int_{\left\{ \left|\boldsymbol{m}_{N}\right|<\epsilon\right\} }\mu_{\gamma u}^{\otimes N}\leq\gamma\delta_{\gamma}(u)-\frac{1}{N}\inf_{\mu\in\mathcal{S}}\log\int_{\left\{ \left|N^{-1}\sum_{i=1}^{N}\boldsymbol{x}_{i}\right|<\epsilon\right\} }\mu^{\otimes N}+|\boldsymbol{\beta}_{0}|\epsilon,
\]
 where $|\boldsymbol{\beta}_{0}|$ may be bounded in terms of a constant
$C_{\boldsymbol{m}}$ only depending on $\boldsymbol{m}$ as follows:
\[
|\boldsymbol{\beta}_{0}|\leq-\frac{2}{R_{\mathrm{in}}}\gamma\left(\mathcal{D}_{\gamma}(u)+\delta_{\gamma}(u)+\sup_{X}u-\mathcal{E}(u)-\frac{C_{\boldsymbol{m}}}{\gamma}\right).
\]
\end{prop}

Indeed, proceeding as in the proof of Theorem \ref{thm:effective LD bound with Ding},
we can then bound the term involving $\mu_{\gamma(1-\epsilon)u}^{\otimes N}$
in the inequality \ref{eq:ineq in concl of pf LDP Ding}, by applying
the previous proposition with $\gamma$ replaced by $(1-\epsilon)\gamma(=\gamma_{k})$
to deduce that $\frac{1}{\gamma N(1-\gamma/k)(1-\epsilon)}\log Z_{N,\epsilon}(-\gamma)$
is bounded from above by 
\[
\mathcal{D}_{\gamma_{k}}(u)+\sigma_{k}(\phi_{0})\left(\mathcal{J}(u)+C_{\text{sup}}\right)+\delta_{\gamma_{k}}(u)
\]
\[
-\frac{2\epsilon_{N}}{R_{\text{in}}}\left(\mathcal{D}_{\gamma_{k}}(u)+\delta_{\gamma_{k}}(u)+\sup_{X}u-\mathcal{E}(u)-\frac{C_{\boldsymbol{m}}}{\gamma}\right)+\frac{A_{N}}{\gamma_{k}}+\frac{\log N!}{N(k-\gamma)(1-\epsilon)}.
\]
Finally, assuming that $\mu\coloneqq\omega_{u}^{n}/V\in\mathcal{P}(X)_{0}$
and using that $\mathcal{D}_{\gamma_{k}}(u)+\delta_{\gamma_{k}}(u)\leq F_{\gamma_{k}}(\omega_{u}^{n}/V)$
(by Lemma \ref{lem:constrained F geq penalized Ding}) rearranging
the terms above and using the energy comparison inequalities \ref{eq:comparison I minus J and J}
concludes the proof of Theorem \ref{thm:effective LD bound with F}.

\subsubsection{Proof of Prop \ref{prop:upper bound mean entr}}

We start with the following
\begin{lem}
Let $X$ be a Fano manifold and $\K$ a compact Lie group acting holomorphically
and faithfully on $X.$ Fix a $\K$-invariant smooth metric on $-K_{X}$
with positive curvature and denote by $\boldsymbol{m}$ the corresponding
moment map from $X$ into the dual $\mathfrak{k}^{*}$ of the Lie
algebra $\mathfrak{k}$ of $\K.$ Then $0$ is contained in the interior
of the convex hull of $\boldsymbol{m}(X).$ As a consequence, if $\nu$
is a finite measure on $X$ of the form $\nu=e^{-v}dV$ for $v\in L^{1}(X)$
and $dV$ a smooth volume form on $X,$ then $0$ is contained in
the interior of the convex hull of the support of $\boldsymbol{m}_{*}\nu.$
\end{lem}

\begin{proof}
We first claim that any linear projection on $\R$ of $\boldsymbol{m}(X)$
is an interval containing $0$ in its interior. To prove this, consider,
for a fixed element $\boldsymbol{\beta}$ in $\mathfrak{k},$ the
smooth function $h\coloneqq\left\langle \boldsymbol{\beta},\boldsymbol{m}(x)\right\rangle $
on $X.$ This is the Hamiltonian of the vector field on $X$ induced
by $\boldsymbol{\beta}$ and it is well-known that $h(X)$ is an interval
containing $0$ in its interior (as follows directly from the identity
\ref{eq:average of Ham vanishes}). This proves the claim in question.
It follows that $0$ is contained in the interior of any linear projection
of the convex hull $C$ of $\boldsymbol{m}(X)$ on $\R.$ Indeed,
since the image of the projection of $\boldsymbol{m}(X)$ under $\left\langle \boldsymbol{\beta},\cdot\right\rangle $
is an interval the image coincides, by linearity, with the projection
of $C.$ Finally, assume, to get a contradiction, that $0$ is not
contained in the interior of $C.$ Then, since $C$ is a closed convex
set, there exists (by the supporting hyperplane theorem) $\boldsymbol{\beta}\in\mathfrak{k}$
such that $\left\langle \beta,C\right\rangle \leq0.$ But this contradicts
that $0$ is contained in the interior of the linear projection of
$C$ under $\left\langle \boldsymbol{\beta},\cdot\right\rangle $.
Finally, to conclude the proof of the lemma it will be enough to show
that $\boldsymbol{m}(X)$ is the support of $\boldsymbol{m}_{*}\nu.$
To this end, assume for contradiction, that there exists $x\in X$
such that $\boldsymbol{m}(x)$ is in the complement of the support
of $\boldsymbol{m}_{*}\nu.$ This means that there exists a small
open ball $B$ containing $\boldsymbol{m}(x)$ such that 
\[
0=(\boldsymbol{m}_{*}\nu)(B)\coloneqq\nu(\boldsymbol{m}^{-1}(B))\coloneqq\int_{\boldsymbol{m}^{-1}(B)}e^{-v}dV.
\]
 But, since, $\boldsymbol{m}^{-1}(B)$ is an open subset of $X$ this
contradicts the assumption that $v\in L^{1}(X).$ 
\end{proof}
We will use the previous lemma in conjunction with the following general
lemma, where $V$ denotes a finite dimensional vector space and $V^{*}$
its dual, and the Legendre transforms of convex lsc functions $\phi$
and $\psi$ on $V$ and $V^{*},$ respectively, are the functions
on $V^{*}$ and $V$ defined by
\[
\phi^{*}(\boldsymbol{e})\coloneqq\sup_{\boldsymbol{\beta}\in V}\left(\left\langle \boldsymbol{\beta},\boldsymbol{e}\right\rangle -\phi(\boldsymbol{\beta})\right),\,\,\psi^{*}(\boldsymbol{\beta})\coloneqq\sup_{\boldsymbol{e}\in V^{*}}\left(\left\langle \boldsymbol{\beta},\boldsymbol{e}\right\rangle -\psi(\boldsymbol{e})\right).
\]
By classical duality, $\phi^{**}=\phi$ and $\psi^{**}=\psi.$ 
\begin{lem}
\label{lem:Legendre duality for m}Let $X$ be a compact topological
space and $\boldsymbol{m}$ a continuous function from $X$ into the
dual $V^{*}$ of a finite dimensional vector space $V.$ Given a probability
measure $\nu$ on $X$ such that $\boldsymbol{m}_{*}\nu$ is not contained
in a proper affine subspace of $V^{*},$ define the functions $F(\boldsymbol{\beta})$
and $D(\boldsymbol{e})$ on $V$ and $V^{*}$ by
\[
F(\boldsymbol{\beta})\coloneqq\log\int_{X}e^{-\left\langle \boldsymbol{\beta},\boldsymbol{m}\right\rangle }d\nu,\,\,\,\,D(\boldsymbol{e})\coloneqq\inf\left\{ D(\mu|\nu)\colon\mu\in\mathcal{P}(X),\,\,\int\boldsymbol{m}\mu=\boldsymbol{e}\right\} ,
\]
respectively. Then $F(\boldsymbol{\beta})$ is smooth and strictly
convex on $V$ and
\begin{equation}
-F(\boldsymbol{\beta})=\inf_{\boldsymbol{e}\in V^{*}}\left(\left\langle \boldsymbol{\beta},\boldsymbol{e}\right\rangle +D(\boldsymbol{e})\right),\label{eq:minus F beta as inf}
\end{equation}
 i.e. $D(\boldsymbol{e})$ is the Legendre transform of $F(-\boldsymbol{\beta}).$
Moreover, $D(\boldsymbol{e})$ is convex and lsc. As a consequence,
if $0$ is contained in the interior $U$ of the convex hull of the
support of the measure $\boldsymbol{m}_{*}\nu,$ then there exists
$\boldsymbol{\beta}_{0}\in V$ such that 
\begin{equation}
\int_{X}\boldsymbol{m}e^{-\left\langle \boldsymbol{\beta}_{0},\boldsymbol{m}\right\rangle }d\nu=0,\,\,\,F(\boldsymbol{\beta}_{0})=-D(0),\,\,\boldsymbol{\beta}_{0}=-(dD)(0),\label{eq:beta zero}
\end{equation}
where $dD$ denotes the differential of $D,$ which is well-defined
on all of $U.$
\end{lem}

The previous lemma is without doubt essentially well-known (for example,
similar results are used when establishing the equivalence between
the canonical and micro-canonical ensembles in statistical mechanics,
where the roles of $\boldsymbol{\beta}$ and $\boldsymbol{e}$ are
played by a multi-dimensional generalizations of inverse temperature
and energy, respectively \cite{e-h-t}). But, for completeness, a
proof is provided in the appendix. 

\selectlanguage{american}%
In the following \foreignlanguage{english}{we will denote by $\boldsymbol{\beta_{0}}$
the vector in $\R^{r}$ satisfying
\[
\int_{X}e^{-\left\langle \boldsymbol{\beta}_{0},\boldsymbol{m}\right\rangle }\boldsymbol{m}\mu_{\gamma u}=0,
\]
 whose existence (and uniqueness) follows from combining the previous
two lemmas.}
\selectlanguage{english}%
\begin{lem}
The following inequality holds,
\[
-\frac{1}{N}\log\int_{\left\{ \left|\boldsymbol{m}_{N}\right|<\epsilon\right\} }\mu_{\gamma u}^{\otimes N}\leq D(0)-\frac{1}{N}\inf_{\mu\in\mathcal{S}}\log\int_{\left\{ \left|N^{-1}\sum_{i=1}^{N}\boldsymbol{x}_{i}\right|<\epsilon\right\} }\mu^{\otimes N}+|\boldsymbol{\beta}_{0}|\epsilon.
\]
\end{lem}

\begin{proof}
Set $\nu\coloneqq\mu_{\gamma u}$ and $\nu_{0}\coloneqq\nu e^{-\left\langle \boldsymbol{\beta}_{0},\boldsymbol{m}\right\rangle }/\int_{X}\nu e^{-\left\langle \boldsymbol{\beta}_{0},\boldsymbol{m}\right\rangle }.$
Then $\int_{X}\boldsymbol{m}\nu_{0}=0$ and 
\[
\int_{\left\{ \left|\boldsymbol{m}_{N}\right|<\epsilon\right\} }\nu^{\otimes N}=\left(\int_{X}\nu e^{-\left\langle \boldsymbol{\beta}_{0},\boldsymbol{m}\right\rangle }\right)^{N}\int_{\left\{ \left|\boldsymbol{m}_{N}\right|<\epsilon\right\} }e^{-N\left\langle \boldsymbol{\beta}_{0},\boldsymbol{m}_{N}\right\rangle }\nu_{0}^{\otimes N}
\]
 Hence, 
\[
\int_{\left\{ \left|\boldsymbol{m}_{N}\right|<\epsilon\right\} }\nu^{\otimes N}\geq\left(\int_{X}\nu e^{-\left\langle \boldsymbol{\beta}_{0},\boldsymbol{m}\right\rangle }\right)^{N}e^{-N|\boldsymbol{\beta}_{0}|\epsilon}\int_{\left\{ \left|\boldsymbol{m}_{N}\right|<\epsilon\right\} }\nu_{0}^{\otimes N},
\]
giving 
\[
-\frac{1}{N}\log\int_{\left\{ \left|\boldsymbol{m}_{N}\right|<\epsilon\right\} }\nu^{\otimes N}\leq-\log\left(\int_{X}\nu e^{-\left\langle \boldsymbol{\beta}_{0},\boldsymbol{m}\right\rangle }\right)+|\boldsymbol{\beta}_{0}|\epsilon-\frac{1}{N}\log\int_{\left\{ \left|\boldsymbol{m}_{N}\right|<\epsilon\right\} }\nu_{0}^{\otimes N}.
\]
By the previous lemma, $-\log\left(\int_{X}\nu e^{-\left\langle \boldsymbol{\beta}_{0},\boldsymbol{m}\right\rangle }\right)=D(0).$
Finally, since the pushforward $\mu$ of $\nu_{0}$ under $\boldsymbol{m}$
is in $\mathcal{S}$ this concludes the proof of the lemma. 
\end{proof}
All that remains to establish Prop \ref{prop:upper bound mean entr}
is thus to prove the following
\begin{lem}
The following inequality holds, where $C_{\boldsymbol{m}}$ is a constant
only depending on $\boldsymbol{m}$: 
\[
\frac{R_{\mathrm{in}}}{2}|\boldsymbol{\beta}_{0}|\leq-\gamma\left(\mathcal{D}_{\gamma}(u)+\delta_{\gamma}(u)+\sup_{X}u-\mathcal{E}(u)+\frac{C_{\boldsymbol{m}}}{\gamma}\right).
\]
\end{lem}

\begin{proof}
Set $\nu\coloneqq\boldsymbol{m}_{*}\mu_{\gamma u}.$ By Lemma \ref{lem:Legendre duality for m},
$\boldsymbol{\beta}_{0}=dD(0)$ and hence $|\boldsymbol{\beta}_{0}|$
is the norm of the gradient $\nabla D$ at $\boldsymbol{e}=0.$ Now
fix a vector $\boldsymbol{e}$ in the interior $U$ of the convex
hull of $\boldsymbol{m}(X)$ that is a positive multiple of $\nabla D(0).$
Since $D(\boldsymbol{e})$ is convex on the segment between $0$ and
$\boldsymbol{e}$,
\[
\left|\nabla D(0)\right|\left|\boldsymbol{e}\right|\leq D(\boldsymbol{e})-D(0)\leq D(\mu|\mu_{\gamma u})-D(0),
\]
 where $\mu$ is any probability measure on $X$ such that $\int\boldsymbol{m}\mu=\boldsymbol{e}$
(appearing in the definition of $D(\boldsymbol{e})).$ We may additionally
choose $\mu$ to be absolutely continuous with continuous density
by choosing $\mu=e^{-\langle\boldsymbol{\beta}'_{0},\boldsymbol{m}\rangle}\mu_{0}$
for suitable $\boldsymbol{\beta}'_{0}$, where $\mu_{0}$ denotes
the volume form corresponding to the fixed $\K$-invariant metric
$\phi_{0}$ on $-K_{X},$ which is assumed to be a probability measure.
The existence of such a $\boldsymbol{\beta_{0}^{'}}$ follows from
Lemma \ref{lem:Legendre duality for m} (with $\nu=\mu_{0}$). Hence,
we may assume that 
\[
\mu\leq e^{C_{\boldsymbol{m}}}\mu_{0},
\]
 where $\mu_{0}$ denotes the volume form corresponding to the fixed
$\K$-invariant metric $\phi_{0}$ on $-K_{X},$ which is assumed
to be a probability measure. Rewriting $\mu/\mu_{\gamma u}=(\mu/\mu_{0})\mu_{0}/\mu_{\gamma u}$
thus gives 
\[
D(\mu|\mu_{\gamma u})\leq C_{\boldsymbol{m}}+\int_{X}\mu\log\frac{\mu_{0}}{\mu_{\gamma u}}\leq C_{\boldsymbol{m}}+\gamma\left(\gamma^{-1}\log\int_{X}e^{-\gamma u}\mu_{0}+\sup_{X}u\right).
\]
The first term in the latter bracket coincides, by definition, with
$-\mathcal{D}_{\gamma}(u)-\mathcal{E}(u).$ Moreover, by Lemma \ref{lem:Legendre duality for m},
$D(0)=\gamma\delta_{\gamma}(u).$ Hence, combining the previous inequalities
and using that, by construction, $\boldsymbol{e}$ can be taken to
be any vector in the open ball of radius $R_{\text{out}}$ centered
at the origin (and in particular, so that $|\boldsymbol{e}|=R_{\text{out}}/2$),this
concludes the proof.
\end{proof}

\subsection{Proof of Cor \ref{cor:lower bound on delta A ps text}}
\begin{proof}
Take $\gamma\in]0,\gamma_{\epsilon_{N}}^{(N)}(X)[$ and $\phi_{0}\in\mathcal{H}(X,\omega_{0})^{\K}.$
Then Theorem \ref{thm:effective LD bound with F} implies that $F_{-\tilde{\gamma}}$
is bounded from below on $\mathcal{P}(X)_{0},$ for $\tilde{\gamma}\coloneqq\gamma_{k}-\left(n\sigma_{k}(\phi_{0})+n^{-1}\epsilon_{N}\right)\left(1-\frac{2\epsilon_{N}}{R_{\text{in}}}\right)^{-1}.$
Hence, $\delta_{\text{PS}}^{\mathrm{A}}(X,\Delta)\geq\tilde{\gamma},$
giving 
\[
\delta_{\text{PS}}^{\mathrm{A}}(X,\Delta)\geq\frac{1+k^{-1}}{1/\gamma_{\epsilon_{N}}^{(N)}(X)+k^{-1}}-C_{2}k^{-1}
\]
when $k\geq k_{2}$ for effective constants $k_{2}$ and $C_{2},$
using the effective bounds in Prop \ref{prop:bound on Bergman} and
formula \ref{eq:effective RR bound}.
\end{proof}

\subsection{Proof of Cor \ref{cor:asymptotic upper bound on minus log Z text}}
\begin{proof}
To prove the first inequality in Cor \ref{cor:asymptotic upper bound on minus log Z text}
we may as well assume that $\inf_{\mathcal{P}(X)_{0}}F_{-\gamma}$
is finite. Given $\delta>0$ we can thus take $\mu\in\mathcal{P}(X_{0})$
such that $F_{-\gamma}(\mu)\leq\inf_{\mathcal{P}(X)_{0}}F_{-\gamma}+\delta.$
Since $\sigma_{k}(\phi_{0})\leq O(k^{-1})$ Theorem \ref{thm:effective LD bound with F}
yields, using that and $\widetilde{\gamma}=\gamma+O(\epsilon_{N})+O(k^{-1}),$
\[
-\frac{\log Z_{N,\epsilon_{N}}(-\gamma)}{N\gamma(1+k^{-1})}\leq\left(1-O(\epsilon_{N})\right)\left(F_{-\gamma}(\mu)+\left(O(\epsilon_{N})+O(k^{-1})\right)E(\mu)\right)+O(k^{-1}\log N).
\]
 Letting first $N\rightarrow\infty$ and then $\delta\rightarrow0$
thus concludes the proof. Moreover, the second inequality in Cor \ref{cor:asymptotic upper bound on minus log Z text}
also follows by taking $\delta=0.$ Finally, taking $\gamma=1$ and
assuming that $X$ admits a KE metric, the second inequality in Cor
\ref{cor:asymptotic upper bound on minus log Z text} specializes
the third one, since $\inf_{\mathcal{P}(X)_{0}}F_{-1}$ then coincides
with $\inf_{\mathcal{H}(X,\omega_{0})}\mathcal{M}.$
\end{proof}

\section{Appendix}

\subsection{Proof of Lemma \ref{lem:Legendre duality for m} }
\begin{proof}
The smoothness (and even real-analyticity) of $F(\boldsymbol{\beta})$
follows directly from the dominated convergence theorem, since $X$
is compact and $\boldsymbol{m}$ is thus bounded --- the strict convexity
of $F(\boldsymbol{\beta})$ follows from Hölder's inequality. Furthermore,
\[
-F(\boldsymbol{\beta})=\inf_{\mu\in\mathcal{P}(X)}F_{\boldsymbol{\beta}}(\mu),\,\,F_{\boldsymbol{\beta}}(\mu)\coloneqq\left\langle \boldsymbol{\beta},\int_{X}\boldsymbol{m}\mu\right\rangle +D(\mu|\nu),
\]
 using that $D(\mu|e^{-\left\langle \boldsymbol{\beta},\boldsymbol{m}\right\rangle }\nu/\int_{X}e^{-\left\langle \boldsymbol{\beta},\boldsymbol{m}\right\rangle })\geq0.$
Hence, by decomposing $\mathcal{P}(X)$ as the disjoint union of the
subspaces $\mathcal{P}_{\boldsymbol{e}}(X)$ where $\int_{X}\boldsymbol{m}\mu=\boldsymbol{e},$
we get 
\[
-F(\boldsymbol{\beta})=\inf_{\boldsymbol{e}\in V^{*}}\left(\left\langle \boldsymbol{\beta},\boldsymbol{e}\right\rangle +\inf_{\mathcal{P}_{\boldsymbol{e}}(X)}D(\mu|\nu)\right),
\]
 proving formula \ref{eq:minus F beta as inf}. Next, the convexity
and lower semi-continuity of $D(\boldsymbol{e})$ on $V$ follows
directly from the convexity and lower semi-continuity of $\mu\mapsto D(\mu|\nu)$
on $\mathcal{P}(X).$ Note that, by Legendre duality, formula \ref{eq:minus F beta as inf}
equivalently means that $\phi(\boldsymbol{\beta})\coloneqq F(-\boldsymbol{\beta})$
is the Legendre transform of $D(\boldsymbol{e})$: 
\[
\phi=D^{*}.
\]
Next, setting $\eta\coloneqq\boldsymbol{m}_{*}\nu$ we can express
\begin{equation}
F(-\boldsymbol{\beta})=\phi(\boldsymbol{\beta})=\log\int_{V^{*}}e^{\left\langle \boldsymbol{\beta},\boldsymbol{e}\right\rangle }d\eta(\boldsymbol{e}).\label{eq:phi  as Laplace transform}
\end{equation}
In general, by basic convex analysis, a smooth strictly convex function
$\phi$ on a finite dimensional vector space $V$ is proper iff $0$
is contained in the interior of $\overline{(d\phi)(V)}.$ Moreover,
when $\phi$ is of the form \ref{eq:phi  as Laplace transform}, $\overline{(d\phi)(V)}$
coincides with the convex hull $C$ of the support of the measure
$\eta$ and $(d\phi)(V)$ coincides with the interior $\mathring{C}$
of $C.$ Hence, $\phi$ admits a critical point $\boldsymbol{\beta}_{0},$
i.e. $d\phi(\boldsymbol{\beta}_{0})=0,$ which equivalently means
that $\boldsymbol{\beta}_{0}$ satisfies the first equation in formula
\ref{eq:beta zero}. Moreover, since $\phi$ is smooth and strictly
convex, $d\phi$ defines a diffeomorphism between $V$ and $\mathring{C}.$
Since $\phi^{*}=D,$ it follows that the differential $dD$ of $D$
is well-defined on $\mathring{C}$ and coincides with the inverse
of $\boldsymbol{\beta}\mapsto d\phi(\boldsymbol{\beta}).$ Moreover,
on $\mathring{C}$
\[
\phi(dD(\boldsymbol{e}))=\left\langle \boldsymbol{e},dD(\boldsymbol{e})\right\rangle -D(\boldsymbol{e}).
\]
 Setting $\boldsymbol{e}=0$ thus concludes the proof of formula \ref{eq:beta zero}.
\end{proof}

\end{document}